\documentclass[12pt]{amsart}
\usepackage[]{amscd,amsfonts,amsmath,amsxtra,amssymb}
\usepackage{hyperref}
\usepackage{graphics}
\usepackage[all]{xy}
\usepackage[left]{lineno}
\newtheorem{sub}{}[section]
\newtheorem{subsub}{}[sub]

\def\hfl#1#2{\smash{\mathop{\ \hbox to 12mm{\rightarrowfill}}
\limits^{\scriptstyle#1}_{\scriptstyle#2} \ }}
\def\hflb#1#2{\smash{\mathop{\hbox to 12mm{\leftarrowfill}}
\limits^{\scriptstyle#1}_{\scriptstyle#2}}}

\def\og{\leavevmode\raise.3ex\hbox{$\scriptscriptstyle\langle\!\langle$}}
\def\fg{\leavevmode\raise.3ex\hbox{$\scriptscriptstyle\,\rangle\!\rangle$}}

\def\Ssect#1#2{\pagebreak[3]\begin{sub}\label{#2}{\sc\small\small
#1}\rm\medskip}

\def\xmat#1{\[\xymatrix{#1}\]}
\def\flinc{\ar@{^{(}->}}
\def\fleq{\ar@{=}}
\def\flon{\ar@{->>}}
\def\fmaps{\ar@{|-{>}}}

\def\AA{\mathop{\mathcal A}\nolimits}
\def\BB{\mathop{\mathcal B}\nolimits}

\def\LL{\mathop{\mathcal L}\nolimits}

\def\SS{\mathop{\mathcal S}\nolimits}

\def\WW{\mathop{\mathcal W}\nolimits} 
\def\XX{\mathop{\mathcal X}\nolimits}

\def\ZZ{\mathop{\mathcal Z}\nolimits}

\newcommand{\N}{{\mathbb N}}
\newcommand{\M}{{\mathbb M}}
\newcommand{\A}{{\mathbb A}}
\newcommand{\B}{{\mathbb B}}
\newcommand{\Z}{{\mathbb Z}}

\newcommand{\C}{{\mathbb C}}
\newcommand{\Q}{{\mathbb Q}}
\renewcommand{\P}{{\mathbb P}}

\newcommand{\F}{{\mathbb F}}
\newcommand{\E}{{\mathbb E}}

\newcommand{\K}{{\mathbb K}}

\renewcommand{\L}{{\mathbb L}}

\newcommand{\ka}{{\mathcal A}}
\newcommand{\kb}{{\mathcal B}}

\newcommand{\ke}{{\mathcal E}}
\newcommand{\kf}{{\mathcal F}}
\newcommand{\kg}{{\mathcal G}}

\newcommand{\kk}{{\mathcal K}}
\newcommand{\kl}{{\mathcal L}}
\newcommand{\km}{{\mathcal M}}
\newcommand{\kn}{{\mathcal N}}
\newcommand{\ko}{{\mathcal O}}

\newcommand{\kq}{{\mathcal Q}}

\newcommand{\ku}{{\mathcal U}}
\newcommand{\kv}{{\mathcal V}}

\newcommand\bigzero{\makebox(0,0){\text{\huge0}}}

\begin{document}

\def\refname{R\'ef\'erences}
\def\contentsname{ }
\def\proofname{D\'emonstration}
\def\abstractname{R\'esum\'e}

\author{Mohamed Bahtiti}
\address{Institut de Math\'{e}matiques de Jussieu,
 Case 247, 4 place Jussieu,\newline
F-75252 Paris, France}
\email{mohamed.bahtiti@imj-prg.fr }
\title{ Fibr\'e vectoriel de rang $2n+1$ sur~l'espace $\P^{2n+2}$}
 
\maketitle

\begin{abstract}
Nous construisons dans cet article des nouvelles familles de fibr\'es vectoriels alg\'ebriques de rang $2n+1$ sur l'espace projectif complexe $\P^{2n+2}$ \`a partir de deux fibr\'es \mbox{de rang $2n$} sur $\P^{2n+1}$, le fibr\'e de corr\'elation nulle pond\'er\'e et le fibr\'e de Tango pond\'er\'e tout en utilisant la m\'ethode de Kumar-Peterson-Rao. Nous construisons un autre exemple de fibr\'e vectoriel de rang $3$ sur $\P_{\K}^{4}$ diff\'erent de ce qui est pr\'esent\'e dans l'article de Kumar-Peterson-Rao, o\`u $\K$ est un corps quelconque.
$$$$

{\scriptsize ABSTRACT}. We build in this paper new families of algebraic vector bundles of rank $2n+1 $ on the complex projective space $ \P ^ {2n +2} $ from two bundles of rank $2n$ on $\P^{2n+1}$, the weighted null correlation bundles and the weighted Tango bundles while using the method of Kumar-Peterson-Rao. We construct another example of vector bundle of rank $ 3 $ on $ \P_{\K}^{4} $ different than the one presented in the paper of Kumar-Peterson-Rao, where $\K$ is any field.
\end{abstract}

{\scriptsize {\em Date}: March, 2016.}
{\scriptsize Mots-cl\'es. Fibr\'e de Tango pond\'er\'e, Fibr\'e vectoriel de $0$-corr\'elation pond\'er\'e.}

\vspace{1.8cm}

\section{\large\bf Introduction} 
\vspace{1cm}

Soient $n, \alpha ,\gamma \in \N$ et $\beta \in \Z$ tels que $n>0$, $\gamma >0$, $\alpha \geq \beta$, $\alpha +\beta\geq 0$ et $\gamma > (2n+1)\mid \beta \mid$. Soient $U$ un $\C$-espace vectoriel de dimension $2$, et $Q$ le fibr\'e de quotient et $N(f)$ le fibr\'e de corr\'elation nulle classique sur $\P(\SS^{2n+1}U)$ pour $f\in (\bigwedge^{2}\SS^{2n+1}U)^{*}=H^{0}(Q^{*}(1))$ une $(2n+2)\times (2n+2)$-matrice antisym\'etrique symplectique non d\'eg\'en\'er\'ee. Ces fibr\'es sont d\'efinis par les suites exactes suivantes
$$0\longrightarrow \ko_{\P(\SS^{2n+1}U)}(-1)\stackrel{g}{\longrightarrow} \SS^{2n+1}U\otimes \ko_{\P(\SS^{2n+1}U)} \longrightarrow Q \longrightarrow 0,$$

et
$$0\longrightarrow N(f) \longrightarrow  Q  \stackrel{^{T}g\circ f}{\longrightarrow} \ko_{\P(\SS^{2n+1}U)}(1) \longrightarrow 0.$$

Soit $F(W)$ le fibr\'e de Tango pour $W \subset \bigwedge^{2}\SS^{2n+1}U=(H^{0}(Q^{*}(1)))^{*} $ un sous-espace vectoriel tel que
$ dim(W)= \left(\begin{array}{c}2n+2\\2\end{array}\right)- 4n-1 $ et que $W$ ne contient pas d'\'el\'ement d\'ecomposable non nul de $  \bigwedge^{2}\SS^{2n+1}U$. Ce fibr\'e est d\'efini par la suite exacte suivante

$$0\longrightarrow Q(-1)\stackrel{\varpi_{W}}{\longrightarrow} \left(    (\bigwedge^{2}\SS^{2n+1}U)\diagup W \right)\otimes\ko_{\P(\SS^{2n+1}U)} \longrightarrow  F(W)(1)   \longrightarrow 0.$$ 

Tous les sous-espaces comme $W$, qui sont $\C^{*}$-invariants par rapport \`a la $\C^{*}$-action suivante
$$\sigma_{1}=\left(\begin{array}{cc } t^{\alpha}&0 \\0&    t^{\beta} \\ \end{array} \right) \text{ avec } t\in \C^{*}, $$

forment l'ensemble $\WW$ qui \'etait \'etudi\'e dans le th\'eor\`eme $3.4$ \cite{bah1}. Dans le th\'eor\`eme \ref{3.2.3}, on trouve $W \in \WW$ et $f:=f_{W}\in (\bigwedge^{2}\SS^{2n+1}U)^{*}=H^{0}(Q^{*}(1))$ une $(2n+2)\times (2n+2)$-matrice antisym\'etrique symplectique non d\'eg\'en\'er\'ee tels que 
$$^{T}\sigma_{1} f\sigma_{1}= af \text{ et } \langle  f,w \rangle=0, \hspace{0,2 cm} a\in \C^{*}, \text{ pour tout } w\in W $$

o\`u $ \langle , \rangle $ est la forme quadriatique associ\'ee \`a $\bigwedge^{2}\SS^{2n+1}U$ dans \ref{3.2.2}. On trouve aussi $D_{W}$ un sous-espace vectoriel $\C^{*}$-invariant de $\bigwedge^{2}\SS^{2n+1}U$ tel que $\bigwedge^{2}\SS^{2n+1}U\simeq D_{W}\oplus W $. D'apr\`es la proposition \ref{3.2.1}, on obtient la suite exacte suivante 
$$0\longrightarrow N(f)(-1)\longrightarrow  \left( D_{W} \diagup  \C.f  \right)  \otimes\ko_{\P^{2n+1}}  \longrightarrow  F(W)(1)   \longrightarrow 0.$$

On consid\`ere la transform\'ee de Horrocks \cite{bah1}
\xmat{{\bf Iminvg}:  \kf\kv(\P(\SS^{2n+1} U),\overline{ \sigma})   \ar[rr] &  &  \kf\kv(\P^{2n+1}),  \\
 }
 
o\`u $\kf\kv(\P(\SS^{2n+1} U),\overline{\sigma})$ est la cat\'egorie de fibr\'es vectoriels sur $\P(\SS^{2n+1} U)$ qui sont $\C^{*}$-invariants au-dessus de l'action $\overline{ \sigma}= t^{\gamma}.\sigma_{1}$, et $ \kf\kv(\P^{2n+1})$ est la cat\'egorie de fibr\'es vectoriels sur $\P^{2n+1}$ qui sont $\C^{*}$-invariants au-dessus de la multiplication usuelle sur $\C^{2n+2}$. Dans la proposition \ref{3.2.5}, on applique cette transform\'ee sur la suite pr\'ec\'edente pour obtenir la suite exacte suivante
$$0\longrightarrow \kn(-\gamma)\stackrel{\pi}{\longrightarrow} \Upsilon_{1}  \stackrel{\iota}{\longrightarrow}  \kf(\gamma) \longrightarrow 0,$$

o\`u $\kn$ (resp. $\kf$) est l'image inverse g\'en\'eralis\'ee du fibr\'e $N(f)$ dans la proposition \ref{1.4.8} \mbox{(resp. $F(W)$} dans la proposition 3.7 \cite{bah1}), et
$$  \Upsilon_{1}=\bigoplus_{i=1}^{4n} \ko_{\P^{2n+1}}(\zeta_{i}),$$

avec 
$$\zeta_{1}=\alpha(4n+1)+\beta,\zeta_{2}=4n\alpha+2\beta,\ldots,\zeta_{2n}=\alpha(2n+2)+2n\beta,\zeta_{2n+1}=2n\alpha+(2n+2)\beta,$$
$$\ldots,\zeta_{4n}=\alpha+(4n+1)\beta .$$

Dans la proposition \ref{3.2.61}, pour la suite pr\'ec\'edente il existe deux morphismes $\phi$ et $\psi$ tels que le diagramme suivant est commutatif et exacte en $\Upsilon_{1}$

\xmat{ \Upsilon_{1}^{*} (-\hslash_{1})  \ar[rr]^{\phi}\ar[rd]_{}  &  & \Upsilon_{1} 
\ar[rr]^{\psi}\ar[rd]^{\iota}  & & \bigwedge ^{2n-1}\Upsilon_{1}^{*}(\hslash_{2})\\
&   \kn(-\gamma) \ar[ru]^{\pi}  \ar[rd]& &  \kf(\gamma)   \ar[rd]\ar[ru]^{}    \\
 0\ar[ru] &&   0\ar[ru]    && 0 ,}
  
o\`u $\hslash_{1}=2\gamma- (2n+1)(\alpha+\beta )$ et $\hslash_{2}=2n\gamma+3n(2n+1)(\alpha+\beta )$. De plus, si 
$$\gamma> max\lbrace \mid \alpha+\beta (4n+1)\mid + (2n+1)(\alpha+\beta ), \alpha(4n+1)+\beta\rbrace  ,$$

alors les deux morphismes $\phi$ et $\psi$ sont des matrices de rang $2n$ et $\hslash_{1}>0$.

Soient $z \in H^{0}(\P^{2n+2},\ko_{\P^{2n+2}}(1))\setminus \lbrace 0 \rbrace $ et $\P^{2n+1}   \subset  \P^{2n+2} $ l'hyperplan de $\P^{2n+2}$ d\'efini par l'\'equation $z=0$. Soient $\epsilon\geq 1$ un entier et $Y_{(\epsilon)}\subset \P^{2n+2}$ le voisinage infinit\'esimal de l'ordre $\epsilon$ de $\P^{2n+1} $ dans $\P^{2n+2}$, d\'efini par l'\'equation $z^{\epsilon}=0$. Soit $J_{(\epsilon)}:Y_{(\epsilon)}\longrightarrow \P^{2n+1}$ la projection d'un point, en dehors de $Y_{(\epsilon)}$, sur un sous-espace projectif $\P^{2n+1}$ (\cite{ha} page $22$). En appliquant l'image inverse de $ J_{(\epsilon)}$ sur la suite exacte pr\'ec\'edente, on obtient la suite exacte suivante 
$$0\longrightarrow \kn_{(\epsilon)}(-\gamma)\stackrel{\pi}{\longrightarrow}  \Upsilon_{(\epsilon)}  \stackrel{\iota}{\longrightarrow}  \kf_{(\epsilon)}(\gamma)   \longrightarrow 0,$$

o\`u $\kn_{(\epsilon)}=J_{(\epsilon)}^{*}\kn$ et $\kf_{(\epsilon)}=J_{(\epsilon)}^{*}\kf$, et 
$$ \Upsilon_{(\epsilon)}:=J_{(\epsilon)}^{*}\Upsilon_{1}= \bigoplus_{i=1}^{4n} \ko_{Y_{(\epsilon)} }(\zeta_{i}).$$

Nous utilisons la m\'ethode de Kumar-Peterson-Rao dans l'article \cite{ku-ra-pe} (voir \ref{3.1}) pour construire un fibr\'e vectoriel sur $\P^{2n+2}$ \`a partir de la suite exacte pr\'ec\'edente qui est d\'efinie sur $Y_{(\epsilon)}$.\\
Soit $\kg$ un fibr\'e vectoriel sur $\P^{2n+2}$ d\'efini par le diagramme commutatif suivant
\xmat{ &0\ar[d]&0\ar[d]&&\\
&\Upsilon(-\epsilon)\ar@{=}[r] \ar[d]& \Upsilon(-\epsilon)\ar[d]&&\\
0 \ar[r]&\kg \ar[r]\ar[d]&\Upsilon \ar[r]\ar[d]& \kf_{(\epsilon)}(\gamma) \ar[r]\ar@{=}[d]&0\\
0 \ar[r]& \kn_{(\epsilon)}(-\gamma) \ar[r]^{\pi}\ar[d]&\Upsilon_{(\epsilon)} \ar[r]^{\iota}\ar[d]& \kf_{(\epsilon)}(\gamma) \ar[r]&0\\
&0&0&&}

o\`u 
$$ \Upsilon:= \bigoplus_{j=1}^{4n} \ko_{\P^{2n+2}}(\zeta_{j}).$$ 

Dans le th\'eor\`eme \ref{3.2.7}, pour $n>1$, soit 
$$\Gamma:=\bigoplus_{i=1}^{2n-1} \ko_{\P^{2n+2}}(\zeta_{b_{i}})  \text{ avec } 1\leq b_{1}<b_{2}< \ldots<b_{2n-1}\leq 4n.$$
On a:

- Si $\epsilon=\epsilon_{1}:=2\gamma+(2n+1)(\alpha+\beta )$, alors il existe un fibr\'e vectoriel $\ke$ de rang $2n+1$ sur $\P^{2n+2}$ qui est d\'efini par la suite exacte suivante
$$0 \longrightarrow \Gamma(-\epsilon_{1}) \longrightarrow \kg \longrightarrow \ke \longrightarrow 0. $$
- Si $\epsilon=\epsilon_{2}:=2n\gamma+n(2n+1)(\alpha+\beta )$, alors il existe un fibr\'e vectoriel $\kk$ de rang $2n+1$ sur $\P^{2n+2}$ qui est d\'efini par la suite exacte suivante
$$0 \longrightarrow \kk \longrightarrow \kg \longrightarrow \Gamma\longrightarrow 0. $$
Dans la proposition \ref{3.2.8}, on trouve des $3$-fibr\'es sur $\P^{4}$ pour certaines valeurs de $\epsilon$.

Dans la proposition \ref{4.2.4}, on trouve un autre exemple de $3$-fibr\'es sur $\P_{\K}^{4}$ diff\'erent de l'exemple 4.1 de Kumar-Peterson-Rao \cite{ku-ra-pe} tout en utilisant l'identit\'e de Binet-Cauchy, o\`u $\K$ est un corps quelconque. Parmi ces fibr\'es Il y a des fibr\'es qui sont diff\'erents du fibr\'e de Tango pond\'er\'e provenant d'une image inverse g\'en\'eralis\'ee sur $\P^{2n+2}$ \cite{bah1}.

\section{ \large\bf Pr\'eliminaires} 

\vspace{1cm}

\subsection{Fibr\'e de 0-corr\'elation pond\'er\'e sur $\P^{2n+1}$ provenant d'une image inverse g\'en\'eralis\'ee} \label{1.4} 
Nous allons \'etudier dans cette partie un fibr\'e de 0-corr\'elation 
pond\'er\'e sur $\P^{2n+1}$ qui provient d'une image inverse g\'en\'eralis\'ee de fibr\'e de 
corr\'elation nulle classique sp\'ecial sur $\P(\SS^{2n+1}U)$.

\subsubsection{ D\'efinition du fibr\'e de corr\'elation nulle classique sur 
$\P^{2n+1}$}\label{1.4.5} 
Soient $V$ un espace \mbox{vectoriel} complexe de dimension $dim(V)=2n+2$,  
$\P^{2n+1}=\P(V)$ l'espace projectif \mbox{complexe} associ\'e dont les points sont les droites de $V$, et $v\in V$ un point correspondant \`a la droite $x=\C.v \in P(V)$. On a l'inclusion canonique
$$g_{x}:(\ko_{\P^{2n+1}}(-1))_{x}:=\C.v\hookrightarrow V$$
$$\hspace{5 cm}av\longmapsto g_{x}(av)=av ,$$

o\`u $a\in \C$. Alors on a la suite exacte suivante
$$ 0\longrightarrow \ko_{\P^{2n+1}}(-1)\stackrel{g}{\longrightarrow} V\otimes \ko_{\P^{2n+1}}\longrightarrow Q\longrightarrow 0.$$

On a aussi le morphisme surjectif
$$\hspace{1 cm}{}^{T}g_{x}:V^{*}\longrightarrow \C.v^{*}:=\ko_{\P^{2n+1},x}(1)$$
$$h\longmapsto {}^{T}g_{x}(h),$$ 

o\`u $^{T}g_{x}(h)$ est la restriction de h \`a $\C.v$. Ce morphisme nous donne le morphisme surjectif suivant
$$V^{*}\otimes\ko_{\P^{2n+1}} \stackrel{^{T}g}{\longrightarrow} \ko_{\P^{2n+1}}(1)\longrightarrow 0.$$

Soit $f:V\longrightarrow V^{*}$ un morphisme. On sait que f est antisym\'etrique si et seulement si, pour tout 
$u \in V$, $f(u)(u)=0$. Alors $(^{T}g_{x}\circ f \circ g_{x})(v) = {}^{T}g_{x}(f(v))$, la restriction de $f(v)$ \`a $\C.v$, est nul si et seulement si $f(v)(v)=0$. 

Soit $f:V\longrightarrow V^{*}$ un isomorphisme antisym\'etrique. On peut d\'efinir la monade suivante
$$0\longrightarrow \ko_{\P^{2n+1}}(-1)\stackrel{g}{\longrightarrow} V\otimes \ko_{\P^{2n+1}} \stackrel{f}{\longrightarrow} V^{*}\otimes \ko_{\P^{2n+1}} \stackrel{^{T}g}{\longrightarrow}\ko_{\P^{2n+1}}(1) \longrightarrow 0 $$ 

qui peut s'\'ecrire \'egalement
$$0\longrightarrow \ko_{\P^{2n+1}}(-1)\stackrel{g}{\longrightarrow} V\otimes \ko_{\P^{2n+1}} \stackrel{^{T}g\circ f}{\longrightarrow}\ko_{\P^{2n+1}}(1) \longrightarrow 0 .$$

On d\'efinit {\em le fibr\'e de corr\'elation nulle classique} $N(f)$ sur $\P^{2n+1}$ par la cohomologie de la monade pr\'ec\'edente. 

Maintenant, on consid\`ere $f:V\longrightarrow V^{*}$ un isomorphisme antisym\'etrique d\'efini relativement \`a la base canonique de $V$ par la $(2n+2)\times (2n+2)$-matrice antisym\'etrique symplectique suivante

$$
\left(
\begin{array}{cccccccccc}
  & &  &  &   & & & & & -r_{0} \\
  & &  &\bigzero  &   & & & & .&  \\
  & &  &  &   & & &. & &  \\
  & &  &  &   & &.& & &  \\
  & &  &  &   & -r_{n}& & & &  \\
 & &  &  & r_{n} & & & & &  \\
  & &  &. &   & & & & &  \\
  & &. &  &   & & & & &  \\
  &.&  &  &   & & \bigzero& & &   \\
r_{0} & &  &  &  & & & & &   \\
\end{array}
\right).$$
 
o\`u $r_{0},r_{1}, \ldots  ,r_{n}\in \C^{*}$. Par d\'efinition, {\em le fibr\'e de corr\'elation nulle classique sp\'ecial} $N(f)$ sur $\P^{2n+1}$ est la cohomologie de la monade suivante

\begin{equation}
0\longrightarrow \ko_{\P^{2n+1}}(-1)\stackrel{g}{\longrightarrow} V\otimes \ko_{\P^{2n+1}}  \stackrel{^{T}g\circ f}{\longrightarrow}\ko_{\P^{2n+1}}(1) \longrightarrow 0.
\end{equation}
 
Donc on a les suites exactes suivantes

\begin{equation}
0\longrightarrow \ko_{\P^{2n+1}}(-1)\stackrel{g}{\longrightarrow} V\otimes \ko_{\P^{2n+1}} \longrightarrow Q \longrightarrow 0,
\end{equation}

o\`u $Q$ est le fibr\'e de quotient vectoriel de rang $2n+1$ sur $\P^{2n+1}$,
et

\begin{equation}
0\longrightarrow N(f) \longrightarrow  Q  \stackrel{^{T}g\circ f}{\longrightarrow} \ko_{\P^{2n+1}}(1) \longrightarrow 0.
\end{equation}
   
\subsubsection{\bf Proposition }\label{1.4.6}
{\em
1- Soient $f_{1}, f_{2}: V\longrightarrow V^{*}$ des isomorphismes antisym\'etriques. Alors on a 
$N(f_{1})\simeq N(f_{2})$ si et seulement s'il existe $a \in \C^{*}$ tel que $f_{1}=a.f_{2}$.

2- Soient $\kappa \in GL(V)$, $\overline{\kappa } \in PGL(V) $ son \'el\'ement correspondant \`a $\kappa$ et $f: V\longrightarrow V^{*}$ un \mbox{isomorphisme} antisym\'etrique. Alors on a
$$\overline{\kappa }^{*}(N(f))\simeq N(^{T}\kappa \circ f\circ \kappa).$$

3-  Soient $N(f)$ le fibr\'e de corr\'elation nulle classique sp\'ecial sur $\P^{2n+1}$ et $\overline{ \sigma(t)}\in PGL(V)$ une $\C^{*}$-action sur $\P(V)$, pour tout $t\in \C^{*}$, qui est de la forme 
$$ \overline{ \sigma(t)}:=\left(
\begin{array}{ccccccccc}
 t^{a_{0}} &  &   & & &            &   \\
     & t^{a_{1}} &  &   & &   &\bigzero \\
  &  & \ddots &   &    &   &   \\
\bigzero &  &   &    &  t^{a_{2n}}&  &   \\
 &  &   &    &    &  &      t^{a_{2n+1}}\\
\end{array}
\right)
$$

o\`u  $a_{0}\geq a_{1}\geq  \ldots  \geq a_{2n}\geq a_{2n+1}$ sont des entiers tels que $a_{i}+ a_{2n+1-i}=d \in \Z$, $ i=0,1, \ldots , 2n+1 $. Alors on a un isomorphisme canonique
$$s_{t}:\overline{ \sigma(t)}^{*}N(f)\simeq N(f),$$

tel que pour $t_{1},t_{2}\in \C^{*}$, on a $ s_{t_{1}.t_{2}}=s_{t_{2}}\circ \overline{ \sigma(t_{2})}^{*}(s_{t_{1}})$.}

\begin{proof}
1- Soit la suite exacte suivante

$$0\longrightarrow \ko_{\P^{2n+1}}(-1)\stackrel{g}{\longrightarrow} V\otimes \ko_{\P^{2n+1}}  
 \longrightarrow Q \longrightarrow 0. $$

Pour $f_{1}, f_{2}: V\longrightarrow V^{*}$ des isomorphismes antisym\'etriques, on a les suites exactes suivantes
$$0\longrightarrow \ko_{\P^{2n+1}}(-1)\stackrel{^{T}b_{1} }{\longrightarrow} Q^{*} \stackrel{d_{1} }{\longrightarrow} N^{*}(f_{1}) \longrightarrow 0 $$ 

$$0\longrightarrow \ko_{\P^{2n+1}}(-1)\stackrel{ ^{T}b_{2} }{\longrightarrow} Q^{*} \stackrel{d_{2} }{\longrightarrow} N^{*}(f_{2}) \longrightarrow 0 ,$$ 

o\`u $b_{1}:=(^{T}g \circ f_{1})$ et $b_{2}:=(^{T}g \circ f_{2})$. S'il existe $a \in \C^{*}$ tel que $f_{1}=a.f_{2}$, alors on a $ker(^{T}g \circ f_{2})=ker(^{T}g \circ f_{1})$. Donc on obtient $N(f_{1})\simeq N(f_{2})$.

Si $N(f_{1})\simeq N(f_{2})$, alors on a $\psi:N^{*}(f_{1})\stackrel{\sim}{\longrightarrow} N^{*}(f_{2}) $. On a la suite suivante
$$0\longrightarrow Hom(Q^{*},\ko_{\P^{2n+1}}(-1))\longrightarrow Hom(Q^{*},Q^{*})\longrightarrow Hom(Q^{*},N^{*}(f_{2}))$$
$$\longrightarrow Ext^{1}(Q^{*},\ko_{\P^{2n+1}}(-1))\longrightarrow ,$$

mais $Hom(Q^{*},\ko_{\P^{2n+1}}(-1))=0$ et $Ext^{1}(Q^{*},\ko_{\P^{2n+1}}(-1)=0$. On obtient que 
$$Hom(Q^{*},Q^{*})\stackrel{\sim }{\longrightarrow} Hom(Q^{*},N^{*}(f_{2}))$$

Comme $\psi$ est non trivial, donc $\psi\circ d_{1}\neq 0$. Il existe un morphisme non trivial $\varphi:Q^{*}\longrightarrow Q^{*}$ tel que $\psi\circ d_{1}=d_{2}\circ \varphi $. D'apr\`es (\cite{ok-sc-sp}, Lemma 1.2.8), on obtient que $\varphi$ est un isomorphisme. Comme le fibr\'e $Q^{*}$ est simple, on obtient que $\varphi$ est une homoth\'etie. L'isomorphisme $\varphi$ d\'efinit, localement en $x\in \P^{2n+1}$ o\`u $x=\C.v$ pour un $v\in V$, un isomorphisme
$$\varphi_{1,x}:\C.v \stackrel{a \times }{\longrightarrow} \C.v$$
 
tel que $ a. ^{T}b_{2}= ^{T}b_{2}\circ \varphi_{1}=\varphi\circ ^{T}b_{1}=^{T}b_{1}$, pour $a \in \C^{*}$. Donc on a le diagramme commutatif suivant

\xmat{ 0 \ar[r] & \ko_{\P^{2n+1}}(-1) \ar[r]^{^{T } b_{1}}\ar[d]_{\varphi_{1}} & 
Q^{*}\ar[r]^{d_{1}}\ar[d]_{\varphi}  & N^{*}(f_{1}) \ar[r]\ar[d]^{\psi}& 0\\
0 \ar[r] & \ko_{\P^{2n+1}}(-1) \ar[r]^{^{T } b_{2}} & Q^{*} \ar[r]^{d_{2}}  &  N^{*}(f_{2}) \ar[r]& 0\\
}

et on a $(a. f_{2})\circ g =f_{1}\circ g $, pour tout $g:X\hookrightarrow V$ o\`u $X$ est une droite dans $V$, ce qui donne $f_{1}=a.f_{2}$.

2- Pour la monade suivante
$$0\longrightarrow \ko_{\P^{2n+1}}(-1)\stackrel{g}{\longrightarrow} V\otimes \ko_{\P^{2n+1}} \stackrel{^{T}g\circ f}{\longrightarrow}\ko_{\P^{2n+1}}(1) \longrightarrow 0 $$

et d'apr\`es la d\'efinition de l'image inverse, on a $\overline{\kappa }^{*}(g)=\kappa \circ g$ et $\overline{\kappa }^{*}(^{T}g)= ^{T}(\kappa \circ g)$ qui donnent la monade suivante
$$0\longrightarrow \ko_{\P^{2n+1}}(-1)\stackrel{\kappa \circ g}{\longrightarrow} V\otimes \ko_{\P^{2n+1}} \stackrel{^{T}(\kappa \circ g) \circ f}{\longrightarrow}\ko_{\P^{2n+1}} (1)\longrightarrow 0 .$$

Donc on a
$$\overline{\kappa }^{*}(N(f))\simeq N(^{T}\kappa\circ f\circ \kappa).$$

Soit $f:V \longrightarrow V^{*}$ un isomorphisme antisym\'etrique d\'efini relativement \`a la base canonique de ~$V$ par la $(2n+2)\times (2n+2)$-matrice antisym\'etrique symplectique, comme dans la d\'efinition \ref{1.4.5}. Il suffit de remarquer que $^{T} \sigma(t)\circ f\circ \sigma(t)= t^{d} f$, pour tout $t\in C^{*}$. On obtient donc un isomorphisme canonique
$$s_{t}:\overline{ \sigma(t)}^{*}N(f)\simeq N(f),$$

tel que pour $t_{1},t_{2}\in \C^{*}$ on a $ s_{t_{1}.t_{2}}=s_{t_{2}}\circ \overline{ \sigma(t_{2})}^{*}(s_{t_{1}})$, 
\xmat{   \overline{\sigma(t_{2})}^{*} \overline{\sigma(t_{1})}^{*}N(f)=\overline{\sigma( t_{1}.t_{2} )}^{*} N(f) \ar[rr]^{s_{t_{1}. t_{2}}}\ar[rd]_{\overline{\sigma(t_{2})}^{*}(s_{t_{1}) }} &  
&   N(f)  \\
&    \overline{\sigma(t_{2})}^{*} N(f)  \ar[ru]_{s_{t_{2}}}.}

\end{proof} 

\subsubsection{\bf Proposition}\label{1.4.8} 
{\em Soient $U$ un $\C$-espace vectoriel de dimension $2$, et $\lbrace x,y\rbrace$ sa base. Soient $i,n, \alpha ,\gamma \in \N$ et $\beta \in \Z$ tels que $n>1$, $\gamma >0$, $\alpha \geq \beta$, $\alpha +\beta\geq 0$ et $\gamma > (2n+1)\mid \beta \mid$, et $g_{0}, \ldots ,g_{2n+1}$ des formes homog\`enes sans un z\'ero commun sur $\P(\SS^{2n+1}U)$ telles que 
$$ deg(g_{i})=\gamma + (2n+1)\alpha+i(\beta -\alpha),\hspace{0.2 cm} i=0,1, \ldots ,2n+1.$$

Soient
$$H_{0}:=\lbrace x^{2n+1-j}y^{j}  ,  \hspace{0.2 cm} 0\leq j\leq 2n+1 \rbrace $$

la base canonique de $\SS^{2n+1}U$, $Q$ le fibr\'e de quotient et $N(f)$ le fibr\'e de corr\'elation nulle classique sp\'ecial sur $\P(\SS^{2n+1}U)$. Les fibr\'es $Q$ et $N(f)$ sont d\'efinis par les suites exactes
$$ 0\longrightarrow \ko_{\P(\SS^{2n+1}U)}(-1)\stackrel{g}{\longrightarrow} \SS^{2n+1}U\otimes \ko_{\P(\SS^{2n+1}U)} \longrightarrow Q \longrightarrow 0,$$ 

$$0\longrightarrow N(f) \longrightarrow  Q  \stackrel{ ^{T}g\circ f}{\longrightarrow} \ko_{\P(\SS^{2n+1}U)}(1) \longrightarrow 0.$$

o\`u $g$ est l'inclusion canonique, et $f:\SS^{2n+1}U \longrightarrow \SS^{2n+1}U^{*}$ est un isomorphisme antisym\'etrique d\'efini relativement \`a la base canonique pr\'ec\'edente de $\SS^{2n+1}U$ par la $(2n+2)\times (2n+2)$-matrice antisym\'etrique symplectique, comme dans la d\'efinition \ref{1.4.5}. 

Alors le fibr\'e $Q$ (resp. $N(f)$) a une image invers\'ee g\'en\'eralis\'ee $Q_{\gamma,\alpha,\beta}$ (resp. $N_{\gamma,\alpha,\beta}$) d\'efinie par
\begin{equation}\label{12} 
 0\longrightarrow \ko_{\P^{2n+1}}(-\gamma)\stackrel{g(-\gamma)}{\longrightarrow} \SS^{2n+1}\ku\longrightarrow \kq_{\gamma,\alpha,\beta} \longrightarrow 0
\end{equation}

\begin{equation}\label{13} 
 \text{( resp.} \hspace{1 cm}0\longrightarrow \kn_{\gamma,\alpha,\beta} \longrightarrow \kq_{\gamma,\alpha,\beta}  \stackrel{(^{T}g\circ f)(-\gamma)}{\longrightarrow} \ko_{\P^{2n+1}}(\gamma+(\alpha+\beta)(2n+1)) \longrightarrow 0),
\end{equation}

o\`u $\ku =\ko_{\P^{2n+1}}(\alpha )\oplus \ko_{\P^{2n+1}}(\beta)$, $\kn_{\gamma,\alpha,\beta}:=N_{\gamma,\alpha,\beta}(-\gamma)$ et $\kq_{\gamma,\alpha,\beta}:=Q_{\gamma,\alpha,\beta}(-\gamma)$. On appelle le fibr\'e $\kq_{\gamma,\alpha,\beta}$ le fibr\'e de quotient 
pond\'er\'e par les poids $\gamma,\alpha,\beta $, provenant d'une image inverse g\'en\'eralis\'ee sur $\P^{2n+1}$. On appelle le fibr\'e  $\kn_{\gamma,\alpha,\beta} $ le fibr\'e de 0-corr\'elation 
pond\'er\'e par les poids $\gamma,\alpha, \beta$, provenant d'une image inverse g\'en\'eralis\'ee sur $\P^{2n+1}$}.

\begin{proof} On consid\`ere l'application 
$$\omega:= (g_{0},\ldots ,g_{2n+1}): \C^{2n+2}\setminus  \lbrace 0 \rbrace \longrightarrow \SS^{2n+1}U \setminus  \lbrace 0 \rbrace \hspace{4 cm}$$
$$\hspace{2 cm}v\longmapsto (g_{0}(v),\ldots ,g_{2n+1}(v))$$

On consid\`ere l'action de $\C^{*}$ sur $\SS^{2n+1}U $
$$\sigma:\C^{*}\times  \SS^{2n+1} U \longrightarrow \SS^{2n+1} U   $$
$$(t,u)\longmapsto t^{\gamma}.\left(\begin{array}{cc } t^{\alpha}&0 \\0&    t^{\beta} \\ \end{array} \right).u.$$

Cette action est repr\'esent\'ee par la matrice 
$$ t^{\gamma}. \left(
\begin{array}{ccccccccc}
 t^{(2n+1)\alpha} &  &   & & &            &   \\
     & t^{(2n+1)\alpha+(\beta -\alpha)}  &  &   & &   &\bigzero \\
  &  & \ddots &   &    &   &   \\
\bigzero &  &   &    &  t^{(2n+1)\alpha+2n(\beta -\alpha)} &  &   \\
 &  &   &    &    &  &      t^{(2n+1)\beta} \\
\end{array}
\right)\in PGL(\SS^{(2n+1)} U).
$$

On consid\`ere aussi l'action de $\C^{*}$ (multiplication usuelle sur $\C^{2n+2} $)
$$T:\C^{*}\times  \C^{2n+2} \longrightarrow \C^{2n+2}   $$
$$(t,u)\longmapsto t.u.$$

Alors $\omega $ est une $\C^{*}$-application par rapport \`a ces deux actions. L'action $\sigma$ induit une action $\overline{ \sigma}\in PGL(\SS^{2n+1} U) $ de $\C^{*}$ sur $\P(\SS^{2n+1} U)$, et l'action $T$ induit une action triviale de $\C^{*}$ sur $\P^{2n+1}$. On a donc une transform\'ee de Horrocks \cite{bah1}
\xmat{{\bf Iminvg}:  \kf\kv(\P(\SS^{2n+1} U),\overline{ \sigma})   \ar[rr] &  &  \kf\kv(\P^{2n+1}),  \\
 }

o\`u $\kf\kv(\P(\SS^{2n+1} U),\overline{\sigma})$ est la cat\'egorie de fibr\'es vectoriels sur $\P(\SS^{2n+1} U)$ qui sont $\C^{*}$-invariants \mbox{au-dessus} de l'action $\overline{ \sigma}$, et $ \kf\kv(\P^{2n+1})$ est la cat\'egorie de fibr\'es vectoriels sur $\P^{2n+1}$ qui sont $\C^{*}$-invariants au-dessus de l'action $T$. On consid\`ere le morphisme

$$g:={}^{T}\omega:\ko_{\P(\SS^{2n+1}U)}(-1) \longrightarrow \SS^{2n+1}U \otimes \ko_{\P(\SS^{2n+1}U)} ,$$

avec $\ko_{\P(\SS^{2n+1}U)}(-1)$, $\SS^{2n+1}U \otimes \ko_{\P(\SS^{2n+1}U)}$ qui sont munis de l'action canonique $\sigma(t)$. Alors $g$ est un $\C^{*}$-morphisme. Comme on a, pour tout $t\in \C^{*}$,
$$\sigma(t).(x^{2n+1-i}.y^{i}) =t^{\gamma + (2n+1)\alpha+i(\beta -\alpha)}.(x^{2n+1-i}.y^{i}), $$
 
alors le sous-fibr\'e $(x^{2n+1-i}.y^{i}.\C)\otimes \ko_{\P(\SS^{2n+1} U)} $ de $\SS^{2n+1}U \otimes \ko_{\P(\SS^{2n+1}U)}$ est $\C^{*}$-invariant. On obtient
$$\ko_{\P(\SS^{2n+1} U)}^{(\gamma + (2n+1)\alpha+i(\beta -\alpha))}\simeq(x^{2n+1-i}.y^{i}.\C)\otimes \ko_{\P(\SS^{2n+1} U)} , $$

qui est d\'efini localement, pour tout $v\in \SS^{2n+1} U $, par
$$(\ko_{\P(\SS^{2n+1} U)}^{(\gamma + (2n+1)\alpha+i(\beta -\alpha))})_{v} \simeq \C \stackrel{\simeq}{\longrightarrow}((x^{2n+1-i}.y^{i}.\C) \otimes \ko_{\P(\SS^{2n+1} U)})_{v}\simeq x^{2n+1-i}.y^{i}.\C   $$
$$a\longmapsto a x^{2n+1-i}.y^{i}.\hspace{3 cm}$$

Donc on a un $\C^{*}$-isomorphisme
$$\SS^{2n+1}U \otimes \ko_{\P(\SS^{2n+1}U)}\simeq \bigoplus_{i=0}^{2n+1}\ko_{\P(\SS^{2n+1} U)}^{(\gamma + (2n+1)\alpha+i(\beta -\alpha))}. $$

Alors on a un $\C^{*}$-morphisme
$$g:={}^{T}\omega:\ko_{\P(\SS^{2n+1}U)}(-1) \longrightarrow \bigoplus_{i=0}^{2n+1}\ko_{\P(\SS^{2n+1} U)}^{(\gamma + (2n+1)\alpha+i(\beta -\alpha))} .$$

On consid\`ere le morphisme 
$$B:={}^{T}g\circ f:\SS^{2n+1}U \otimes \ko_{\P(\SS^{2n+1}U)} \stackrel{ }{\longrightarrow} \ko_{\P(\SS^{2n+1}U)}(1),$$ 

avec $\SS^{2n+1}U \otimes \ko_{\P(\SS^{2n+1}U)}$ (resp. $\ko_{\P(\SS^{2n+1}U)}(1)$) qui est muni de l'action canonique $\sigma(t)$ (\mbox{resp. de $\widehat{\sigma(t)}$} l'action duale de $\sigma(t)$). Pour que le morphisme $B$ soit un $\C^{*}$-morphisme il ~faut avoir
$$B_{v}(\sigma(t).v_{1})=t^{q}\widehat{\sigma(t)}.B_{v}(v_{1}),$$

pour tout $v, v_{1}\in \SS^{2n+1}U$, $q$  un entier et
$$B_{v}:\SS^{2n+1}U   \stackrel{ }{\longrightarrow} (\ko_{\P(\SS^{2n+1}U)}(1))_{v}=(\C.v)^{-1}.$$

Alors on obtient $q=2 \gamma + (2n+1)(\beta +\alpha)$. Donc il faut que $\ko_{\P(\SS^{2n+1}U)}(1)$ soit muni de l'action duale $\widehat{\sigma(t)}$ multipli\'ee par le caract\`ere $t^{2 \gamma + (2n+1)(\beta +\alpha)}$. Donc $B$ est le $\C^{*}$-morphisme suivant
$$B:\bigoplus_{i=0}^{2n+1}\ko_{\P(\SS^{2n+1} U)}^{(\gamma + (2n+1)\alpha+i(\beta -\alpha))} \stackrel{ }{\longrightarrow} \ko_{\P(\SS^{2n+1}U)}(1)^{(2 \gamma + (2n+1)(\beta +\alpha))}.$$

Donc on a la $\C^{*}$-monade suivante
$$0\longrightarrow \ko_{\P(\SS^{2n+1}U)}(-1) \stackrel{g }{\longrightarrow} \bigoplus_{i=0}^{2n+1}\ko_{\P(\SS^{2n+1} U)}^{(\gamma + (2n+1)\alpha+i(\beta -\alpha))} \stackrel{B}{\longrightarrow} \ko_{\P(\SS^{2n+1}U)}(1)^{(2 \gamma + (2n+1)(\beta +\alpha))}\longrightarrow 0.$$  

D'apr\`es la d\'efinition 2.2 et la proposition 2.3 en \cite{bah1}, on obtient que
$${\bf Iminvg}(\bigoplus_{i=0}^{2n+1}\ko_{\P(\SS^{2n+1} U)}^{(\gamma + (2n+1)\alpha+i(\beta -\alpha))})=\bigoplus_{i=0}^{2n+1}\ko_{\P^{2n+1}}(\gamma + (2n+1)\alpha+i(\beta -\alpha))=\SS^{2n+1}\ku(\gamma)$$

et
$${\bf Iminvg}(\ko_{\P(\SS^{2n+1}U)}(1)^{(2 \gamma + (2n+1)(\beta +\alpha))})=\ko_{\P^{2n+1}}(2 \gamma + (2n+1)(\beta +\alpha)),$$

et
$$\hspace{0,2 cm} {\bf Iminvg}(\ko_{\P(\SS^{2n+1}U)}(-1))=\ko_{\P^{2n+1}},$$

o\`u $\ku =\ko_{\P^{2n+1}}(\alpha )\oplus \ko_{\P^{2n+1}}(\beta)$. Comme le fibr\'e $Q$ a une $GL(\SS^{2n+1} U)$-action donc il a une $\C^{*}$-action, ce qui donne qu'il est $\C^{*}$-invariant au-dessus de l'action $\overline{\sigma(t)}$. D'apr\`es la proposition \ref{1.4.6}, on obtient
$$\overline{ \sigma}(t)^{*}N(f) \simeq N(f). $$ 

On obtient donc ${\bf Iminvg}(Q)=Q_{\gamma,\alpha,\beta} $ et ${\bf Iminvg}(N(f))=N_{\gamma,\alpha,\beta} $. Alors on ~a la monade suivante
$$0\longrightarrow \ko_{\P^{2n+1}} \stackrel{g }{\longrightarrow} \SS^{2n+1}\ku(\gamma) \stackrel{B}{\longrightarrow} \ko_{\P^{2n+1}}(2 \gamma + (2n+1)(\beta +\alpha))\longrightarrow 0$$

qui nous donne les deux suites exactes suivantes
$$ 0\longrightarrow \ko_{\P^{2n+1}}\stackrel{g}{\longrightarrow} \SS^{2n+1}\ku(\gamma)\longrightarrow Q_{\gamma,\alpha,\beta} \longrightarrow 0,$$

$$0\longrightarrow N_{\gamma,\alpha,\beta} \longrightarrow Q_{\gamma,\alpha,\beta}  \stackrel{^{T}g\circ f}{\longrightarrow} \ko_{\P^{2n+1}}(2 \gamma + (2n+1)(\beta +\alpha)) \longrightarrow 0.$$

Donc on obtient
$$ 0\longrightarrow \ko_{\P^{2n+1}}(-\gamma)\stackrel{g(-\gamma)}{\longrightarrow} \SS^{2n+1}\ku\longrightarrow \kq_{\gamma,\alpha,\beta} \longrightarrow 0,$$

$$0\longrightarrow \kn_{\gamma,\alpha,\beta} \longrightarrow \kq_{\gamma,\alpha,\beta}   \stackrel{(^{T}g\circ f)(-\gamma)}{\longrightarrow} \ko_{\P^{2n+1}}(\gamma + (2n+1)(\beta +\alpha)) \longrightarrow 0,$$

o\`u $\kq_{\gamma,\alpha,\beta}:=Q_{\gamma,\alpha,\beta}(-\gamma)$, et $\kn_{\gamma,\alpha,\beta}:=N_{\gamma,\alpha,\beta}(-\gamma)$.

\end{proof}

\subsubsection{\bf Proposition}\label{1.4.9} 
{\em Soit $\kn:=\kn_{\gamma,\alpha,\beta}$ le fibr\'e de 0-corr\'elation pond\'er\'e 
provenant d'une image inverse g\'en\'eralis\'ee sur $\P^{2n+1}$ qui est d\'efini dans la proposition \ref{1.4.8}.\
Alors le fibr\'e $\kn$ est un fibr\'e symplectique ($\kn^{*}((\alpha+\beta )(2n+1))=\kn $)}.

\begin{proof}
voir l'article \cite{bah2}.
\end{proof}

\subsection{Technique de construction d'un fibr\'e sur une vari\'et\'e lisse}\label{3.1} 

Kumar-Peterson-Rao ont pr\'esent\'e dans leur article \cite{ku-ra-pe} une m\'ethode pour construire un fibr\'e vectoriel sur une vari\'et\'e projective lisse $X$ \`a partir de deux fibr\'es vectoriels sur une sous-vari\'et\'e de $X$.

Soient $X$ une vari\'et\'e projective lisse, $Y$ un diviseur de $X$ correspondant 
\`a une section $s\in~H^{0}(X,\ko_{X}(Y))$ et $i:Y\hookrightarrow X$ le morphisme d'inclusion. 
On notera $D$ le faisceau $i_{*}D$ pour tout faisceau coh\'erent $D$ sur la vari\'et\'e $Y$.

Supposons que l'on ait la suite exacte suivante de fibr\'es vectoriels sur $Y$
$$  0 \longrightarrow \A \stackrel{a }{\longrightarrow} \F \stackrel{b }{\longrightarrow} \B\longrightarrow 0$$ 

et que l'on ait $\hat F$ un fibr\'e vectoriel sur $X$ tel que la restriction 
sur $Y$ est $\hat F|_{Y}=\F$. On d\'efinit $G$ un fibr\'e vectoriel sur $X$ par 
le diagramme commutatif suivant

\xmat{ &0\ar[d]&0\ar[d]&&\\
&\hat F(-Y) \ar@{=}[r] \ar[d]&\hat F(-Y)\ar[d]&&\\
0 \ar[r]&G \ar[r]\ar[d]&\hat F \ar[r]\ar[d]&\B \ar[r] \ar@{=}[d] &0\\
0 \ar[r]&\A \ar[r]\ar[d]&\F \ar[r]\ar[d]&\B \ar[r]&0\\
&0&0&&}

On obtient que le rang du fibr\'e $G$ est $rg(G)=rg(\A)+rg(\B)$. Supposons maintenant qu'il y ait aussi deux fibr\'es vectoriels $ \hat L_{1}$, $\hat L_{2}$ sur $X$ avec $\L_{1}$, $\L_{2}$ leurs restrictions sur $Y$ telles 
qu'il existe un morphisme surjectif de fibr\'es $f: \L_{1}\longrightarrow \A$ et un morphisme 
injectif de fibr\'es $ g: \B \longrightarrow \L_{2}$. Soit le morphisme $f$ (resp. $g$) tel que le compos\'e
$\phi :\L_{1}\stackrel{f }{\longrightarrow} \A \stackrel{a }{\longrightarrow}  \F$ (resp. $\psi: \F \stackrel{b }{\longrightarrow} \B \stackrel{g }{\longrightarrow}  \L_{2}$) sur $Y$ se rel\`eve en un morphisme $\Phi:\hat L_{1}\longrightarrow \hat F $ (resp. $\Psi:\hat F\longrightarrow  \hat L_{2}$ ) sur $X$. Donc on a le diagramme commutatif suivant
\xmat{ \hat L_{1} \ar[rr]^{\Phi}\ar[d]&&\hat F \ar[rr]^{\Psi}\ar[d]&&\hat L_{2} 
\ar[d]\\
 \L_{1} \ar[rr]^{\phi}\ar[rd]_{f} &  & \F \ar[rr]^{\psi}\ar[rd]_{b}  & &  \L_{2}\\
&   \A \ar[ru]_{a}  \ar[rd] & &  \B   \ar[rd]\ar[ru]_{g}   \\
 0\ar[ru] &&   0\ar[ru]    && 0.}

\subsubsection{\bf Lemme} {\em (\cite{ku-ra-pe}, proposition 2.1).  \label{3.1.1} 
On a un morphisme de fibr\'es vectoriels $\Delta$ qui est d\'efini de la fa\c{c}on suivante
$$\Delta:\hat F(-Y)\oplus \hat L_{1}\longrightarrow \hat F\oplus \hat L_{2}(-Y) $$ 
$$
{\Delta} = \left(
\begin{array}{c|c}
sI    &    \Phi \\
 \hline
\Psi  &    s^{-1}\Psi\Phi \\
\end{array}
\right)
.$$
 
Le morphisme $\Delta$ se factorise de la fa\c{c}on suivante
\xmat{  \hat F(-Y)\oplus \hat L_{1}   \ar[rr]^{\Delta}\ar[rd]^{\Delta_{1}} &  
&  
\hat F\oplus \hat L_{2}(-Y)  \\
&   G \ar[ru]^{\Delta_{2}}\ar[rd]\\
0\ar[ru] &&0 }

o\`u $\Delta_{2}:G \longrightarrow \hat F\oplus \hat L_{2}(-Y)$ est un morphisme 
injectif de fibr\'es, et $\Delta_{1}:\hat F(-Y)\oplus \hat L_{1} \longrightarrow G$ 
est un morphisme surjectif de fibr\'es}.

\begin{proof}
 voir \cite{ku-ra-pe}.
\end{proof}
 
\subsubsection{\bf Proposition } {\em (\cite{ku-ra-pe}, proposition 2.2).\label{3.1.2} 
 Soient $\hat N_{1}, \hat N_{2}$ des fibr\'es vectoriels sur $X$ et $ \N_{1}, \N_{2}$ 
leurs restrictions sur $Y$ respectivement tels que  $ \hat F=\hat N_{1}\oplus \hat 
N_{2} $ et $\F=\N_{1}\oplus \N_{2}$. Soit 
$$\pi_{1}:\N_{1}(-Y) \longrightarrow \A$$
un morphisme qui se rel\`eve en un morphisme
$$\pi: \hat N_{1}(-Y) \longrightarrow \hat L_{1}$$

tel que $Im(\Phi \pi)\subset \hat N_{2}$. Le morphisme $\mu : \hat N_{1}(-Y)\longrightarrow G$ est injectif 
si et seulement si \mbox{sa restriction} \`a Y l'est.}

\begin{proof}
On a le diagramme commutatif suivant   
\xmat{& & \hat N_{1}(-Y) \ar@{^{(}->}[lldd]_{i_{1}}\ \ar@{-->}[dddd]^{\mu}\ar[rrdd]^{\mu_{1}}\\
\\
\hat N_{1}(-Y)\oplus \hat N_{2}(-Y) \oplus \hat L_{1} \ar[rrrr]_{\Delta}\ar[rrdd]^{\Delta_{1}} & & & & \hat N_{1}\oplus \hat N_{2}\oplus \hat L_{2}(-Y)\\
\\
& & G\ar[uurr]^{\Delta_{2}}\ar[dr]^{}\\
& 0\ar[ur] & & 0
 }

o\`u $i_{1} $ est un morphisme d'inclusion et le morphisme 
$$\mu_{1} : \hat N_{1}(-Y)\longrightarrow \hat N_{1}\oplus \hat N_{2}\oplus \hat L_{2}(-Y)  $$

est d\'efini par  
$$\mu_{1}: =\left(
\begin{array}{ccccccc}
 
 sI\\
\Phi \pi \\
s^{-1}\Psi\Phi \pi+ \Psi i\\
\end{array}
\right).
$$

o\`u $i:\hat N_{1}\hookrightarrow \hat N_{1}\oplus \hat N_{2}$ est une inclusion. On en d\'eduit que $\mu=\mu_{1}$ et que le morphisme $\mu_{1}$ \mbox{est injectif} sur $X-Y$. Donc le morphisme $\mu : \hat N_{1}(-Y)\longrightarrow G$ est injectif si et seulement si \mbox{sa restriction} \`a Y l'est.

\end{proof}

\subsubsection{\bf Remarque} \label{3.1.3} Dans l'application de la proposition \ref{3.1.2}, on va construire directement le morphisme injectif de fibr\'es $\mu: \hat N_{1}(-Y)\longrightarrow G$ tout en utilisant le fait que: si $ rg( \N_{1})<  rg( \N_{2})$ et $rg(\A)= rg(\B)$ alors on a $rg( \N_{1})\leq rg(\A)-1 $.

\subsubsection{\bf Proposition } {\em (\cite{ku-ra-pe}, proposition 2.3). \label{3.1.4}
 Soient $\hat M_{1}, \hat M_{2}$ des fibr\'es vectoriels sur $X$ et $ \M_{1}, \M_{2}$ 
leurs restrictions sur $Y$ respectivement tels que  $ \hat F=\hat M_{1}\oplus \hat 
M_{2} $ et $\F=\M_{1}\oplus \M_{2}$. Soit 
$$d_{1}:\B(-Y) \longrightarrow \M_{1}$$

un morphisme qui se rel\`eve en un morphisme
$$d: \hat L_{2}(-Y) \longrightarrow \hat M_{1}$$

tel que $d\Psi $ s'annule sur $ \hat M_{1}(-Y)$. Le morphisme $\delta :G \longrightarrow \hat M_{1}$ est surjectif si et seulement si sa restriction \`a Y l'est.}

\begin{proof}
On a le diagramme commutatif suivant   
 \xmat{& & \hat M_{1} \\
\\
\hat M_{1}(-Y)\oplus \hat M_{2}(-Y) \oplus \hat L_{1} \ar[rrrr]_{\Delta}\ar[rrdd]^{\Delta_{1}}\ar[rruu]^{\delta_{1}} & & & & \hat M_{1}\oplus \hat M_{2}\oplus \hat L_{2}(-Y)\ar@{->>}[lluu]_{p_{1}}\\
\\
& & G\ar[uurr]^{\Delta_{2}}\ar[dr]^{} \ar@{-->}[uuuu]_{\delta}\\
& 0\ar[ur] & & 0
 }

o\`u $p_{1} $ est un morphisme projection et le morphisme 
$$\delta_{1} : \hat M_{1}(-Y)\oplus \hat M_{2}(-Y)\oplus \hat L_{1} \longrightarrow  \hat M_{1}   $$

est d\'efini par  
$$\delta_{1}: =\left(
\begin{array}{ccccccc}
 
 sI && d\Psi  & & d s^{-1}\Psi\Phi  + p\Phi
\end{array}
\right).
$$

o\`u $p:\hat M_{1}\oplus \hat M_{2}\longrightarrow \hat M_{1}$ est une projection. On en d\'eduit que $\delta=\delta_{1}$ et que le morphisme $\delta_{1}$ est surjectif sur $X-Y$. Donc le morphisme $\delta :G \longrightarrow \hat M_{1}$ est surjectif si et seulement si sa restriction \`a Y l'est. 
\end{proof}

\subsubsection{\bf Th\'eor\`eme de  Hirzebruch-Riemann-Roch} {\em (Hartshorne, 
\cite{ha}).  
Soient $X$ une vari\'et\'e projective lisse sur un corps $\K$ et $F$ un faisceau coh\'erent sur $X$. On a donc
le \mbox{diagramme} commutatif suivant

\xmat{ K(X) \ar[r]^{ch(\bullet).td(X)}  \ar[d]_{\chi(X,\bullet)} &  
CH^*(X)\ar[d]^{deg(\bullet)} \\
\Z \ar[r] &  \Q}
$\chi(X,F) = deg(ch(F).td(X))$

o\`u $deg(\bullet)$ est l'application de degr\'e, $ch(F)$ le caract\`ere de Chern de $F$, $K(X)$ l'anneau de Grothendieck de $X$, $CH^*(X)$ l'anneau de Chow de $X$ et $td(X)$ les classes de Todd de $X$}.

\begin{proof}
 voir \cite{ha}.
\end{proof}

Soit $f: X\longrightarrow Y$ un morphisme de vari\'et\'es alg\'ebriques. A l'aide du th\'eor\`eme de Hirzebruch-Riemann-Roch, on peut calculer les classes de Chern de l'image directe d'un fibr\'e $F$ d\'efini sur la vari\'et\'e $X$.

\subsubsection{\bf D\'efinition} 
{\em Un morphisme $f: X\longrightarrow Y$ est un morphisme fini si et seulement s'il existe un recouvrement des ouverts affines $V_{i}=Spec(B_{i})$ de $Y$ tels que $f^{-1}(V_{i})=Spec(A_{i})$ \mbox{o\`u $A_{i}$ est une} $B_{i}$-alg\`ebre de type fini, pour tout $i$}.  

\subsubsection{\bf Lemme }{\em (Le Potier \cite{lep}). \label{3.1.7} 
 Soient  $f:X\longrightarrow Y$ un morphisme fini de vari\'et\'es alg\'ebriques
 et $F$ un faisceau coh\'erent  sur $X$.

1 - Le faisceau $f_*F$ est coh\'erent et $R^qf_*F=0$ pour $q>0$.

2 - On a $H^q(X,F)=H^q(Y,f_*F)$ pour tout int\`egre $q$. En particulier si $X$ et $Y$ 
sont des vari\'et\'es projectives, alors on a
$$\chi(X,F(t))=\chi(Y,f_*F(t)) \hspace{0.2 cm} pour \hspace{0.2 cm} t\in \Z.$$}

\begin{proof}
 voir \cite{lep}.
\end{proof}

\subsubsection{\bf Proposition }\label{3.1.8} 
{\em Soient $z \in H^{0}(\P^{n},\ko_{\P^{n}}(1))\setminus \lbrace 0 \rbrace $ et $\P^{n-1}   \subset  \P^n $ l'hyperplan de $\P^n$ d\'efini par l'\'equation $z=0$. Soient $\epsilon\geq 1$ un entier positif et $Y_{(\epsilon)}\subset \P^n$ le voisinage infinit\'esimal de l'ordre $\epsilon$ de $\P^{n-1} $ dans $\P^n$ d\'efini par l'\'equation $z^{\epsilon}=0$.
Soient les morphismes d'inclusion suivants
 
 \xmat{  \P^{n-1} :=Y_{(1)}  \ar@{^{(}->}[rr]^{b_{(2)}} &  & Y_{(2)} \ar@{^{(}->}[rr]^{b_{(3)}} &  &  \ldots \ar@{^{(}->}[rr]^{b_{(\epsilon)}} &  &   Y_{(\epsilon)}\ar@{^{(}->}[rr]^{d} &  &  \P^n . \\
 }
 
Soient $\widehat{F}$ un fibr\'e vectoriel sur $\P_{n-1}$ de rang $r$ dont les classes de Chern sont $c_1,c_2,c_3, \ldots ,c_{n-1}$, $F_{(\epsilon)}$ un faisceau coh\'erent sur $Y_{(\epsilon)}$ tel que $\widehat{F}=F_{(\epsilon)| \P_{n-1}}:=F_{(1)}$ et $F=d_{*}(F_{(\epsilon)})$ un faisceau coh\'erent sur $\P_{n}$ dont les classes de Chern sont $c^{'}_{1},c^{'}_{2},c^{'}_{3}, \ldots ,c^{'}_{n}$. Alors    
$$c^{'}_{1}=r\epsilon,\hspace{0.2 cm}  c^{'}_{2}=\epsilon^{2}(\frac{r(r+1)}{2})-\epsilon c_1$$

$$c^{'}_{3}=\epsilon^{3}(\frac{r(r+1)(r+2)}{6})-\epsilon^{2}(r+1)c_1+\epsilon (c_{1}^{2}-2c_2)$$

$$c^{'}_{4}=\epsilon^{4}(\frac{r(r+1)(r+2)(r+3)}{24})-\epsilon^{3}( \frac{(r+2)(r+1)}{2})c_1+$$

$$\epsilon^{2}((r+2)c_{1}^{2}-(2r+3)c_2)+\epsilon (-c_{1}^{3}+3c_1c_2-3c_3).$$}

\begin{proof}D'apr\`es l'exemple $15.1.2$ dans \cite{fu}: Pour $E$ un faisceau coh\'erent de rang $r$ dont les classes de Chern sont $c_{1},c_{2},\ldots ,c_{n}$ sur $\P^{n}=\P(V)$, on a
$$ch(E)=r+h_{1}+\frac{1}{2}h_{2}+ \frac{1}{3!}h_{3}+ \ldots +\frac{1}{n!}h_{n}+ \ldots$$ 

o\`u 
 $$
{h_i} =det \left(
\begin{array}{ccccc}
c_1    & 1      &   0    & \cdots & 0 \\
2 c_2  & c_1    &  1     &\cdots  & 0 \\
\vdots & \vdots &\vdots  &\ddots  & \vdots \\
i c_i & c_{i-1} & c_{i-2}& \cdots & c_1 \\
\end{array}
\right)
.$$

et 
$$td(\P^{n})=1+T_{1}+T_{2}+T_{3}+ \ldots  +T_{n}+ \ldots$$

o\`u 
$$ T_{0}=1, T_{1}= \frac{1}{2} \alpha_{1}, T_{2}=\frac{1}{12}(\alpha_{1}^{2}+ \alpha_{2}), T_{3}= \frac{1}{24} (\alpha_{1} \alpha_{2}), \ldots$$

o\`u $\alpha_{1},\alpha_{2},\ldots ,\alpha_{n}$ sont les classes de Chern du fibr\'e $Q(1)$ qui est d\'efini par la suite exacte suivante
$$0\longrightarrow \ko_{\P^{n}}(-1) \longrightarrow V\otimes \ko_{\P^{n}} \longrightarrow Q \longrightarrow 0.$$

D'apr\`es le th\'eor\`eme de Hirzebruch-Riemann-Roch sur la vari\'et\'e $\P^{n}$, on a 
$$\chi(\P^{n},E(t))=\frac{r }{n!} .t^{n}+\sum_{i=1}^{n}\frac{A_{i} }{(n-i)!}. t^{n-i}$$

o\`u $A_{i}=rT_{i}+\sum _{j=0}^{i-1} T_{i-j-1}.h_{j+1}$, $ i=1,2,\ldots ,n$ et pour tout $t\in \Z$. De la suite exacte sur $Y_{(\epsilon)}$
$$0\longrightarrow F_{(\epsilon-1)}(-1)\longrightarrow F_{(\epsilon)}\longrightarrow F_{(1)}\longrightarrow 0,$$

on obtient alors que 
$$\chi(Y_{(\epsilon)},F_{(\epsilon)}(t))=\sum_{k=0}^{\epsilon-1}\chi(\P^{n-1},F_{(1)}
(t-k)).$$

Comme l'immersion ferm\'ee affine $d$ est un morphisme fini, en utilisant le lemme \ref{3.1.7}, on obtient que 
$$\sum_{k=0}^{\epsilon-1}\chi(\P^{n-1},F_{(1)}(t-k))=\chi(Y_{(\epsilon)},F_{(\epsilon)}(t))=\chi(\P^{n},F(t))  .$$

Alors on en d\'eduit
$$\sum_{k=0}^{\epsilon-1}\lbrace \frac{ r}{(n-1)!}.(t-k)^{n-1}+\sum_{i=1}^{n-1}
\frac{A^{'}_{i}}{(n-1-i)!}. (t-k)^{n-1-i}\rbrace=$$ 

$$\sum_{k=0}^{\epsilon-1}\lbrace \frac{r }{(n-1)!} .t^{n-1}+ 
(\frac{r.k }{(n-2)!}+\frac{A_{1}^{'}}{(n-2)!}) .t^{n-2}  + $$

$$(\frac{r.k^{2} }{2!.(n-3)!}-\frac{A_{1}^{'}.k}{(n-3)!}+\frac{A_{2}^{'}}{(n-3)!}) .t^{n-3}+$$

$$(-\frac{r.k^{3} }{3!.(n-4)!}+\frac{A_{1}^{'}.k^{2}}{(n-4)!}-\frac{A_{2}^{'}.k }{(n-4)!}+\frac{A_{3}^{'}}{(n-4)!}) .t^{n-4}
\rbrace +\ldots =$$

$$=\frac{rg(F)}{n!} .t^{n}+ \frac{ A_{1}}{ (n-1)! } .t^{n-1}+ \frac{ A_{2}}{ (n-2)! } .t^{n-2}+ \frac{ A_{3}}{ (n-3)! } .t^{n-3}+ \frac{ A_{4}}{ (n-4)! } .t^{n-4}+  \ldots$$

Ces polyn\^omes \'etant les m\^emes, on peut identifer les coefficients et obtenir des
\'equations permettant d'obtenir  $rg(F)=0$ et 
$$c^{'}_{1}=r\epsilon,\hspace{0.2 cm}  c^{'}_{2}=\epsilon^{2}(\frac{r(r+1)}{2})-\epsilon c_1$$

$$c^{'}_{3}=\epsilon^{3}(\frac{r(r+1)(r+2)}{6})-\epsilon^{2}(r+1)c_1+\epsilon(c_{1}^{2}-2c_2)$$

$$c^{'}_{4}=\epsilon^{4}(\frac{r(r+1)(r+2)(r+3)}{24})-\epsilon^{3}( \frac{(r+2)(r+1)}{2})c_1+$$

$$\epsilon^{2}((r+2)c_{1}^{2}-(2r+3)c_2)+\epsilon(-c_{1}^{3}+3c_1c_2-3c_3).$$

\end{proof}

\section{ \large\bf Fibr\'e vectoriel de rang $2n+1$ sur $\P^{2n+2}$ }

\vspace{1cm}
 
\subsection{ Fibr\'e vectoriel de rang $2n+1$ sur $\P^{2n+2}$ } 

Pour obtenir un fibr\'e vectoriel de rang $2n+1$ sur la vari\'et\'e projective $\P^{2n+2}$, nous allons \'etablir la relation entre le fibr\'e de Tango pond\'er\'e provenant d'une image inverse g\'en\'eralis\'ee en \cite{bah1} et le fibr\'e de 0-corr\'elation pond\'er\'e provenant d'une image inverse g\'en\'eralis\'ee sur $\P^{2n+1}$ en \ref{1.4.5}.

\subsubsection{\bf Proposition}\label{3.2.1} 

{\em Soient $V$ un $\C$-espace vectoriel de dimension $dim(V)=2n+2$, et $f \in \bigwedge^{2}V^{*} =H^{0}(Q^{*}(1))$ une forme bilin\'eaire antisym\'etrique non d\'eg\'en\'er\'ee qui correspond \`a l'hyperplan $H\subset \bigwedge^{2} V=(H^{0}(Q^{*}(1)))^{*} $, et $W\subset \bigwedge^{2}V=(H^{0}(Q^{*}(1)))^{*} $ un sous-espace vectoriel v\'erifiant la condition (*) en \cite{bah1}. Soient $N(f) \hspace{0.2cm} et\hspace{0.2cm} F(W)$ le fibr\'e de corr\'elation nulle classique et le fibr\'e de Tango correspondants respectivement \`a $f$ et \`a $W$. Si $W\subset H$, Alors on a la suite exacte suivante
$$0\longrightarrow N(f)(-1)\longrightarrow  \left(    H\diagup W \right)\otimes\ko_{\P^{2n+1}}  \longrightarrow  F(W)(1)   \longrightarrow 0.$$}
 
\begin{proof} Le fibr\'e $N(f)$ est d\'efini par la suite exacte suivante
$$0\longrightarrow \ko_{\P^{2n+1}}\stackrel{^{T} f\circ g }{\longrightarrow}Q^{*}(1)\longrightarrow N(f)(1) \longrightarrow 0 .$$ 

Alors on a la suite exacte de la cohomologie
$$0\longrightarrow \C f \stackrel{H^{0}(^{T} f\circ g) }{\longrightarrow} H^{0}(Q^{*}(1))\simeq \bigwedge^{2}V^{*} \longrightarrow H^{0}(N(f)(1)) \longrightarrow 0 $$

ou bien
$$0\longrightarrow \C  \stackrel{H^{0}(^{T} f\circ g) }{\longrightarrow} H^{0}(Q^{*}(1))\simeq \bigwedge^{2}V^{*} \longrightarrow H^{0}(N(f)(1)) \longrightarrow 0 $$

O\`u $H^{0}(^{T} f\circ g)(a):= a.f$ pour tout $a \in \C$. Alors on a $H^{0}(^{T} f\circ g)(1)= f$, c'est-\`a-dire on a le diagramme commutatif suivant

\xmat{   && &&0 \\
0 \ar[rr]&& \ko_{\P^{2n+1}}  \ar[rr]^{^Tf\circ g }\ar[rrd]_{-^Tf} &&Q^{*}(1)\ar[u]   \\
 &&    &&\bigwedge^{2}V^{*}  \otimes\ko_{\P^{2n+1}}\ar[u]_{ev_{Q^{*}(1)}}  \\
   }

dont le dual est

\xmat{   0\ar[d]&& && \\
Q(-1)  \ar[rr]^{^Tg\circ f}\ar[d] && \ko_{\P^{2n+1}}  \ar[rr]&&0 \\
   \bigwedge^{2}V \otimes\ko_{\P^{2n+1}}\ar[rru]^{f}.}
   
Comme $f:\bigwedge^{2}V \longrightarrow \C$ est une application lin\'eaire surjective, on consid\`ere la suite exacte suivante
$$0\longrightarrow H \longrightarrow \bigwedge^{2}V \stackrel{f  }{\longrightarrow} \C \longrightarrow 0 .$$ 

Si $W\subset H$, alors on obtient le diagramme commutatif suivant 

\xmat{   0 && && \\
 \left(    (\bigwedge^{2}V)\diagup W \right) \ar[rr]^{ \overline{f}}\ar[u] &&\C  \ar[rr]&&0 \\
   \bigwedge^{2}V  \ar[rru]^{f} \ar[u],  \\
   }

et $ker(\overline{f})=H\diagup W$. Alors on a la suite exacte suivante
 $$0\longrightarrow H\diagup W\otimes\ko_{\P^{2n+1}}\longrightarrow \left(    ( \bigwedge^{2}V)\diagup W\right) \otimes\ko_{\P^{2n+1}}  \stackrel{\overline{f}}{\longrightarrow}   \otimes\ko_{\P^{2n+1}} \longrightarrow 0.$$

On obtient donc le diagramme commutatif suivant   

\xmat{   0\ar[d]&& && \\
Q(-1)  \ar[rr]^{^Tg\circ f}\ar[d] && \ko_{\P^{2n+1}}\ar@{=}[d] \ar[rr]&&0 \\
\bigwedge^{2}V \otimes\ko_{\P^{2n+1}}\ar[rr]^{f}\ar[d] &&  \ko_{\P^{2n+1}} \ar[rr]&& 0 \\
  \left(   (\bigwedge^{2}V)\diagup W  \right) \otimes\ko_{\P^{2n+1}}\ar[rru]^{\overline{f}} }

O\`u $Q(-1)\longrightarrow  \left(    (\bigwedge^{2}V)\diagup W\right) \otimes\ko_{\P^{2n+1}} $ est le morphisme injectif $\varpi_{W}$ dans la suite exacte
$$0\longrightarrow Q(-1)\stackrel{\varpi_{W}}{\longrightarrow}  \left(    (\bigwedge^{2}V)\diagup W\right) \otimes \ko_{\P^{2n+1}}  \longrightarrow  F(W)(1)   \longrightarrow 0.$$

Donc on a le diagramme commutatif suivant

\xmat{ &0\ar[d]&0\ar[d]&&\\
0 \ar[r]&N(f)(-1) \ar[r]\ar[d]^{h}&Q(-1) \ar[r]^{^{T}g\circ f}\ar[d]^{\varpi_{W}}&   \ko_{\P^{2n+1}} \ar[r]\ar@{=}[d]&0\\
0 \ar[r]&H\diagup W\otimes\ko_{\P^{2n+1}} \ar[r]\ar[d]&  \left(    (\bigwedge^{2}V)\diagup W\right) \otimes\ko_{\P^{2n+1}} \ar[r] ^{\overline{f}}\ar[d]&  \ko_{\P^{2n+1}} \ar[r]&0\\
&coker(h)\ar[r] \ar[d]&F(W)(1)\ar[d]&&\\
&0&0&&}

D'apr\`es le lemme de serpent, on obtient que $coker(h)\simeq F(W)(1)$ qui donne la suite exacte suivante
$$0\longrightarrow N(f)(-1)\longrightarrow H\diagup W\otimes\ko_{\P^{2n+1}}  \longrightarrow  F(W)(1)   \longrightarrow 0.$$
  
\end{proof}

A partir de maintenant, nous allons nous int\'eresser au fibr\'e de corr\'elation nulle classique sp\'ecial $N(f)$ sur $\P^{2n+1}$.
 
\subsubsection{\bf Dualit\'e }\label{3.2.2}  Soient $U= \C^{2}$ un espace vectoriel, et $\lbrace x,y\rbrace$ sa base. Il existe un isomorphisme $ \nu :U\longrightarrow U^{*}$ d\'efini \`a une constante multiplicative pr\`es par la matrice antisym\'etrique 
$$
\left(
\begin{array}{cccccccccc}
   0 & -1 \\
 
1 & 0   \\
\end{array}
\right),$$

on obtient que $x^{*}:=^Tx\nu=-y$, $y^{*}:=^Ty\nu=x $ est la base de $U^{*}$. La forme symplectique $\langle ,\rangle$ sur $U$ correspondant \`a $\nu$ est donc d\'efini par 
$$\langle x,x \rangle=\langle y,y \rangle=0,\hspace{0.2cm} \langle x,y \rangle=-1.$$

On appelle l'espace $(U;\langle ,\rangle)$ le plan symplectique. On en d\'eduit l'isomorphisme $\SS^{m}U\simeq \SS^{m}U^{*}$ et la forme quadriatique (resp. la forme symplectique) associ\'ee, si $m$ est pair (resp. si $m$ est impair), d\'efinie par 
$$\langle x^{i}y^{m-i},x^{j}y^{m-j}  \rangle:= \left\{
\begin{array}{ccc}

0&:& i+j\neq m \\
&&\\

 (-1)^{i}i!. j!&:& i+j= m. \\

 \end{array}
\right.
$$

qui est une $(m+1)\times(m+1)$-matrice $B_{m}$. Ensuite on en d\'eduit l'isomorphisme $\bigwedge^{2}\SS^{2n+1}U \simeq \bigwedge^{2}\SS^{2n+1}U^{*}$ et la forme quadriatique sym\'etrique associ\'ee d\'efinie par: soient $0\leq i< j\leq 2n+1$, $0\leq s< t\leq 2n+1$, et
$$z_{i,j}=x^{i}y^{2n+1-i}\wedge x^{j}y^{2n+1-j}, \hspace{0.2cm} z_{s,t}=x^{s}y^{2n+1-s}\wedge x^{t}y^{2n+1-t} $$

alors on a
$$\langle z_{i,j},z_{s,t}  \rangle:= \left\{
\begin{array}{ccc}

0&:& j+s\neq 2n+1\hspace{0.2cm}  ou  \hspace{0.2cm} i+t\neq 2n+1  \\
&&\\

 (-1)^{i+j+1}i!. t!. j!.s!&:& j+s= 2n+1\hspace{0.2cm}  et \hspace{0.2cm} i+t= 2n+1. \\

 \end{array}
\right.
$$

Pour toute forme lin\'eaire $f \in (\bigwedge^{2}\SS^{2n+1}U)^{*}$, on a que l'action de cette forme est
$$f(g):=\langle f,g  \rangle$$ 

o\`u $g \in \bigwedge^{2}\SS^{2n+1}U$. Ensuite on en d\'eduit une forme non d\'eg\'en\'er\'ee sur l'espace vectoriel $\bigoplus_{i=0}^{m}\SS^{i}U$ d\'efinie par la matrice
$$B = diag(B_{m},B_{m-1},B_{m-2},\ldots ,B_{0}).$$

 \subsubsection{\bf Th\'eor\`eme} \label{3.2.3} 
{\em En gardant les notations de la proposition \ref{1.4.8}. Soient $U= \C^{2}$ le plan symplectique et $n>0$ un entier, et $\P(\SS^{2n+1}U)$ l'espace projectif associ\'e \`a l'espace vectoriel $\SS^{2n+1}U$. Pour tout $W \in \WW$, comme dans le th\'eor\`eme 3.4 en \cite{bah1} pour $\P(\SS^{2n+1}U)$, il existe une forme antisym\'etrique symplectique non d\'eg\'en\'er\'ee $f:=f_{W}:\bigwedge^{2}\SS^{2n+1}U \longrightarrow \C$ telle que 
$$^{T}\sigma_{1} f\sigma_{1}= af, \text{ et } \langle  f,w \rangle=0 \text{ pour tout } w\in W,$$

o\`u $a\in \C^{*}$ et $ \sigma_{1}=\left(\begin{array}{cc } t^{\alpha}&0 \\0&    t^{\beta} \\ \end{array} \right),\hspace{0.2 cm} t\in \C^{*} $. Il existe aussi $D_{W}$ un sous-espace vectoriel $\C^{*}$-invariant de $\bigwedge^{2}\SS^{2n+1}U$ tel que $\bigwedge^{2}\SS^{2n+1}U\simeq D_{W}\oplus W $. Dans ce cas on obtient la suite exacte suivante 
$$0\longrightarrow N(f)(-1)\longrightarrow  \left( D_{W} \diagup  \C.f  \right)  \otimes\ko_{\P(\SS^{2n+1}U)}  \longrightarrow  F(W)(1)   \longrightarrow 0,$$

o\`u $N(f)$ est le fibr\'e de corr\'elation nulle classique sp\'ecial associ\'e \`a la forme $f$, et $F(W)$ est le fibr\'e de Tango associ\'e au sous-espace $W$.}

\begin{proof} On va utiliser les m\^emes notations de 3.3 et 3.4 en \cite{bah1}. Soit 
$$\BB:=\lbrace  z_{p,q}:=x^{2n+1-p}y^{p}\wedge x^{2n+1-q}y^{q}  ,  \hspace{0.2 cm} 0\leq p< q \leq 2n+1 \rbrace, $$
 
la base de l'espace vectoriel $\bigwedge^{2}\SS^{2n+1}U$. Comme dans la remarque 3.3 en \cite{bah1} pour $\P(\SS^{2n+1}U)$, on a  
$$\bigwedge^{2}\SS^{2n+1}U \simeq \bigoplus_{1\leq k\leq 4n+1}E_{k},$$

o\`u $E_{k}$ est engendr\'e par les \'el\'ements $z_{p,q}$ tels que $k=p+q$. On pose que
$$\E_{k}:=\left\{
\begin{array}{ccc}

E_{k}\oplus E_{4n+2-k}& :&1\leq k\leq 2n \\
\\
E_{2n+1}
 \end{array}
\right.
$$

qui donne
$$\bigwedge^{2}\SS^{2n+1}U \simeq \bigoplus_{1\leq k\leq 2n+1}\E_{k}.$$

On obtient que $\langle u_{i},u_{j}  \rangle=0$ pour tout $u_{i}\in \E_{i}$ et $u_{j}\in \E_{j}$ avec $i\neq j$. Donc si $ u,h \in \bigwedge^{2}\SS^{2n+1}U$ on peut les \'ecrire sous la forme
$$u=u_{1}+u_{2}+ \ldots +u_{2n}+u_{2n+1},  \hspace{0.2 cm}  h=h_{1}+h_{2}+ \ldots +h_{2n}+h_{2n+1},$$

o\`u $u_{k},h_{k} \in \E_{k}$ pour $1\leq k\leq 2n+1$. Apr\`es la dualit\'e \ref{3.2.2}, on obtient  que
$$\langle u,h \rangle=\langle u_{1},h_{1}\rangle+ \langle u_{2},h_{2}\rangle+\langle u_{3},h_{3}\rangle+ \ldots +\langle u_{2n},h_{2n}\rangle+ \langle u_{2n+1},h_{2n+1}\rangle.$$

On cherche une forme lin\'eaire $f\in (\bigwedge^{2}\SS^{2n+1}U)^{*} $ telle que 
$$^{T}\sigma_{1} f\sigma_{1}= af, \text{ et } f(w):=\langle f,w\rangle=0$$

pour tout $w\in W$ et $a\in \C^{*}$. Ce qui est donc \'equivalent \`a dire que l'on cherche une forme lin\'eaire telle que $f\in E_{2n+1} $, et 
$$f(w_{2n+1}):=\langle f,w_{2n+1} \rangle=0$$

pour tout $w_{2n+1}\in W_{2n+1}$. Comme $W_{2n+1}$ est un hyperplan de $E_{2n+1}$, on a dans ce cas 
$$\C.f=W_{2n+1}^{\vee},$$

o\`u $W_{2n+1}^{\vee}$ est l'orthogonal de $W_{2n+1}$ dans $E_{2n+1}$. Donc $f$ est uniquement d\'etermin\'ee \`a une constante multiplicative pr\`es. On obtient que la forme lin\'eaire $f$ est non d\'eg\'en\'er\'ee, car si le coefficient de $z_{p,2n+1-p}$ dans $f$ est nul, pour $0\leq p\leq n$, alors on a
$$f(z_{p,2n+1-p}):=\langle f,z_{p,2n+1-p} \rangle=0$$

qui donne $z_{p,2n+1-p} \in W_{2n+1}$; ce qui est une contradiction \`a $W_{2n+1}\in \ZZ_{2n+1}$. Comme la forme quadriatique $ \langle ,\rangle$ est une forme sym\'etrique non d\'eg\'en\'er\'ee d\'efinie sur $\bigwedge^{2}\SS^{2n+1}U$ et comme $f\neq 0$, alors on a $ f\notin W_{2n+1}$. Pour tout $W\in \WW$, on choisit donc $d_{k}\in E_{k}\diagdown W_{k}$ non nuls pour tout $1\leq k\leq 4n+1$ tel que $\C.d_{2n+1}=W_{2n+1}^{\vee}$ et on d\'efinit $f=a.d_{2n+1}$ o\`u $a \in \C^{*}$. En consid\'erant le sous-espace vectoriel de $\bigwedge^{2}\SS^{2n+1}U$  
$$D_{W}= \bigoplus_{1\leq k\leq 4n+1}\C.d_{k},$$ 

on obtient que $\bigwedge^{2}\SS^{2n+1}U\simeq D_{W}\oplus W $ est un $\C^{*}$-isomorphisme. Si l'hyperplan $H :=ker(f)\subset  \bigwedge^{2}\SS^{2n+1}U$, alors on obtient $W\subset H$. D'apr\`es la proposition \ref{3.2.1}, on obtient la suite exacte suivante 
$$0\longrightarrow N(f)(-1)\longrightarrow H\diagup W \otimes\ko_{\P(\SS^{2n+1}U)}  \longrightarrow  F(W)(1)   \longrightarrow 0,$$

ou bien 
$$0\longrightarrow N(f)(-1)\longrightarrow \left( \bigoplus_{1\leq k \leq 4n+1 }^{k\neq 2n+1}  \C.d_{k}  \right)  \otimes \ko_{\P(\SS^{2n+1}U)}  \longrightarrow  F(W)(1)   \longrightarrow 0,$$

o\`u $N(f)$ est le fibr\'e de corr\'elation nulle classique sp\'ecial associ\'e \`a la forme $f$, et $F(W)$ est le fibr\'e de Tango associ\'e au sous-espace $W$.

\end{proof}

\subsubsection{\bf Proposition}\label{3.2.5} 
{\em En gardant les notations de la proposition \ref{1.4.8}. Soit $U= \C^{2}$ le plan symplectique comme dans la dualit\'e \ref{3.2.2}. Soient $g_{0}, \ldots ,g_{2n+1}$ des formes homog\`enes sans un z\'ero commun sur $\P(\SS^{2n+1}U)$ telles que 
$$ deg(g_{i})=\gamma + (2n+1)\alpha+i(\beta -\alpha),\hspace{0.2 cm} i=0,1, \ldots ,2n+1.$$

Soient $N(f)$ le fibr\'e de corr\'elation nulle classique sp\'ecial sur $\P(\SS^{2n+1}U)$ et $F(W)$ le fibr\'e de Tango sur $\P(\SS^{2n+1}U)$ comme dans le th\'eor\`eme \ref{3.2.3}
$$0\longrightarrow N(f)(-1)\longrightarrow D_{1}  \otimes \ko_{\P(\SS^{2n+1}U)}  \longrightarrow  F(W)(1)   \longrightarrow 0.$$  

o\`u $D_{1}:= \bigoplus_{1\leq k \leq 4n+1 }^{k\neq2n+1}  \C.d_{k}$. Alors on a la suite exacte suivante

\begin{equation}  
0\longrightarrow \kn(-\gamma)\longrightarrow \Upsilon_{1} \longrightarrow  \kf(\gamma) \longrightarrow 0.
\end{equation} 

o\`u $\Upsilon_{1}:=\bigoplus_{1\leq k \leq 4n+1 }^{k\neq 2n+1} \ko_{\P^{2n+1}}(2\alpha(2n+1)+k(\beta -\alpha))$, et $\kn$ (resp. $\kf$) est le fibr\'e d\'efini sur $\P^{2n+1}$ comme dans la proposition \ref{1.4.8} (resp. la Proposition 3.7 en \cite{bah1}).}

\begin{proof}
D'apr\`es la proposition \ref{1.4.6} et le th\'eor\`eme \ref{3.2.3}, on obtient que le fibr\'e $N(f)$ est $\C^{*}$-invariant par rapport \`a l'action $\overline{\sigma(t) }$ et
$$\overline{ \sigma}(t)^{*}N(f) \simeq N(f). $$ 

D'apr\`es la proposition $3.2$ et le th\'eor\`eme $3.4$ en \cite{bah1}, on obtient que l'on obtient que le fibr\'e $F(W)$ est $\C^{*}$-invariant par rapport \`a l'action $\overline{\sigma(t) }$ et 
$$\overline{ \sigma}(t)^{*}F(W) \simeq F(W). $$ 

On consid\`ere la transform\'ee de Horrocks $\bf Iminvg$, qui est d\'efinie dans la proposition \ref{1.4.8}. D'apr\`es la d\'efinition $2.2$ en \cite{bah1}, on obtient ${\bf Iminvg}(F(W))= \kf_{\gamma,\alpha,\beta}$ et ${\bf Iminvg}(N(f))=\kn_{\gamma,\alpha,\beta} $. 
Le fibr\'e $\left( \bigoplus_{1\leq k \leq 4n+1 }^{k\neq2n+1}  \C.d_{k}  \right)  \otimes \ko_{\P(\SS^{2n+1}U)}$ est muni de l'action canonique $\overline{\sigma(t)}$. Comme on a, pour tout $t\in \C^{*}$,
$$\sigma(t).(d_{k}) =t^{ ((2\alpha(2n+1)+k(\beta -\alpha) +2\gamma)}.(d_{k}), $$
 
alors le sous-fibr\'e $(d_{k}.\C)\otimes \ko_{\P(\SS^{2n+1} U)} $ de $\left( \bigoplus_{1\leq k \leq 4n+1 }^{k\neq2n+1}  \C.d_{k}  \right)  \otimes \ko_{\P(\SS^{2n+1}U)}$ est $\C^{*}$-invariant, et on a
$$\ko_{\P(\SS^{2n+1} U)}^{((2\alpha(2n+1)+k(\beta -\alpha) +2\gamma)}\simeq (d_{k}.\C)\otimes \ko_{\P(\SS^{n} U)} . $$

Ce qui est d\'efini localement, pour tout $v\in \SS^{2n+1} U $, par
$$(\ko_{\P(\SS^{2n+1} U)}^{((2\alpha(2n+1)+k(\beta -\alpha) +2\gamma)})_{v} \simeq \C \stackrel{\simeq}{\longrightarrow}((d_{k}.\C)\otimes \ko_{\P(\SS^{2n+1} U)})_{v}\simeq d_{k}.\C  . $$
$$a\longmapsto a d_{k}$$

Donc on a un $\C^{*}$-isomorphisme
$$\left( \bigoplus_{1\leq k \leq 4n+1 }^{k\neq 2n+1}  \C.d_{k}  \right)  \otimes \ko_{\P(\SS^{2n+1}U)} \simeq  \bigoplus_{1\leq k \leq 4n+1 }^{k\neq 2n+1}\ko_{\P(\SS^{2n+1} U)}^{((2\alpha(2n+1)+k(\beta -\alpha) +2\gamma)}. $$

Donc on obtient   
$${\bf Iminvg}\left(\left( \bigoplus_{1\leq k \leq 4n+1 }^{k\neq 2n+1}  \C.d_{k}  \right)  \otimes \ko_{\P(\SS^{2n+1}U)}\right)=\bigoplus_{1\leq k \leq 4n+1 }^{k\neq 2n+1} \ko_{\P^{2n+1}}(2\alpha(2n+1)+k(\beta -\alpha) +2\gamma):=\Upsilon_{1}(2\gamma),$$
 
o\`u $\Upsilon_{1}:=\bigoplus_{1\leq k \leq 4n+1 }^{k\neq 2n+1} \ko_{\P^{2n+1}}(2\alpha(2n+1)+k(\beta -\alpha))$. Nous appliquons la transform\'ee de Horrocks $\bf Iminvg$ sur la suite exacte suivante
$$0\longrightarrow N(f)(-1)\longrightarrow D_{1}  \otimes \ko_{\P(\SS^{2n+1}U)}  \longrightarrow  F(W)(1)   \longrightarrow 0,$$ 
 
on obtient la suite exacte suivante 

$$0\longrightarrow \kn(-\gamma)\longrightarrow \Upsilon_{1} \longrightarrow  \kf(\gamma) \longrightarrow 0.$$ 

o\`u $\kn:=\kn_{\alpha,\gamma}(-\gamma),\hspace{0.2 cm} \kf(\gamma):=\kf_{\alpha,\gamma}(-2\gamma)$. 

\end{proof}

On consid\`ere
$$  \Upsilon_{1}=\bigoplus_{i=1}^{4n} \ko_{\P^{2n+1}}(\zeta_{i}),$$

o\`u 
$$\zeta_{1}=\alpha(4n+1)+\beta,\zeta_{2}=4n\alpha+2\beta,\ldots,\zeta_{2n}=\alpha(2n+2)+2n\beta,\zeta_{2n+1}=2n\alpha+(2n+2)\beta,$$
$$\ldots,\zeta_{4n}=\alpha+(4n+1)\beta.$$ 

On a
$$ \sum_{i=1}^{4n}\zeta_{i}=c_{1}(\kn(-\gamma))+c_{1}(\kf(\gamma))=c_{1}(\kn) +c_{1}(\kf)$$
$$= n(2n+1)(\beta +\alpha)+3n(2n+1)(\beta +\alpha)=4n(2n+1)(\alpha+\beta ),$$

et
$$\zeta_{i}+\zeta_{4n+1-i}=2(2n+1)(\alpha+\beta ).$$

\subsubsection{\bf Proposition}\label{3.2.61} 
{\em Nous gardons les m\^emes notations de la proposition \ref{3.2.5}. Pour la suite exacte suivante
$$0\longrightarrow \kn(-\gamma)\stackrel{\pi}{\longrightarrow} \Upsilon_{1}  \stackrel{\iota}{\longrightarrow}  \kf(\gamma) \longrightarrow 0,$$ 

il existe deux morphismes $\phi$ et $\psi$ tels que le diagramme suivant est commutatif et exacte en $\Upsilon_{1}$
\xmat{ \Upsilon_{1}^{*} (-\hslash_{1})  \ar[rr]^{\phi}\ar[rd]_{}  &  & \Upsilon_{1} 
\ar[rr]^{\psi}\ar[rd]^{\iota}  & & \bigwedge ^{2n-1}\Upsilon_{1}^{*}(\hslash_{2})\\
&   \kn(-\gamma) \ar[ru]^{\pi}  \ar[rd]& &  \kf(\gamma)   \ar[rd]\ar[ru]^{}    \\
 0\ar[ru] &&   0\ar[ru]    && 0 ,}
  
o\`u $\hslash_{1}=2\gamma- (2n+1)(\alpha+\beta )$ et $\hslash_{2}=2n\gamma+3n(2n+1)(\alpha+\beta )$. De plus, si 
$$\gamma> max\lbrace \mid \alpha+\beta (4n+1)\mid + (2n+1)(\alpha+\beta ), \alpha(4n+1)+\beta\rbrace  ,$$

alors les deux morphismes $\phi$ et $\psi$ sont des matrices de rang $2n$ et $\hslash_{1}>0$.}

\begin{proof}
On consid\`ere 
$$ \BB_{1}=\lbrace e_{1}=d_{1},e_{2}=d_{2}, \ldots,e_{2n}=d_{2n},e_{2n+1}=d_{2n+2},e_{2n+2}=d_{2n+3},\ldots,e_{4n}=d_{4n+1} \rbrace $$

la base canonique de l'espace vectoriel $D_{1}$, et $ \BB_{1}{}^{*}=\lbrace e_{i}^{*} \rbrace_{1\leq i\leq 4n}$ la base canonique de l'espace vectoriel dual $D_{1}^{*}$. On consid\`ere aussi que 
$$ \bigwedge ^{2n-1}(\BB)^{*}=\lbrace e_{l_{1}}^{*}\wedge e_{l_{2}}^{*}\wedge \ldots  \wedge e_{l_{2n-1}}^{*} | 1\leq l_{1}<l_{2}< \ldots<l_{2n-1}\leq 4n  \rbrace, $$

est la base canonique de l'espace vectoriel $\bigwedge ^{2n-1}D_{1}^{*}$. D'apr\`es la proposition \ref{3.2.5} et que le fibr\'e $\kn$ est symplectique comme dans la proposition \ref{1.4.9}; le morphisme injectif $\pi: \kn(-\gamma)\longrightarrow \Upsilon_{1}$ donne
$$\Upsilon_{1}^{*}    \stackrel{^{T}\pi}{\longrightarrow} \kn^{*}(\gamma)\simeq \kn(\gamma-(2n+1)(\alpha+\beta )) \longrightarrow 0.$$

Alors on obtient le diagramme commutatif suivant
\xmat{ \Upsilon_{1}^{*}(-\hslash_{1})  \ar[rr]^{\phi}\ar[rd]_{^{T}\pi(-\hslash_{1})}  &  & \Upsilon_{1} \\
&   \kn(-\gamma) \ar[ru]_{\pi}  \ar[rd]& &     \\
 0\ar[ru] &&     0 ,}

o\`u $\hslash_{1}=2\gamma- (2n+1)(\alpha+\beta )$. D'apr\`es la proposition \ref{3.2.5}, le morphisme surjectif $\iota: \Upsilon_{1}\longrightarrow  \kf(\gamma)$ donne le morphisme injectif suivant 
$$\bigwedge ^{2n-1}{}^{T}\iota :\bigwedge ^{2n-1}(\kf(\gamma))^{*} \longrightarrow \bigwedge ^{2n-1}\Upsilon_{1}^{*}. $$

Comme $c_{1}(\kf(\gamma))=2n\gamma+3n(2n+1)(\alpha+\beta )=\hslash_{2}$ et $\bigwedge^{2n-1}((\kf(\gamma))^{*}) \otimes \bigwedge^{2n}(\kf(\gamma))=\kf(\gamma)$, alors on obtient le morphisme injectif
$$(\bigwedge ^{2n-1}{}^{T}\iota)(\hslash_{2}) : \kf(\gamma)\longrightarrow \left(\bigwedge ^{2n-1}\Upsilon_{1}^{*} \right)(\hslash_{2}) . $$

Donc on a le diagramme commutatif suivant
 \xmat{   \Upsilon_{1} \ar[rr]^{\psi}\ar[rd]_{\iota}  & & \left(\bigwedge ^{2n-1}\Upsilon_{1}^{*} \right)(\hslash_{2})\\
&     \kf(\gamma)  \ar[rd]\ar[ru]_{(\bigwedge ^{2n-1}{}^{T}\iota)(\hslash_{2})}    \\
 0\ar[ru] &&   0   .}

Comme on a
$$0\longrightarrow \kn(-\gamma)\stackrel{\pi}{\longrightarrow} \Upsilon_{1}  \stackrel{\iota}{\longrightarrow}  \kf(\gamma) \longrightarrow 0,$$
 
alors on obtient $\psi.\phi=0$ et que le diagramme suivant est commutatif 
\xmat{ \Upsilon_{1}^{*} (-\hslash_{1})  \ar[rr]^{\phi}\ar[rd]_{^{T}\pi(-\hslash_{1})}  &  & \Upsilon_{1} 
\ar[rr]^{\psi}\ar[rd]^{\iota}  & & \bigwedge ^{2n-1}\Upsilon_{1}^{*}(\hslash_{2})\\
&   \kn(-\gamma) \ar[ru]^{\pi}  \ar[rd]& &  \kf(\gamma)   \ar[rd]\ar[ru]_{(\bigwedge ^{2n-1}{}^{T}\iota)(\hslash_{2})}    \\
 0\ar[ru] &&   0\ar[ru]    && 0 .}

Si 
$$\gamma> max\lbrace \mid \alpha+\beta (4n+1)\mid + (2n+1)(\alpha+\beta ), \alpha(4n+1)+\beta\rbrace  ,$$

alors on peut d\'efinir le morphisme
$$\Upsilon_{1}^{*}    \stackrel{^{T}\pi}{\longrightarrow}\kn(\gamma-(2n+1)(\alpha+\beta )) \longrightarrow 0,$$

par une ligne $[S_{1},S_{2},\ldots ,S_{4n}]\otimes I_{\ko_{\P^{2n+1}}}$ o\`u $S_{i}\in H^{0}(\ko_{\P^{2n+1}}(\gamma+\zeta_{i}-(2n+1)(\alpha+\beta ))) $ sont des sections du fibr\'e $\kn(\gamma+\zeta_{i}-(2n+1)(\alpha+\beta ))$ pour $i=1,\ldots ,4n$. Dans ce cas, on peut d\'efinir $\phi$ par 
$$\phi= \left(
\begin{array}{ccccccccccccc}
0&S_{1,2}&\cdots&&S_{1,4n}&\\
-S_{1,2}&0&&&&\\
 \vdots&&\ddots&&\vdots& \\
&  &&0&S_{4n-1,4n}&\\
-S_{1,4n} &\cdots&&-S_{4n-1,4n}&0& \\
  \\ 
\end{array}
\right)=2 \sum_{1\leq i<j \leq 4n } S_{i,j} .e_{i}\wedge e_{j}$$ 

qui est une $(4n\times 4n)$-matrice antisym\'etrique de rang $2n$ d\'efinie relativement aux bases $\BB_{1}$ et $\BB_{1}{}^{*}$, et
$ S_{i,j}:=S_{i}\wedge S_{j}\in H^{0}(\ko_{\P^{2n+1}}(\zeta_{i}+\zeta_{j}+\hslash_{1}-(2n+1)(\alpha+\beta ))) $ pour $1\leq i<j\leq 4n $. On peut d\'efinir aussi le morphisme
$$\iota: \Upsilon_{1}\longrightarrow  \kf(\gamma)\longrightarrow 0,$$

par une ligne $[q_{1},q_{2},\ldots ,q_{4n}]\otimes I_{\ko_{\P^{2n+1}}}$ o\`u $q_{i}\in H^{0}(\ko_{\P^{2n+1}}(\gamma-\zeta_{i})) $ sont des sections du fibr\'e $\kf(\gamma-\zeta_{i})$ pour $i=1,\ldots ,4n$. Dans ce cas, on peut d\'efinir $\psi$ par

$$ \psi=\left(
\begin{array}{cccccccccc}
U_{1,1}&&\cdots&&U_{1,4n} \\
U_{2,1}&&\cdots&&U_{2,4n}  \\
 && &&            \\
 &&& &            \\
\vdots &&\cdots&& \vdots \\ 
 &&&&            \\
\vdots& &\cdots&&\vdots\\
&& &&           \\
 &&& &            \\
 U_{\left(\begin{array}{c}4n \\2n-1\end{array}\right) ,1}& &\cdots&& U_{\left(\begin{array}{c}4n \\2n-1\end{array}\right) ,4n}\\ 
\end{array}
\right)
$$ 

qui est une $(\left(\begin{array}{c}4n \\2n-1\end{array}\right)\times 4n)$-matrice de rang $2n$ d\'efinie relativement aux bases $\BB_{1}$ et $\bigwedge ^{2n-1}\BB_{1}{}^{*}$, et $ U_{l,j} \in H^{0}(\ko_{\P^{2n+1}}(2n\gamma-\zeta_{l_{1}}-\zeta_{l_{2}}-\ldots -\zeta_{l_{2n-1}}-\zeta_{j}))$ pour tout $1\leq l\leq \left(\begin{array}{c}4n \\2n-1\end{array}\right)$, $1\leq l_{1}< \ldots < l_{2n-1}\leq  4n$ et $1\leq j\leq  4n$. On a 
$$\psi= \sum_{1\leq l_{1}\neq \ldots \neq l_{2n-1}\leq 4n }^{1\leq j\leq 4n} q_{l_{1},\ldots , l_{2n-1},j} .e_{l_{1}}^{*}\wedge \ldots  \wedge e_{l_{2n-1}}^{*}\wedge e_{j}^{*},$$ 

o\`u $ U_{l,j}:= q_{l_{1},\ldots , l_{2n-1},j} :=q_{l_{1}}\wedge \ldots  \wedge q_{l_{2n-1}}\wedge q_{j}\in H^{0}(\ko_{\P^{2n+1}}(2n\gamma-\zeta_{l_{1}}-\zeta_{l_{2}}-\ldots -\zeta_{l_{2n-1}}-\zeta_{j})) $.  

Nous allons d\'ecrire les entr\'ees de la matrice $\psi$ par rapport aux entr\'ees de la matrice $ \phi$ et la construction de la matrice $\psi$ telle que $\psi.\phi=0$. Comme la matrice $\phi$ est antisym\'etrique, alors on a 
$$rg(\phi)=2n \text{ donne } Pf(\phi_{(l_{1},\ldots ,l_{2n-2})})=0,\hspace{0.2 cm} \text{ pour tout } 1\leq l_{1}<\ldots <l_{2n-2}\leq 4n, $$

o\`u $Pf(\phi_{(l_{1},\ldots ,l_{2n-2})})$ est le pfaffien de la $((2n+2)\times (2n+2)) $-sous-matrice antisym\'etrique $\phi_{(l_{1},\ldots ,l_{2n-2})} $ de $\phi$ qui est obtenue de $\phi$ par une \'elimination des colonnes et des lignes de num\'ero $ l_{1} , \ldots \text{ et } l_{2n-2}$ (simplement le $(2n+2)$-pfaffien de $\phi$). 

Pour tout $ i\in \lbrace 1,2, \ldots ,4n\rbrace \diagdown\lbrace  l_{1} , \ldots,  l_{2n-2}\rbrace$, on a
 
$$ Pf(\phi_{(l_{1},\ldots ,l_{2n-2})})= \sum_{j_{k}\in \lbrace 1,2, \ldots ,4n\rbrace \diagdown\lbrace i, l_{1} , \ldots,  l_{2n-2}\rbrace}^{1\leq k\leq 2n+1} (-1)^{k}S_{i,j_{k}}. Pf(\phi_{(l_{1},\ldots ,l_{2n-2},i,j_{k})})=0,$$

o\`u $Pf(\phi_{(l_{1},\ldots ,l_{2n-2},i,j_{k})})$ est le pfaffien de la $(2n\times 2n) $-sous-matrice antisym\'etrique $\phi_{(l_{1},\ldots ,l_{2n-2},i,j_{k})} $ de $\phi_{(l_{1},\ldots ,l_{2n-2})}$ (resp. de $\phi$) qui est obtenue de $\phi_{(l_{1},\ldots ,l_{2n-2})}$ (resp. de $\phi$) par une \'elimination des colonnes et des lignes de num\'ero $i,j_{k}$ (resp. $ l_{1} , \ldots,  l_{2n-2}, i \text{ et } j_{k}$). 
 
Par la d\'efinition de $\psi$, pour une ligne de la matrice $\psi$, il y a $2n-1$ entr\'ees qui sont nulles, et on peut \'ecrire cette matrice sous la forme suivante
$$ \psi=\left(
\begin{array}{cc}
\eta_{2n+2} \\
\eta_{2n+1}  \\
   \\
\vdots \\ 
   \\
 \eta_{4}\\ 
\end{array}
\right),
$$
 
o\`u chaque $\eta_{r}$, pour $4\leq r\leq 2n+2$, est une $(h_{r}\times4n)$-sous-matrice de $\psi$ telle que $$\sum_{r=4}^{2n+2}h_{r}= \left(\begin{array}{c}4n \\2n-1\end{array}\right),$$

et chaque $\eta_{r}$ contient des $(r\times4n)$-sous-matrices qui sont compos\'ees d'une $(r\times r)$-sous-matrice antisym\'etrique et de $2n-2$ colonnes nulles et de $2n+2-r$ colonnes. Soit $\AA(r,q)$ une $(r\times4n)$-sous-matrice de la $(h_{r}\times4n)$-matrice $\eta_{r}$ pour $4\leq r\leq 2n+2$ et pour $q\geq 1$
$$
\AA(r,q)=  \left(
\begin{array}{c|c|c|c|c|c|c|c|c|c|c|c}
0&0&\theta(r,q,m_{1})&0&\theta(r,q,m_{2})    &  \cdots &0& \XX(r,q)&0&\theta(r,q,m_{2n-4-r}) &\cdots &0\\
\end{array}
\right),
$$

o\`u $\XX(r,q)$ est une $(r\times r)$-sous-matrice antisym\'etrique dont les colonnes sont de num\'ero $1\leq w_{1}<w_{2}< \ldots<w_{r}\leq 4n$, et $\theta(r,q,m_{1}) ,\cdots , \theta(r,q,m_{2n+2-r})$ sont des colonnes de num\'ero $1\leq m_{1}<m_{2}< \ldots<m_{2n+2-r}\leq 4n$, et le reste des colonnes sont des $2n-2$ colonnes nulles de num\'ero $1\leq l_{1}<l_{2}< \ldots<l_{2n-2}\leq 4n$ avec $w_{1},w_{2}, \ldots,w_{r}$ et $m_{1},m_{2}, \ldots,m_{2n+2-r}$ et $l_{1},l_{2}, \ldots,l_{2n-2}$ sont des num\'eros diff\'erents. On consid\`ere $\phi_{(l_{1},\ldots ,l_{2n-2})}$ la $((2n+2)\times (2n+2)) $-sous-matrice antisym\'etrique de $\phi$ qui est obtenue de $\phi$ par une \'elimination des colonnes et des lignes de num\'ero $ l_{1} , \ldots \text{ et } l_{2n-2}$, et $Pf(\phi_{(l_{1},\ldots ,l_{2n-2})})$ est son pfaffien. Pour $w_{1}$, on a 
$$ Pf(\phi_{(l_{1},\ldots ,l_{2n-2})})= \sum_{j_{k}\in \lbrace 1,2, \ldots ,4n\rbrace \diagdown\lbrace w_{1}, l_{1} , \ldots,  l_{2n-2}\rbrace}^{1\leq k\leq 2n+1} (-1)^{k} S_{w_{1},j_{k}}. Pf(\phi_{(l_{1},\ldots ,l_{2n-2},w_{1},j_{k})})=0,$$

et on consid\`ere $(-1)^{k} Pf(\phi_{(l_{1},\ldots ,l_{2n-2},w_{1},j_{k})})$ comme les entr\'ees de la premi\`ere ligne de $\AA(r,q)$ \`a l'exception de $2n-1$ entr\'ees qui sont d\'ej\`a nulles, pour tout $1\leq k\leq 2n+1 $. Pour $w_{2}$, on a
$$ Pf(\phi_{(l_{1},\ldots ,l_{2n-2})})= \sum_{j_{k}\in \lbrace 1,2, \ldots ,4n\rbrace \diagdown\lbrace w_{2},  l_{1} , \ldots,  l_{2n-2}\rbrace}^{1\leq k\leq 2n+1} (-1)^{k} S_{w_{2},j_{k}}. Pf(\phi_{(l_{1},\ldots ,l_{2n-2},w_{2},j_{k})})=0,$$

et on consid\`ere $Pf(\phi_{(l_{1},\ldots ,l_{2n-2},w_{2},j_{1})})$ et $(-1)^{k} Pf(\phi_{(l_{1},\ldots ,l_{2n-2},w_{2},j_{k})})$ comme les entr\'ees de la deuxi\`eme ligne de $\AA(r,q)$ \`a l'exception de $2n-1$ entr\'ees qui sont d\'ej\`a nulles, pour tout $2\leq k\leq 2n+1 $. Pour $w_{3}$, on a
$$ Pf(\phi_{(l_{1},\ldots ,l_{2n-2})})= \sum_{j_{k}\in \lbrace 1,2, \ldots ,4n\rbrace \diagdown\lbrace w_{3},  l_{1} , \ldots,  l_{2n-2}\rbrace}^{1\leq k\leq 2n+1} (-1)^{k} S_{w_{3},j_{k}}. Pf(\phi_{(l_{1},\ldots ,l_{2n-2},w_{3},j_{k})})=0,$$

et on consid\`ere $Pf(\phi_{(l_{1},\ldots ,l_{2n-2},w_{3},j_{1})})$ et $-Pf(\phi_{(l_{1},\ldots ,l_{2n-2},w_{3},j_{2})})$ et $(-1)^{k} Pf(\phi_{(l_{1},\ldots ,l_{2n-2},w_{3},j_{k})})$ comme les entr\'ees de la deuxi\`eme ligne de $\AA(r,q)$ \`a l'exception de $2n-1$ entr\'ees qui sont d\'ej\`a nulles, pour tout $3\leq k\leq 2n+1 $. Pour $w_{4}$, on a
$$ Pf(\phi_{(l_{1},\ldots ,l_{2n-2})})= \sum_{j_{k}\in \lbrace 1,2, \ldots ,4n\rbrace \diagdown\lbrace w_{4},  l_{1} , \ldots,  l_{2n-2}\rbrace}^{1\leq k\leq 2n+1} (-1)^{k} S_{w_{4},j_{k}}. Pf(\phi_{(l_{1},\ldots ,l_{2n-2},w_{4},j_{k})})=0,$$

et on consid\`ere $Pf(\phi_{(l_{1},\ldots ,l_{2n-2},w_{4},j_{1})})$ et $-Pf(\phi_{(l_{1},\ldots ,l_{2n-2},w_{4},j_{2})})$ et $Pf(\phi_{(l_{1},\ldots ,l_{2n-2},w_{4},j_{2})})$ et $(-1)^{k} Pf(\phi_{(l_{1},\ldots ,l_{2n-2},w_{3},j_{k})})$ comme les entr\'ees de la deuxi\`eme ligne de $\AA(r,q)$ \`a l'exception de $2n-1$ entr\'ees qui sont d\'ej\`a nulles, pour tout $4\leq k\leq 2n+1 $. On continue comme \c{c}a jusqu'\`a $w_{r}$, on a
$$ Pf(\phi_{(l_{1},\ldots ,l_{2n-2})})= \sum_{j_{k}\in \lbrace 1,2, \ldots ,4n\rbrace \diagdown\lbrace w_{r}, l_{1} , \ldots,  l_{2n-2}\rbrace}^{1\leq k\leq 2n+1} (-1)^{k} S_{w_{r},j_{k}}. Pf(\phi_{(l_{1},\ldots ,l_{2n-2},w_{r},j_{k})})=0,$$

et on consid\`ere $(-1)^{k+1}Pf(\phi_{(l_{1},\ldots ,l_{2n-2},w_{r},j_{k})})$ comme les entr\'ees de la \mbox{$r^{\text{\`eme}}$-ligne} de $\AA(r,q)$ \`a l'exception de $2n-1$ entr\'ees qui sont d\'ej\`a nulles, pour tout $1\leq k\leq 2n+1 $. Donc on obtient que 
$$\AA(r,q).\varsigma_{\upsilon}(\phi)=0,$$

o\`u $\varsigma_{\upsilon}(\phi)$ est une colonne du num\'ero $\upsilon$ de $\phi$, pour tout $\upsilon \in \lbrace w_{1},w_{2}, \ldots,w_{r}; m_{1},m_{2}, \ldots,m_{2n+2-r} \rbrace $.

Pour la colonne $l_{1}$ de $\phi $, on a
$$ Pf(\phi_{(l_{2},\ldots ,l_{2n-2},w_{1})})= \sum_{j_{k}\in \lbrace 1,2, \ldots ,4n\rbrace \diagdown\lbrace w_{1}, l_{1},  l_{2} , \ldots,  l_{2n-2}\rbrace}^{1\leq k\leq 2n+1} (-1)^{k} S_{l_{1},j_{k}}. Pf(\phi_{(l_{1},\ldots ,l_{2n-2},w_{1},j_{k})})=0,$$

et
$$ Pf(\phi_{(l_{2},\ldots ,l_{2n-2},w_{2})})= \sum_{j_{k}\in \lbrace 1,2, \ldots ,4n\rbrace \diagdown\lbrace w_{2}, l_{1},  l_{2} , \ldots,  l_{2n-2}\rbrace}^{1\leq k\leq 2n+1} (-1)^{k} S_{l_{1},j_{k}}. Pf(\phi_{(l_{1},\ldots ,l_{2n-2},w_{2},j_{k})})=0.$$

On continue comme \c{c}a jusqu'\`a $w_{r} $
$$ Pf(\phi_{(l_{2},\ldots ,l_{2n-2},w_{r})})= \sum_{j_{k}\in \lbrace 1,2, \ldots ,4n\rbrace \diagdown\lbrace w_{r}, l_{1}, l_{2} , \ldots,  l_{2n-2}\rbrace}^{1\leq k\leq 2n+1} (-1)^{k} S_{l_{1},j_{k}}. Pf(\phi_{(l_{1},\ldots ,l_{2n-2},w_{r},j_{k})})=0.$$

On a la m\^eme chose pour toutes les colonnes $l_{i}$ de $\phi $ pour $2\leq i\leq 2n-2$. Donc on obtient que
$$\AA(r,q).\phi=0,$$

et une matrice $\psi$ dont ses entr\'ees sont les $2n$-pfaffiens de $\phi$ et telle que $\psi.\phi=0$.

\end{proof}

\subsubsection{\bf Remarque}\label{3.2.63} En gardant les m\^emes notations de la proposition \ref{3.2.61}, on a

- Dans la sous-matrice $\AA(r,q)$, chaque ligne de num\'ero $w_{i}$, $1\leq i\leq r$, contient toutes les entr\'ees de $\phi$ sauf les entr\'ees des lignes de num\'ero $w_{i}, l_{1} , \ldots \text{ et } l_{2n-2}$. Deux lignes quelconques de $\AA(r,q)$ contiennent toutes les entr\'ees de la matrice $\phi_{(l_{1},\ldots ,l_{2n-2})}$. 

- Dans la sous-matrice antisym\'etrique $\XX(r,q)$, chaque colonne de num\'ero $w_{i}$, $1\leq i\leq r$, contient toutes les entr\'ees de $\phi$ sauf les entr\'ees des colonnes de num\'ero $w_{i}, l_{1} , \ldots \text{ et } l_{2n-2}$. Deux colonnes quelconques de $\XX(r,q)$ contiennent toutes les entr\'ees de la matrice $\phi_{(l_{1},\ldots ,l_{2n-2})}$. Chaque colonne de $\psi$ est compos\'e des colonnes de toutes les matrices $\AA(r,q)$: soit une colonne nulle ou bien une colonne de la sous-matrice antisym\'etrique $\XX(r,q)$ ou bien une colonne comme $\theta(r,q)$. Donc chaque colonne de la matrice $\psi$ contient toutes les entr\'ees de la matrice $\phi$.
 
- D'apr\`es la d\'efinition de $\psi$ et la base de l'espace vectoriel $\bigwedge ^{2n-1}D_{1}^{*}$
$$ \bigwedge ^{2n-1}(\BB)^{*}=\lbrace e_{l_{1}}^{*}\wedge e_{l_{2}}^{*}\wedge \ldots  \wedge e_{l_{2n-1}}^{*} | 1\leq l_{1}<l_{2}< \ldots<l_{2n-1}\leq 4n  \rbrace,$$

on obtient qu'il existe une ligne de $\psi$ telle que: 

Soient $1\leq l_{1}<l_{2}< \ldots<l_{2n-1}\leq 4n$ les num\'eros des entr\'ees nulles dans cette ligne. Tous les $2n-2$ num\'eros des entr\'ees nulles de cette ligne v\'erifient $4n+1-l_{i}=l_{l}$ avec $ 1\leq i\neq l\leq 2n-1$. 
Pour $1\leq a\leq 4n$ un num\'ero d'une entr\'ee nulle, on a $\left(\begin{array}{c}2n \\n-1\end{array}\right)-1$ lignes comme la ligne pr\'ec\'edente.

- Comme $rg(\phi)=2n$, alors les $a$-mineurs de la matrice $\phi$ n'ont pas un z\'ero commun sur $\P^{2n+1}$ o\`u $1\leq a\leq 2n$, en particulier les entr\'ees $S_{i,j}$ n'ont pas un z\'ero commun sur $\P^{2n+1}$. 

- On donne ici un exemple de la sous-matrice $\eta_{r}$
$$ \eta_{2n+2}=\left(
\begin{array}{cc}
\AA(2n+2,1) \\
\AA(2n+2,2)  \\
  \\
\vdots \\ 
   \\
 \AA(2n+2,n)\\ 
\end{array}
\right),
$$
 
o\`u $\AA(2n+2,q)$ est une $((2n+2)\times4n)$-sous-matrice de la $(n(2n+2)\times4n)$-matrice $\eta_{2n+2}$. On a, pour $1 \leq q\leq n$,
$$
\AA(2n+2,q)= \overbrace{[ 0|0|0|0|0|  \cdots |0  }^{\mbox{$2(n-q)$ colonnes}}|  \XX(2n+2,q)| \overbrace{0|0|0|0   ]}^{\mbox{$2(q-1)$ colonnes}} $$

o\`u $\XX(2n+2,q)$ est une $((2n+2)\times (2n+2))$-sous-matrice antisym\'etrique dont ses entr\'ees sont d\'efinies comme dans \ref{3.2.61}. 

\subsubsection{\bf Remarque}\label{3.2.6} 
Soient $z \in H^{0}(\P^{2n+2},\ko_{\P^{2n+2}}(1))\setminus \lbrace 0 \rbrace $ et $\P^{2n+1}   \subset  \P^{2n+2} $ l'hyperplan de $\P^{2n+2}$ d\'efini par l'\'equation $z=0$. Soient $\epsilon\geq 1$ un entier et $Y_{(\epsilon)}\subset \P^{2n+2}$ le voisinage infinit\'esimal de l'ordre $\epsilon$ de $\P^{2n+1} $ dans $\P^{2n+2}$, d\'efini par l'\'equation $z^{\epsilon}=0$. On a les morphismes d'inclusion suivants

 \xmat{  \P^{2n+1} :=Y_{(1)}  \ar@{^{(}->}[rr]^{b_{(2)}} &  & Y_{(2)} \ar@{^{(}->}[rr]^{b_{(3)}} &  &  \ldots \ar@{^{(}->}[rr]^{b_{(\epsilon)}} &  &   Y_{(\epsilon)}\ar@{^{(}->}[rr]^{d} &  &  \P^{2n+2}  .\\
 }

Soit $e_{0}\in \P^{2n+2}\setminus Y_{(\epsilon)}$, on d\'efinit la projection d'un point (voir \cite{ha} page $22$) sur un sous-espace projectif $\P^{2n+1}$ par
$$\pi:\P^{2n+2}\setminus \{ e_{0}\}\longrightarrow \P^{2n+1}$$
$$ e\longmapsto \pi(e)$$

o\`u $\pi(e)$ est l'intersection de $\P^{2n+1}$ avec la droite unique passant par les points $e_{0}$ et $ e$. On consid\`ere la restriction de la projection $\pi$ sur le voisinage infinit\'esimal $Y_{(\epsilon)}$ et on obtient la projection d'un point sur un sous-espace projectif $\P^{2n+1}$
$$J_{(\epsilon)}:Y_{(\epsilon)}\longrightarrow \P^{2n+1}.$$

On consid\`ere $\kn_{(\epsilon)}=J_{(\epsilon)}^{*}\kn$ et $\kf_{(\epsilon)}=J_{(\epsilon)}^{*}\kf$ o\`u $ \kn$ est le fibr\'e de 0-corr\'elation pond\'er\'e provenant d'une image inverse g\'en\'eralis\'ee sur $\P^{2n+1}$, et $ \kf$ est le fibr\'e de Tango pond\'er\'e provenant d'une image inverse g\'en\'eralis\'ee sur $\P^{2n+1}$ comme dans la proposition \ref{3.2.5}. En prenant l'image inverse de la suite exacte suivant
$$0\longrightarrow \kn(-\gamma)\stackrel{\pi}{\longrightarrow}   \Upsilon_{1}  \stackrel{\iota}{\longrightarrow}  \kf(\gamma) \longrightarrow 0,$$ 

on obtient donc la suite exacte suivante

\begin{equation}
0\longrightarrow \kn_{(\epsilon)}(-\gamma)\stackrel{\pi}{\longrightarrow}  \Upsilon_{(\epsilon)}  \stackrel{\iota}{\longrightarrow}  \kf_{(\epsilon)}(\gamma)   \longrightarrow 0.
\end{equation}

o\`u 
$$ \Upsilon_{(\epsilon)}:=J_{(\epsilon)}^{*}\Upsilon_{1}= \bigoplus_{i=1}^{4n} \ko_{Y_{(\epsilon)} }(\zeta_{i}).$$ 

On a aussi les deux suites exactes suivantes sur $Y_{(\epsilon)}$
$$0\longrightarrow \kf_{(\epsilon-1)}(-1)\longrightarrow \kf_{(\epsilon)}\longrightarrow \kf_{(1)}\longrightarrow 0,$$

et
$$0\longrightarrow \kn_{(\epsilon-1)}(-1)\longrightarrow \kn_{(\epsilon)}\longrightarrow \kn_{(1)}\longrightarrow 0.$$

\subsubsection{\bf Th\'eor\`eme}\label{3.2.7}
{\em  Nous gardons les m\^emes notations de la proposition \ref{3.2.61}. Soit $n$ un entier positif tel que $n>1$. Soient $\kn_{(\epsilon)}$ et $\kf_{(\epsilon)}$ comme dans la remarque \ref{3.2.6} et $\kg$ un fibr\'e vectoriel sur $\P^{2n+2}$ d\'efini par le diagramme commutatif suivant

\xmat{ &0\ar[d]&0\ar[d]&&\\
&\Upsilon(-\epsilon)\ar@{=}[r] \ar[d]& \Upsilon(-\epsilon)\ar[d]&&\\
0 \ar[r]&\kg \ar[r]\ar[d]&\Upsilon \ar[r]\ar[d]& \kf_{(\epsilon)}(\gamma) \ar[r]\ar@{=}[d]&0\\
0 \ar[r]& \kn_{(\epsilon)}(-\gamma) \ar[r]^{\pi}\ar[d]&\Upsilon_{(\epsilon)} \ar[r]^{\iota}\ar[d]& \kf_{(\epsilon)}(\gamma) \ar[r]&0\\
&0&0&&}

o\`u 
$$ \Upsilon:= \bigoplus_{j=1}^{4n} \ko_{\P^{2n+2}}(\zeta_{j}).$$ 

Soit $\Gamma:=\bigoplus_{i=1}^{2n-1} \ko_{\P^{2n+2}}(\zeta_{b_{i}})$ avec $1\leq b_{1}<b_{2}< \ldots<b_{2n-1}\leq 4n$. On a 

1 - Si $\epsilon=\epsilon_{1}:=2\gamma+(2n+1)(\alpha+\beta )$, alors il existe un fibr\'e vectoriel $\ke$ de rang $2n+1$ sur $\P^{2n+2}$ qui est d\'efini par la suite exacte suivante
$$0 \longrightarrow \Gamma(-\epsilon_{1}) \longrightarrow \kg \longrightarrow \ke \longrightarrow 0. $$
2 - Si $\epsilon=\epsilon_{2}:=2n\gamma+n(2n+1)(\alpha+\beta )$, alors il existe un fibr\'e vectoriel $\kk$ de rang $2n+1$ sur $\P^{2n+2}$ qui est d\'efini par la suite exacte suivante
$$0 \longrightarrow \kk \longrightarrow \kg \longrightarrow \Gamma\longrightarrow 0. $$

Parmi ces fibr\'es Il y a des fibr\'es qui sont diff\'erents du fibr\'e de Tango pond\'er\'e provenant d'une image inverse g\'en\'eralis\'ee sur $\P^{2n+2}$ \cite{bah1}.}
 
\begin{proof}
De la proposition \ref{3.2.61}, on a le diagramme commutatif sur $\P^{2n+1}$ suivant
\xmat{  \Upsilon_{1}^{*} (-\hslash_{1})  \ar[rr]^{\phi}\ar[rd]  &  & \Upsilon_{1} 
\ar[rr]^{\psi}\ar[rd]  & & \LL_{1}(\hslash_{2})\\
&   \kn(-\gamma) \ar[ru]  \ar[rd]& &  \kf(\gamma)   \ar[rd]\ar[ru]    \\
 0\ar[ru] &&   0\ar[ru]    && 0 ,}

o\`u $ \LL_{1}:=\left(\bigwedge ^{2n-1}\Upsilon_{1}^{*}\right)$. En prenant l'image inverse du diagramme pr\'ec\'edent par la projection d'un point sur un sous-espace projectif $\P^{2n+1}$, la remarque \ref{3.2.6} ,
$$J_{(\epsilon)}:Y_{(\epsilon)}\longrightarrow \P^{2n+1} ,$$
 
on obtient le diagramme suivant sur $Y_{(\epsilon)}$
\xmat{  \Upsilon^{*}_{(\epsilon)} (-\hslash_{1})  \ar[rr]^{\phi}\ar[rd] &  & \Upsilon_{(\epsilon)} 
\ar[rr]^{\psi}\ar[rd]  & & \kl_{(\epsilon)}(\hslash_{2})\\
&   \kn_{(\epsilon)}(-\gamma) \ar[ru]  \ar[rd] & &  \kf_{(\epsilon)}(\gamma)   \ar[rd]\ar[ru]   \\
 0\ar[ru] &&   0\ar[ru]    && 0 } 

o\`u $  \Upsilon_{(\epsilon)}:=\bigoplus_{i=1}^{4n} \ko_{Y_{(\epsilon)}}(\zeta_{i})$. En choisissant les m\^emes notations de \ref{3.1} avec $\Phi=\phi$ et $\Psi=\psi$ sur $\P^{2n+2}$, tous les fibr\'es $\Upsilon^{*}_{(\epsilon)} (-\hslash_{1})$, $ \Upsilon_{(\epsilon)} $ et $ \kl_{(\epsilon)}(\hslash_{2})$ peuvent se relever en fibr\'es $\Upsilon^{*} (-\hslash_{1})$, $ \Upsilon$ et $ \kl(\hslash_{2})$ respectivement sur $\P^{2n+2}$. On obtient donc le complexe suivant sur $\P^{2n+2}$
$$\Upsilon^{*} (-\hslash_{1}) \stackrel{ \Phi }{\longrightarrow} \Upsilon \stackrel{ \Psi }{\longrightarrow} \kl(\hslash_{2}), $$

o\`u  $\Psi .\Phi  =\psi.\phi=0$. Alors on peut d\'efinir le morphisme suivant sur $\P^{2n+2}$ 

$$\Delta:\Upsilon(-\epsilon)\oplus \Upsilon^{*}(-\hslash_{1})\longrightarrow \Upsilon \oplus \kl(\hslash_{2}-\epsilon) $$
$$
{\Delta} = \left(
\begin{array}{c|c}
z^{\epsilon}I    &    \Phi \\
 \hline
\Psi  &    0 \\
\end{array}
\right)
,$$ 

qui est la $((4n+\left(\begin{array}{c}4n \\2n-1\end{array}\right) )\times 8n)$-matrice suivante

$$
 \left(
\begin{array}{ccccc|ccccccc}
z^{\epsilon}&&&&        &  0&S_{1,2}&\cdots&\cdots&S_{1,4n}\\
&z^{\epsilon}&&&\bigzero&  - S_{1,2}&0&&&\vdots&\\
&&\ddots&&&           \vdots&&\ddots&&\vdots& \\
\bigzero&&&z^{\epsilon}&&     \vdots&&&0&S_{4n-1,4n}&\\
&&&&z^{\epsilon}&            -S_{1,4n} &\cdots&\cdots&-S_{4n-1,4n}&0 \\
  \hline
U_{1,1}&&\cdots&&U_{1,4n}&            &&&&\\
U_{2,1}&&\cdots&&U_{2,4n}&            &&&&\\
\vdots &&\cdots&& \vdots&            && \bigzero && \\ 
U_{\left(\begin{array}{c}4n \\2n-1\end{array}\right) ,1}& &\cdots&& U_{\left(\begin{array}{c}4n \\2n-1\end{array}\right) ,4n}&     &&&& \\ 
\end{array}
\right).
$$ 
 
Soient 
$$\Gamma:=\bigoplus_{i=1}^{2n-1} \ko_{\P^{2n+2}}(\zeta_{b_{i}}) \text{ avec } 1\leq b_{1}<b_{2}< \ldots<b_{2n-1}\leq 4n.$$ 

1- On consid\`ere $\overline{\mu_{b_{i}}}$ la somme de la $b_{i}^{\text{\`eme}}$-colonne de la matrice $\Delta$ avec la $(8n+1-b_{i})^{\text{\`eme}}$-colonne de la matrice $\Delta$. Plus pr\'ecis\'ement, soit

$$(h_{1}, \ldots  ,h_{4n},h_{4n+1}, \ldots ,h_{8n})\in \left( \bigoplus_{j=1}^{4n}H^{0}(\ko_{\P^{2n+2}}( \zeta_{j}-\epsilon))\right)\bigoplus \left( \bigoplus_{j=1}^{4n} H^{0}(\ko_{\P^{2n+2}}(-\zeta_{j}-\hslash_{1}))\right). $$

- On consid\`ere $\overline{\mu_{b_{i}} }$ la somme la $b_{i}^{\text{\`eme}}$-colonne et la $(8n+1-b_{i})^{\text{\`eme}}$-colonne de la matrice $\Delta$ avec la condition d'homog\'en\'eit\'e suivante, sur la $b_{i}^{\text{\`eme}}$-coordonn\'ee et la $(8n+1-b_{i})^{\text{\`eme}}$-coordonn\'ee de $(h_{1}, \ldots  ,h_{4n},h_{4n+1}, \ldots ,h_{8n})$, 
$$deg(h_{b_{i}})=deg(h_{8n+1-b_{i}})$$ 
$$-\epsilon+\zeta_{b_{i}}=-\zeta_{4n+1-b_{i}}-\hslash_{1},$$

o\`u $1\leq i\leq 2n-1$. Donc, pour tout $1\leq i\leq 2n-1$, de telles sommes existent gr\^ace \`a la condition $\epsilon=\epsilon_{1}:=2\gamma+(2n+1)(\alpha+\beta )$ et au fait que $\zeta_{j}=-\zeta_{4n+1-j} +2(2n+1)(\alpha+\beta )$ pour tout $1\leq j\leq 4n$. Donc on obtient la matrice
$$
\overline{\mu } = \left(
\begin{array}{ccccccccccc}
\overline{\mu_{b_{1}}}&,&\overline{\mu_{b_{2}}}&,&\cdots& &,& \overline{\mu_{b_{2n-1}}}
\end{array}
\right),
$$
$$\overline{\mu}:\Gamma(-\epsilon_{1})\longrightarrow  \Upsilon \oplus\kl( \hslash_{2}-\epsilon_{1}). $$

qui est une $((4n+\left(\begin{array}{c}4n \\2n-1\end{array}\right) )\times (2n-1))$-matrice. 

D'apr\`es la proposition \ref{3.2.61} et la remarque \ref{3.2.63}, on a: $p$ un mineur maximal de $\overline{\mu}$ est une forme homog\`ene qui est une somme des termes de la forme
$$ \prod_{i,j}S_{i,j}^{a_{i,j}}. z^{ \epsilon_{1} b}: \text{ $a_{i,j},b \in \N $ tels que $\epsilon_{1} b +\sum_{i,j} a_{i,j}deg (S_{i,j})=deg(p)$}.$$ 

Donc le mineur maximal $p$ est une somme des termes de la forme
$$ X^{r_{0}}_{0}.X^{r_{1}}_{1}. \ldots .X^{r_{2n+1}}_{2n+1}. z^{ \epsilon_{1} b}  : \text{ $r_{0},\ldots,r_{2n+1},b \in \N $ tels que  $ r_{0}+\ldots +r_{2n+1}+\epsilon_{1} b = deg(p)$}.$$

Soient $\overline{\Xi}$ l'id\'eal des mineurs maximaux de $\overline{\mu}$ qui est engendr\'e par tous les mineurs maximaux de $\overline{\mu}$, et $\km=(X_{0},X_{1},\ldots,X_{2n+1},z)$ un id\'eal maximal sur $\P^{2n+2}$, et $\overline{m}:= max\lbrace deg(p): p \in \overline{\Xi}\rbrace$. Alors $\km^{\overline{m}} \subset \overline{\Xi} $ o\`u $\km^{\overline{m}}$ est un exposant de $\km$. Donc $\overline{\Xi}$ est \mbox{$\km$-primaire} sur $\P^{2n+2}$, d'apr\`es la d\'efinition de l'id\'eal primaire. On a aussi $\overline{\Xi}\neq 0$ sur $\P^{2n+2}$. Donc l'id\'eal annulateur de $\overline{\Xi}$ est $Ann(\overline{\Xi})=0$ sur $\P^{2n+2}$. Ce qui entra\^ine que le morphisme $\overline{\mu}$ est un morphisme injectif non nul sur $\P^{2n+2}$. D'apr\`es la proposition \ref{3.1.2} et la remarque \ref{3.1.3}, on a le diagramme commutatif suivant
    
\xmat{ & \Gamma(-\epsilon_{1} ) \ar@{^{(}->}[ld]_{\overline{\aleph}} \ar@{-->}[dd]^{}\ar[rd]^{\overline{\mu}}\\
\Upsilon(-\epsilon_{1})\oplus \Upsilon^{*}(- \hslash_{1})  \ar[rr]_{\Delta}\ar[rd]^{} & & \Upsilon \oplus \kl( \hslash_{2}-\epsilon_{1} )\\
&\kg\ar[ur]^{}\ar[dr]^{}\\
 0\ar[ur] & & 0
 } 
  
o\`u $\overline{\aleph}:\Gamma(-\epsilon_{1} )\hookrightarrow \Upsilon(-\epsilon_{1} )\oplus \Upsilon^{*}(-\hslash_{1})$ est un morphisme d'inclusion, alors le morphisme suivant est repr\'esent\'e par le morphisme $\overline{\mu}$ 
$$\Gamma(-\epsilon_{1})\longrightarrow \kg.$$

Alors on obtient la suite exacte suivante
$$0 \longrightarrow \Gamma(-\epsilon_{1}) \stackrel{ \overline{\mu} }{\longrightarrow} \kg \longrightarrow \ke \longrightarrow 0, $$

o\`u $\ke$ est un fibr\'e vectoriel de rang $2n+1$ sur $\P^{2n+2}$. On utilise la proposition \ref{3.1.8} pour calculer les classes de Chern du faisceau $\kf_{(\epsilon_{1})}(\gamma)$. En utilisant la suite exacte pr\'ec\'edente et la suite exacte suivante
$$0\longrightarrow \kg \longrightarrow \Upsilon \longrightarrow \kf_{(\epsilon_{1})}(\gamma) \longrightarrow 0$$

on obtient que la premi\`ere classe de Chern du fibr\'e $\kg$ est $c_{1}(\kg)=-4n\gamma$, et la premi\`ere classe de Chern du fibr\'e $\ke$ est   
$$c_{1}(\ke)=(4n^{2}-1)(\alpha+\beta)-2\gamma -\sum_{i=1}^{2n-1}\zeta_{b_{i}}.$$

2- Maintenant, pour le fibr\'e $\kk$. D'apr\`es la remarque \ref{3.2.63}, il existe une ligne de la matrice $\Psi$ telle que: 

Soient $1\leq l_{1}<l_{2}< \ldots<l_{2n-1}\leq 4n$ les num\'eros des entr\'ees nulles dans cette ligne. Tous les $2n-2$ num\'eros des entr\'ees nulles de cette ligne v\'erifient $4n+1-l_{r}=l_{l}$ avec $ 1\leq r\neq a\leq 2n-1$.

Pour une ligne de la matrice $\Delta$ de num\'ero exactement plus grand que $4n$, on va appeler la partie de cette ligne dans la sous-matrice $\Psi$ une sous-ligne de cette ligne dans la sous-matrice $\Psi$ (simplement sa sous-ligne dans $\Psi$).

- On choisit une ligne de la matrice $\Delta$ telle que sa sous-ligne dans $\Psi$ poss\`ede un z\'ero de num\'ero $4n+1-b_{1}$. Soient $1\leq l_{1}(1)<l_{2}(1)< \ldots<l_{2n-2}(1)\leq 4n$ les num\'eros des entr\'ees nulles de cette sous-ligne dans $\Psi$, tels qu'ils v\'erifient $4n+1-l_{r}(1)=l_{a}(1)$ avec $ 1\leq r\neq a\leq 2n-2$. Donc cette ligne contient la forme $z^{\epsilon}$ et toutes les entr\'ees de $\Phi$ sauf les entr\'ees des colonnes de num\'ero $4n+1-b_{1}, l_{1}(1) , \ldots \text{ et } l_{2n-2}(1)$. 

On consid\`ere $\mu_{b_{1}}$ la somme de cette ligne avec la $b_{1}^{\text{\`eme}}$-ligne de la matrice $\Delta$ avec la condition d'homog\'en\'eit\'e suivante 
$$\zeta_{b_{1}}=-\zeta_{4n+1-b_{1}}-2(n-1)(2n+1)(\alpha+\beta )-\epsilon+\hslash_{2}.$$

Une telle somme existe gr\^ace \`a la condition $\epsilon=\epsilon_{2}:=2n\gamma+n(2n+1)(\alpha+\beta )$ et au fait que $\zeta_{j}=-\zeta_{4n+1-j} +2(2n+1)(\alpha+\beta )$ pour tout $1\leq j\leq 4n$.

- On choisit maintenant une autre ligne de la matrice $\Delta$ telle que sa sous-ligne dans $\Psi$ poss\`ede un z\'ero de num\'ero $4n+1-b_{2}$. Soient $1\leq l_{1}(2)<l_{2}(2)< \ldots<l_{2n-2}(2)\leq 4n$ les num\'eros des entr\'ees nulles de cette sous-ligne dans $\Psi$, tels qu'ils v\'erifient $4n+1-l_{r}(2)=l_{l}(2)$ avec $ 1\leq r\neq a\leq 2n-2$ et $ l_{1}(2),l_{2}(2), \ldots,l_{2n-2}(2)$ sont diff\'erents des num\'eros $l_{1}(1),l_{2}(1), \ldots,l_{2n-2}(1), 4n+1-b_{1}$. Donc cette ligne contient la forme $z^{\epsilon}$ et toutes les entr\'ees de $\Phi$ sauf les entr\'ees des colonnes de num\'ero $4n+1-b_{2}, l_{1}(2) , \ldots \text{ et } l_{2n-2}(2)$. 

On consid\`ere $\mu_{b_{2}}$ la somme de cette ligne avec la $b_{2}^{\text{\`eme}}$-ligne de la matrice $\Delta$ avec la condition d'homog\'en\'eit\'e suivante 
$$\zeta_{b_{2}}=-\zeta_{4n+1-b_{2}}-2(n-1)(2n+1)(\alpha+\beta )-\epsilon+\hslash_{2}.$$

Une telle somme existe gr\^ace \`a la condition $\epsilon=\epsilon_{2}:=2n\gamma+n(2n+1)(\alpha+\beta )$ et au fait que $\zeta_{j}=-\zeta_{4n+1-j} +2(2n+1)(\alpha+\beta )$ pour tout $1\leq j\leq 4n$. 

- On choisit maintenant une autre ligne de la matrice $\Delta$ telle que sa sous-ligne dans $\Psi$ poss\`ede un z\'ero de num\'ero $4n+1-b_{i}$ pour $ 3\leq i\leq 2n-1$. Soient $1\leq m_{1}(i)<m_{2}(i)< \ldots<m_{2n-2}(i)\leq 4n$ les num\'eros des entr\'ees nulles de cette sous-ligne dans $\Psi$, tels qu'ils v\'erifient $4n+1-m_{r}(i)=m_{l}(i)$ avec $ 1\leq r\neq a\leq 2n-2$. Donc cette ligne contient la forme $z^{\epsilon}$ et toutes les entr\'ees de $\Phi$ sauf les entr\'ees des colonnes de num\'ero $4n+1-b_{i}, m_{1}(i) , \ldots \text{ et } m_{2n-2}(i)$. 

On consid\`ere $\mu_{b_{i}}$ la somme de cette ligne avec la $b_{i}^{\text{\`eme}}$-ligne de la matrice $\Delta$ avec la condition d'homog\'en\'eit\'e suivante 
$$\zeta_{b_{i}}=-\zeta_{4n+1-b_{i}}-2(n-1)(2n+1)(\alpha+\beta )-\epsilon+\hslash_{2}.$$

Une telle somme existe gr\^ace \`a la condition $\epsilon=\epsilon_{2}:=2n\gamma+n(2n+1)(\alpha+\beta )$ et au fait que $\zeta_{j}=-\zeta_{4n+1-j} +2(2n+1)(\alpha+\beta )$ pour tout $1\leq j\leq 4n$. Donc on obtient la $(2n-1)\times(8n)$-matrice suivante
$$ \mu=\left(
\begin{array}{cc}
\mu_{b_{1}} \\
\mu_{b_{2}} \\
  \\
\vdots \\ 
  \\
 \mu_{b_{2n-1}}\\ 
\end{array}
\right),
$$
$$\mu:\Upsilon(-\epsilon_{2})\oplus \Upsilon^{*}(-\hslash_{1}) \longrightarrow \Gamma, $$

qui existe gr\^ace \`a la condition $\epsilon=\epsilon_{2}:=2n\gamma+n(2n+1)(\alpha+\beta )$ et au fait que $$\zeta_{j}=-\zeta_{4n+1-j} +2(2n+1)(\alpha+\beta )\text{ pour tout }  1\leq j\leq 4n.$$

On obtient donc que la matrice $\mu$ contient la forme $z^{\epsilon_{2}}$ et toutes les entr\'ees de $\Phi$ car tous ces num\'eros $ l_{1}(2),l_{2}(2), \ldots,l_{2n-2}(2), l_{1}(1),l_{2}(1), \ldots,l_{2n-2}(1), 4n+1-b_{1}, 4n+1-b_{2}$ sont diff\'erents.
D'apr\`es la proposition \ref{3.2.61} et la remarque \ref{3.2.63}, on a: $p$ un mineur maximal de $\mu$ est une forme homog\`ene qui est une somme des termes de la forme
$$ \prod_{i,j}S_{i,j}^{a_{i,j}}. z^{\epsilon_{2} b}: \text{ $a_{i,j},b \in \N $ tels que $ \epsilon_{2} b +\sum_{i,j} a_{i,j}deg (S_{i,j})=deg(p)$}.$$ 

Donc le mineur maximal $p$ est une somme des termes de la forme
$$ X^{r_{0}}_{0}.X^{r_{1}}_{1}. \ldots .X^{r_{2n+1}}_{2n+1}. z^{ \epsilon_{2} b}  : \text{ $r_{0},\ldots,r_{2n+1},b \in \N $ tels que  $ r_{0}+\ldots +r_{2n+1}+ \epsilon_{2} b = deg(p) $}.$$

Soient $\Xi$ l'id\'eal des mineurs maximaux de $\mu$ est engendr\'e par tous les mineurs maximaux de $\mu$, et $\km=(X_{0},X_{1},\ldots,X_{2n+1},z)$ un id\'eal maximal sur $\P^{2n+2}$, et $m:= max\lbrace deg(p): p \in \Xi\rbrace$. Alors $\km^{m} \subset \Xi$ o\`u $\km^{m}$ est un exposant de $\km$. Donc $\Xi$ est \mbox{$\km$-primaire} sur $\P^{2n+2}$, d'apr\`es la d\'efinition de l'id\'eal primaire. On a aussi $\Xi\neq 0$ sur $\P^{2n+2}$. Donc l'id\'eal annulateur de $\Xi$ est $Ann(\Xi)=0$ sur $\P^{2n+2}$. Ce qui entra\^ine que le morphisme $\mu$ est un morphisme surjectif non nul sur $\P^{2n+2}$. D'apr\`es la proposition \ref{3.1.4}, on a le diagramme commutatif suivant

\xmat{& & \Gamma \\
\\
\Upsilon(-\epsilon_{2})\oplus \Upsilon^{*}(-\hslash_{1}) \ar[rrrr]_{\Delta}\ar[rrdd]^{}\ar[rruu]^{\mu} & & & & \Upsilon \oplus \kl(2n(2n+1)(\alpha+\beta )) \ar@{->>}[lluu]_{ \aleph}\\
\\
& & \kg\ar[uurr]^{}\ar[dr]^{} \ar@{-->}[uuuu]_{}\\
& 0\ar[ur] & & 0} 
 
o\`u $ \aleph:  \Upsilon \oplus \kl(2n(2n+1)(\alpha+\beta )) \longrightarrow \Gamma$ est un morphisme de projection surjectif, alors le morphisme suivant est repr\'esent\'e par le morphisme $\mu$ 
$$ \kg\longrightarrow \Gamma.$$

Alors on obtient la suite exacte suivante
$$0 \longrightarrow  \kk \longrightarrow \kg \stackrel{\mu}{\longrightarrow} \Gamma\longrightarrow 0, $$

o\`u $\kk$ est un fibr\'e vectoriel de rang $2n+1$ sur $\P^{2n+2}$. On utilise la proposition \ref{3.1.8} pour calculer les classes de Chern du faisceau $\kf_{(\epsilon_{2})}(\gamma)$. En utilisant la suite exacte pr\'ec\'edente et la suite exacte suivante
$$0\longrightarrow \kg \longrightarrow \Upsilon \longrightarrow \kf_{(\epsilon_{2})}(\gamma) \longrightarrow 0$$

on obtient que la premi\`ere classe de Chern du fibr\'e $\kg$ est $c_{1}(\kg)=-4n\gamma$, et la premi\`ere classe de Chern du fibr\'e $\kk$ est   
$$c_{1}(\kk)= -4n\gamma-\sum_{i=1}^{2n-1}\zeta_{b_{i}}.$$

\end{proof}

\subsubsection{\bf Proposition}\label{3.2.8}
{\em  Nous gardons les m\^emes notations de la proposition \ref{3.2.61} et du th\'eor\`eme \ref{3.2.7}. Soient $\kn_{(\epsilon)}$ et $\kf_{(\epsilon)}$ comme dans la remarque \ref{3.2.6} et $\kg$ un fibr\'e vectoriel sur $\P^{4}$ d\'efini comme dans le th\'eor\`eme \ref{3.2.7}. On a :
 
- Si $ \epsilon=\epsilon_{3}:=2\gamma+7\alpha -\beta $, alors il existe un fibr\'e vectoriel $\ke_{1}$ (resp. $\kk_{4}$) de rang $3$ sur $\P^{4}$ tel que
$$0 \longrightarrow \ko_{\P^{4}}(\zeta_{1}-\epsilon_{3}) \longrightarrow \kg \longrightarrow \ke_{1} \longrightarrow 0. $$

(resp.
$$0 \longrightarrow \kk_{4} \longrightarrow \kg \longrightarrow\ko_{\P^{4}}(\zeta_{4})  \longrightarrow 0). $$

- Si $ \epsilon=\epsilon_{4}:=2\gamma+5\alpha +\beta $, alors il existe un fibr\'e vectoriel $\ke_{2}$ (resp. $\kk_{3}$) de rang $3$ sur $\P^{4}$ tel que
$$0 \longrightarrow \ko_{\P^{4}}(\zeta_{2}-\epsilon_{4}) \longrightarrow \kg \longrightarrow \ke_{2} \longrightarrow 0. $$

(resp.
$$0 \longrightarrow \kk_{3} \longrightarrow \kg \longrightarrow\ko_{\P^{4}}(\zeta_{3})  \longrightarrow 0). $$

- Si $ \epsilon=\epsilon_{5}:=2\gamma+\alpha +5\beta $, alors il existe un fibr\'e vectoriel $\ke_{3}$ (resp. $\kk_{2}$) de rang $3$ sur $\P^{4}$ tel que
$$0 \longrightarrow \ko_{\P^{4}}(\zeta_{3}-\epsilon_{5}) \longrightarrow \kg \longrightarrow \ke_{3} \longrightarrow 0. $$

(resp.
$$0 \longrightarrow \kk_{2} \longrightarrow \kg \longrightarrow\ko_{\P^{4}}(\zeta_{2})  \longrightarrow 0). $$

- Si $ \epsilon=\epsilon_{6}:=2\gamma-\alpha +7\beta $, alors il existe un fibr\'e vectoriel $\ke_{4}$ (resp. $\kk_{1}$) de rang $3$ sur $\P^{4}$ tel que
$$0 \longrightarrow \ko_{\P^{4}}(\zeta_{4}-\epsilon_{6}) \longrightarrow \kg \longrightarrow \ke_{4} \longrightarrow 0, $$

(resp.
$$0 \longrightarrow \kk_{1} \longrightarrow \kg \longrightarrow\ko_{\P^{4}}(\zeta_{1})  \longrightarrow 0). $$

Parmi ces fibr\'es Il y a des fibr\'es qui sont diff\'erents du fibr\'e de Tango pond\'er\'e provenant d'une image inverse g\'en\'eralis\'ee sur $\P^{4}$ \cite{bah1}.}

\begin{proof} Dans la proposition \ref{3.2.61} Pour $n=1$, on a  
$$ 
\phi = \left(
\begin{array}{cccc}
0 & -S_{1,2}&- S_{1,3} &-S_{1,4} \\
S_{1,2} & 0 &  -S_{2,3} &-S_{2,4} \\
S_{1,3} & S_{2,3} &0  &-S_{3,4} \\
S_{1,4} & S_{2,4} & S_{3,4} & 0 \\
\end{array}
\right)
,
$$

et 
$$ 
\psi = \left(
\begin{array}{cccc}
0 & -S_{3,4}& S_{2,4} &-S_{3,2} \\
S_{3,4} & 0 &  -S_{1,4} &S_{1,3} \\
-S_{2,4} & S_{1,4} &0  &-S_{1,2} \\
S_{3,2} & -S_{1,3} & S_{1,2} & 0 \\
\end{array}
\right)
.
$$

et le fibr\'e $\kf$ (resp. $\kn$) est de rang $2$ si et seulement si $pf(\phi)=S_{1,2}S_{3,4}-S_{1,3}S_{2,4}+S_{1,4}S_{2,3}=0$. On fait les m\^emes choses comme dans le th\'eor\`eme \ref{3.2.7} avec les deux matrices pr\'ec\'edentes pour obtenir le morphisme sur $\P^{4}$
$$\Delta:\Upsilon(-\epsilon)\oplus \Upsilon^{*}(-\hslash_{1})\longrightarrow \Upsilon \oplus \Upsilon^{*}(\hslash_{2}-\epsilon) $$

$$
{\Delta} = \left(
\begin{array}{c|c}
z^{\epsilon}I    &    \Phi \\
 \hline
\Psi  &    0 \\
\end{array}
\right)
,$$

o\`u $\hslash_{1}=2\gamma- 3(\alpha+\beta )>0$ et $\hslash_{2}=2\gamma+9(\alpha+\beta )>0$. On consid\`ere la condition $\epsilon=\epsilon_{3}:=2\gamma+7\alpha -\beta$. Soit 
$$(h_{1},h_{2}, h_{3} ,h_{4},h_{5},h_{6}, h_{7} ,h_{8})\in H^{0}(\ko_{\P^{4}}(\zeta_{1}-\epsilon))\oplus H^{0}(\ko_{\P^{4}}(\zeta_{2}-\epsilon))\oplus H^{0}(\ko_{\P^{4}}(\zeta_{3}-\epsilon))$$
$$\oplus H^{0}(\ko_{\P^{4}}(\zeta_{4}-\epsilon))\oplus H^{0}(\ko_{\P^{4}}(-\zeta_{1}-\hslash_{1}))\oplus H^{0}(\ko_{\P^{4}}(-\zeta_{2}-\hslash_{1}))$$
$$\oplus H^{0}(\ko_{\P^{4}}(-\zeta_{3}-\hslash_{1}))\oplus H^{0}(\ko_{\P^{4}}(-\zeta_{4}-\hslash_{1})). $$

On prend la somme la $1^{\text{\`ere }}$-colonne et la $5^{\text{\`eme}}$-colonne de la matrice $\Delta$ avec la condition d'homog\'en\'eit\'e suivante, sur la $1^{\text{\`eme}}$-coordonn\'ee et la $5^{\text{\`eme}}$-coordonn\'ee de \\
$(h_{1},h_{2}, h_{3} ,h_{4},h_{5},h_{6}, h_{7} ,h_{8})$,
$$deg(h_{1})=deg(h_{5})$$ 
$$-k+\zeta_{1}=-\zeta_{1}-\hslash_{1}.$$ 

Une telle somme existe gr\^ace \`a la condition $\epsilon=\epsilon_{3}:=2\gamma+7\alpha -\beta$. Donc on obtient une colonne $\mu_{1}$
$$\mu_{1}=^{T}[z^{\epsilon_{3}}, S_{1,2},S_{1,3} ,S_{1,4};0,S_{3,4},-S_{2,4} ,S_{3,2}]$$

$$\mu_{1}:(\ko_{\P^{4}}(\zeta_{1}))(-\epsilon_{3})\longrightarrow  \Upsilon \oplus \Upsilon^{*}(\hslash_{2}-\epsilon_{3}). $$ 

Supposons que toutes les formes de la matrice $\mu_{1}$ soient nulles sur $ \P^{4}$. On obtient donc que la matrice $\Phi$ est nulle, ce qui est une contradiction au fait que $rg (\Phi)=2$. Donc toutes les formes de la matrice $\mu_{1}$ n'ont pas un z\'ero commun sur $\P^{4}$. Ce qui donne que le morphisme $\mu_{1}$ est injectif non nul. Mais, le morphisme suivant est repr\'esent\'e par le morphisme $\mu_{1}$
$$\ko_{\P^{4}}(\zeta_{1}-\epsilon_{3})\longrightarrow \kg.$$ 

Alors on obtient la suite exacte suivante
$$0 \longrightarrow \ko_{\P^{4}}(\zeta_{1}-\epsilon_{3}) \stackrel{ \mu_{1} }{\longrightarrow} \kg \longrightarrow \ke_{1} \longrightarrow 0,$$

o\`u $\ke_{1}$ un fibr\'e vectoriel de rang $3$ sur $\P^{4}$. On utilise la proposition \ref{3.1.8} pour calculer les classes de Chern du faisceau $\kf_{(\epsilon_{3})}(\gamma)$. En utilisant la suite exacte pr\'ec\'edente et la suite exacte suivante
$$0\longrightarrow \kg \longrightarrow \Upsilon \longrightarrow \kf_{(\epsilon_{3})}(\gamma) \longrightarrow 0$$

on obtient les classes de Chern du fibr\'e $\ke_{1}$ 
$$c_{1}(\ke_{1})=-2\gamma+12\beta.$$

- Pour le fibr\'e $ \kk_{4}$. Soit 
$$(a_{1},a_{2}, a_{3} ,a_{4},a_{5},a_{6}, a_{7} ,a_{8})\in H^{0}(\ko_{\P^{4}}(\zeta_{1}))\oplus H^{0}(\ko_{\P^{4}}(\zeta_{2}))\oplus H^{0}(\ko_{\P^{4}}(\zeta_{3}))$$
$$\oplus H^{0}(\ko_{\P^{4}}(\zeta_{4}))\oplus H^{0}(\ko_{\P^{4}}(\hslash_{2}-\zeta_{1}-\epsilon))\oplus H^{0}(\ko_{\P^{4}}(\hslash_{2}-\zeta_{2} -\epsilon))$$
$$\oplus H^{0}(\ko_{\P^{4}}(\hslash_{2}-\zeta_{3} -\epsilon))\oplus H^{0}(\ko_{\P^{4}}(\hslash_{2}-\zeta_{4} -\epsilon )). $$

On prend la somme la $4^{\text{\`eme}}$-ligne et la $8^{\text{\`eme}}$-ligne de la matrice $\Delta$ avec la condition d'homog\'en\'eit\'e suivante, sur la $4^{\text{\`eme}}$-coordonn\'ee et la $8^{\text{\`eme}}$-coordonn\'ee de $(a_{1},a_{2}, a_{3} ,a_{4},a_{5},a_{6}, a_{7} ,a_{8})$, 
$$deg(a_{4})=deg(a_{8})$$ 
$$\zeta_{4}=\hslash_{2}-\zeta_{4} -\epsilon.$$ 

Une telle somme existe gr\^ace \`a la condition $\epsilon=\epsilon_{3}:=2\gamma+7\alpha -\beta$. Donc on obtient une ligne $\delta$ 

$$\delta=[ S_{3,2},-S_{1,3} ,S_{1,2},z^{\epsilon_{3}},S_{1,4}, S_{2,4} ,S_{3,4},0]$$

$$\delta: \Upsilon(-\epsilon_{3}) \oplus \Upsilon^{*}(-\hslash_{1}) \longrightarrow  \ko_{\P^{4}}(\zeta_{4}) . $$

Supposons que toutes les formes de la matrice $\delta$ soient nulles sur $ \P^{4}$. On obtient donc que la matrice $\Phi$ est nulle, ce qui est une contradiction au fait que $rg (\Phi)=2$. Donc toutes les formes de la matrice $\delta$ n'ont pas un z\'ero commun sur $\P^{4}$. Ce qui donne que le morphisme $\delta $ est surjectif non nul. Mais, le morphisme suivant est repr\'esent\'e par le morphisme $\delta $
$$ \kg  \longrightarrow  \ko_{\P^{4}}(\zeta_{4}).$$ 

Alors on obtient la suite exacte suivante
$$0 \longrightarrow \kk_{4} \longrightarrow \kg  \stackrel{ \delta  }{\longrightarrow}  \ko_{\P^{4}}(\zeta_{4}) \longrightarrow 0, $$

o\`u $\kk_{4}$ est un fibr\'e vectoriel de rang $3$ sur $\P^{4}$. On utilise la proposition \ref{3.1.8} pour calculer les classes de Chern du faisceau $\kf_{(\epsilon_{3})}(\gamma)$. En utilisant la suite exacte pr\'ec\'edente et la suite exacte suivante
$$0\longrightarrow \kg \longrightarrow \Upsilon \longrightarrow \kf_{(\epsilon_{3})}(\gamma) \longrightarrow 0$$

on obtient les classes de Chern du fibr\'e $\kk_{4}$
$$ c_{1}(\kk_{4})=3(3\beta-\alpha)-4\gamma.$$

On fait la m\^eme chose pour le fibr\'e $\ke_{2}$ (resp. $\ke_{3}$, $\ke_{4}$) en consid\'erant la condition $\epsilon=\epsilon_{4}:=2\gamma+5\alpha +\beta$ (resp. $\epsilon=\epsilon_{5}:=2\gamma+\alpha +5\beta$, $\epsilon=\epsilon_{6}:=2\gamma-\alpha +7\beta$ ). On obtient une colonne qui est la somme de la $2^{\text{\`eme}}$ (resp. $3^{\text{\`eme}}$, $4^{\text{\`eme}}$)-colonne et la $6^{\text{\`eme}}$ (resp.  $7^{\text{\`eme}}$, $8^{\text{\`eme}}$)-colonne de la matrice $\Delta$ sur $\P^{4}$ telle que les formes de cette colonne n'ont pas un z\'ero commun sur $\P^{4}$.

On fait la m\^eme chose aussi pour le fibr\'e $\kk_{3}$ (resp. $\kk_{2}$, $\kk_{1}$) en consid\'erant la condition $\epsilon=\epsilon_{4}:=2\gamma+5\alpha +\beta$ (resp. $\epsilon=\epsilon_{5}:=2\gamma+\alpha +5\beta$, $\epsilon=\epsilon_{6}:=2\gamma-\alpha +7\beta$ ). On obtient une ligne qui est la somme de la $3^{\text{\`eme}}$ (resp. $2^{\text{\`eme}}$, $1^{\text{\`ere}}$)-ligne et la $7^{\text{\`eme}}$ (resp. $6^{\text{\`eme}}$, $5^{\text{\`eme}}$)-ligne de la matrice $\Delta$ sur $\P^{4}$ telle que les formes de cette ligne n'ont pas un z\'ero commun sur $\P^{4}$. Dans ces cas-l\`a, on obtient ce que l'on recherche.

\end{proof}

\section{ \large\bf Exemple de fibr\'e vectoriel de rang $3$ sur $\P^{4}$ }
\vspace{1cm}

Nous construisons un fibr\'e vectoriel de rang $3$ sur l'espace projectif $\P^{4}:=\P_{\K}^{4}$ tout en utilisant l'identit\'e suivante de Binet-Cauchy et la m\'ethode de Kumar-Peterson-Rao \cite{ku-ra-pe}, o\`u $\K$ est un corps quelconque.

Soient $c_{i},x_{i},y_{i}$ et $z_{i}$ des \'el\'ements dans un anneau commutatif et $m_{0}>1$ un entier. On a l'identit\'e suivante
$$\sum_{1\leq i<j\leq m_{0}} 
(x_{i}z_{j}-x_{j}z_{i}).(c_{i}y_{j}-c_{j}y_{i})=\hspace{10cm} $$
$$=\sum_{1\leq i<j\leq m_{0}}(x_{i}z_{j}c_{i}y_{j}+x_{j}z_{i}c_{j}y_{i})- \sum_{1\leq i<j\leq m_{0}}(x_{i}z_{j}c_{j}y_{i}+x_{j}z_{i}c_{i}y_{j})$$
$$ =\sum_{1\leq i<j\leq m_{0}}(x_{i}z_{j}c_{i}y_{j}+x_{j}z_{i}c_{j}y_{i})- \sum_{1\leq i<j\leq m_{0}}(x_{i}z_{j}c_{j}y_{i}+x_{j}z_{i}c_{i}y_{j})$$
$$+ \sum_{1\leq i\leq m_{0}}x_{i}z_{i}c_{i}y_{i}-\sum_{1\leq i\leq m_{0}}x_{i}z_{i}c_{i}y_{i}  $$
$$=(\sum_{i=1}^{m_{0}}x_{i}c_{i}).(\sum_{j=1}^{m_{0}}y_{j}z_{j})-(\sum_{i=1}^{m_{0}}
x_{i}y_{i}).(\sum_{j=1}^{m_{0}}z_{j}c_{j}).$$

Cette identit\'e s'appelle {\em l'identit\'e de Binet-Cauchy}. En particulier pour $m_{0}=2$, on a l'identit\'e suivante

\begin{equation}\label{18} 
(x_{1}c_{1}+x_{2}c_{2}).(y_{1}z_{1}+y_{2}z_{2})=
(x_{1}y_{1}+x_{2}y_{2}).(c_{1}z_{1}+c_{2}z_{2})+
(x_{1}z_{2}-x_{2}z_{1}).(c_{1}y_{2}-c_{2}y_{1})
\end{equation} 

\subsection{ Fibr\'e vectoriel de rang $2$ sur $\P^{3}$}\label{4.1}
Dans cette partie, nous allons construire et \'etudier une famille de fibr\'es 
vectoriels de rang $2$ sur $\P^{3}:=\P_{\K}^{3} $ o\`u $\K$ est un corps quelconque, tout en utilisant les matrices 
antisym\'etriques des formes sur $\P^{3} $.

\subsubsection{Morphisme associ\'e \`a une matrice antisym\'etrique sur  
$\P^{3} $}\label{4.1.1} 

Soient $(S_{i;j})_{1\leq i \neq j\leq 4}$ des polyn\^omes homog\`enes de $4$ 
variables (des formes sur  $\P^3 $). On consid\`ere une matrice $4\times 4$ 
antisym\'etrique
$$ 
{M} = \left(
\begin{array}{cccc}
0 & S_{1,2}& S_{1,3} &S_{1,4} \\
S_{2,1} & 0 &  S_{2,3} &S_{2,4} \\
S_{3,1} & S_{3,2} &0  &S_{3,4} \\
S_{4,1} & S_{4,2} & S_{4,3} & 0 \\
\end{array}
\right)
,
$$

avec $S_{i,j}=-S_{j,i}$ et $d_{ij}=deg(S_{i,j})=deg(-S_{j,i})=d_{ji}$. On d\'efinit le pfaffien de la matrice antisym\'etrique $M$ par $(pfM)^{2}=det(M)$ qui est le polyn\^ome homog\`ene (\'equation du Pl\"ucker pour la grassmannienne $G(2,4)$)
$$Q=pf(M)=S_{1,2}S_{3,4}-S_{1,3}S_{2,4}+S_{1,4}S_{2,3}.$$

Il a un degr\'e  
$$deg(Q)=d=d_{12}+d_{34}=d_{13}+d_{24}=d_{14}+d_{23}.$$  

On suppose que les formes $(S_{i,j})_{1\leq i \neq j\leq 4}$ n'ont pas de z\'ero commun sur $\P^3 $. Donc en calculant les mineurs $2\times 2 $ et $3\times 3$ de la matrice $M$, on constate que $M$ est de rang $2$ en tout point de $\P^3 $ si et seulement si le polyn\^ome $Q=0$. 

On d\'efinit la matrice associ\'ee \`a une matrice antisym\'etrique $M$. En utilisant le fait que $Q=0$, on obtient la matrice antisym\'etrique suivante
$$ 
{N} = \left(
\begin{array}{cccc}
0 & S_{3,4}& -S_{2,4} &S_{3,2} \\
-S_{3,4} & 0 &  S_{1,4} &-S_{1,3} \\
S_{2,4} & -S_{1,4} &0  &S_{1,2} \\
-S_{3,2} & S_{1,3} & -S_{1,2} & 0 \\
\end{array}
\right)
.
$$

On obtient que la matrice $N$ est de rang $2$, et $M.N=N.M=0$. On consid\`ere $\psi_{0}$ le morphisme de fibr\'es suivant associ\'e \`a la matrice $M$
$$\psi_{0} :\F :=\bigoplus_{i=1}^{4}\ko_{\P^3}(a_{i})  \longrightarrow 
 \bigoplus_{i=1}^{4}\ko_{\P^3}(b_{i}).$$
 
On peut supposer que $a_{1}\leq a_{2} \leq a_{3} \leq a_{4}$. Donc on peut trouver facilement 
les relations suivantes 
$$ d_{ij} = b_{i}-a_{j}  \hspace{1cm} si   \hspace{0.2cm} 1\leq i\neq j\leq 4.$$ 

Soit $\Lambda_{ij}=d_{1j}- d_{ij}$ pour $1<j\leq 4 $ et $ 1\leq i\leq 4 $, donc on a
  $$a_{i}=a_{1}  +\Lambda_{ij} \hspace{0.2cm} ,  \hspace{0.2cm}b_{i}=b_{1}  -\Lambda_{ij}.$$

On a aussi $ b_{n}= a_{1}+d_{in}+\Lambda_{ij}$, si n est un entier tel que $i\neq n$ et $ 1\leq n\leq 4 $. On a
$$ b_{1}= a_{1}+d_{12}+d_{13}+d_{14}-d=a_{1}-d_{34}-d_{23}-d_{24}+2d.$$

Donc on a 
$$\F =\bigoplus_{i=1}^{4}\ko_{\P^3}(a_{i})=\ko_{\P^3}(a_{1})\oplus\ko_{\P^3}(a_{1}+\Lambda_{2j})\oplus\ko_{\P^3}(a_{1}+\Lambda_{3j})\oplus\ko_{\P^3}(a_{1}+\Lambda_{4j}).$$

Donc le morphisme $\psi_{0} $ devient 

\begin{equation}\label{19}  
\psi_{0} :\F  \longrightarrow \F^{*}(b_{1}+a_{1}).
\end{equation}

On suppose, dans la suite de notre \'etude, que $Q=0$. On consid\`ere $\varphi_{0}$ le morphisme suivant de fibr\'es associ\'e \`a la matrice $N$
$$\varphi_{0} :\bigoplus_{i=1}^{4}\ko_{\P^3}(q_{i})  \longrightarrow 
 \F =\bigoplus_{i=1}^{4}\ko_{\P^3}(a_{i}).$$

On obtient alors les relations suivantes 
$$q_{m}= a_{i}-d_{lr},$$

o\`u $i,r,l,m$ sont des entiers distincts tels que $1\leq i,r,l,m\leq 4$. Donc on a
$$q_{m}= b_{m}-d_{im}-d_{lr}=b_{1}-d-\Lambda_{mj},$$

o\`u $j$ est un entier tel que $1\leq j\neq m\leq 4$. Alors le morphisme $\varphi_{0}$ devient

\begin{equation} \label{20} 
\varphi_{0} :\F^{*}(b_{1}+a_{1}-d)\longrightarrow \F.  
\end{equation}

Des morphismes \ref{19} et \ref{20}, on obtient un complexe
$$ \F^{*}(b_{1}+a_{1}-d)\stackrel{\varphi_{0} }{\longrightarrow}  \F  
\stackrel{ \psi_{0} }{\longrightarrow}  \F^{*}(b_{1}+a_{1}) .$$

Pour $e \in \Z$ entier, on consid\`ere les morphismes de fibr\'es
$$\varphi_{e}: \F^{*}(b_{1}+a_{1}+(e-1)d) \longrightarrow \F(ed) ,$$
  
$$\psi_{e}: \F(ed)\longrightarrow \F^{*}(b_{1}+a_{1}+ed). $$ 

Donc on obtient un complexe associ\'e infini de deux c\^ot\'es,
$$\ldots \stackrel{\psi_{-3} }{\longrightarrow} 
\F^{*}(b_{1}+a_{1}-3d)\stackrel{\varphi_{-2} }{\longrightarrow} 
\F(-2d)\stackrel{\psi_{-2} 
}{\longrightarrow}\F^{*}(b_{1}+a_{1}-2d)\stackrel{\varphi_{-1} 
}{\longrightarrow} \F(-d) \stackrel{\psi_{-1} }{\longrightarrow} 
\F^{*}(b_{1}+a_{1}-d)$$

$$\stackrel{\varphi_{0} }{\longrightarrow}\F \stackrel{\psi_{0} 
}{\longrightarrow}\F^{*}(b_{1}+a_{1})\stackrel{\varphi_{1} }{\longrightarrow} 
\F(d)\stackrel{\psi_{1} }{\longrightarrow} 
\F^{*}(b_{1}+a_{1}+d)\stackrel{\varphi_{2} }{\longrightarrow} 
\F(2d)\stackrel{\psi_{2} 
}{\longrightarrow}\F^{*}(b_{1}+a_{1}+2d)\stackrel{\varphi_{3} 
}{\longrightarrow} \ldots$$

\subsubsection{\bf D\'efinition du fibr\'e de rang $2$ sur $\P^3 $ }\label{4.1.2} 

{\em Soient les morphismes $\psi_{0} $ et $ \varphi_{0}$ \ref{19} et \ref{20} dans \ref{4.1.1}. On obtient que  $\A=Im\varphi_{0}$ et $\B=Im\psi_{0}$ sont des fibr\'es vectoriels de rang $2$ sur $\P^3 $}. On a le complexe suivant

\xmat{ \ldots  \F^{*}(b_{1}+a_{1}-d)  \ar[rr]^{\varphi_{0}}\ar[rd] &  & \F  
\ar[rr]^{\psi_{0}}\ar[rd]  & & \F^{*}(b_{1}+a_{1}) \ar[rd]    \ar[rr]^{\varphi_{1}} &&\F(d)\ldots    \\
&   \A \ar[ru]  \ar[rd] & &  \B   \ar[rd]\ar[ru] && \A(d) \ar[ru] \ar[rd] \\
 0\ar[ru] &&   0\ar[ru]    && 0\ar[ru] && 0.}
 
 \vspace{0.8cm}

\subsubsection{\bf Proposition }\label{4.1.3} 

{\em  Soient  $\A$ et $\B$ les fibr\'es vectoriels sur $\P^{3}$ qui sont d\'efinis dans \ref{4.1.2}.
Les classes de Chern de fibr\'es vectoriels $\A$ et $\B$ sont
$$c_{1}(\B)=b_{1}+a_{1}=2a_{1}+d_{13}-d_{23}-d_{12}. $$

$$c_{2}(\B)=-\frac{1}{d}\lbrace a_{1}a_{2}a_{3}+a_{1}a_{3}a_{4}+a_{1}a_{2}a_{4}+a_{2}a_{3}a_{4}\hspace{6cm}$$ 

$$-(b_{1}+ a_{1})(a_{1}a_{2}+a_{1}a_{3}+a_{1}a_{4}+a_{2}a_{3}+ a_{2}a_{4}+a_{3}a_{4}-(b_{1}+a_{1})(b_{1}+a_{1}-d)) \rbrace .$$

$$ c_{1}(\A)=b_{1}+a_{1}-d=2a_{1}+d_{13}-d_{23}-d_{34}.$$

$$c_{2}(\A)= \frac{1}{d}\lbrace a_{1}a_{2}a_{3}+a_{1}a_{3}a_{4}+a_{1}a_{2}a_{4}+a_{2}a_{3}a_{4}\hspace{6cm}$$ 

$$-(b_{1}+ a_{1}-d)(a_{1}a_{2}+a_{1}a_{3}+a_{1}a_{4}+a_{2}a_{3}+a_{2}a_{4}+a_{3}a_{4}-(b_{1}+a_{1})(b_{1}+a_{1}-d)) \rbrace .$$ }

\begin{proof}

Comme $\A$, $\B$ sont des fibr\'es vectoriels de rang $2$, alors on a un morphisme surjectif de fibr\'es
$$ \F^{*}(c_{1}(\A))\longrightarrow \A,$$

et un morphisme injectif de fibr\'es
$$ \B \longrightarrow \F^{*}(c_{1}(\B)).$$

Donc on a que

\xmat{   \F^{*}(c_{1}(\A))  \ar[rr]^{\rho}\ar[rd] &  & \F  
\ar[rr]^{\varrho}\ar[rd]  & & \F^{*}(c_{1}(\B))    \\
&   \A \ar[ru]  \ar[rd] & &  \B   \ar[rd]\ar[ru] && \\
 0\ar[ru] &&   0\ar[ru]    && 0.}

et $\rho= \varphi_{0}$ , $\varrho = \psi_{0}$. Alors on a 
$$c_{1}(\B)=b_{1}+a_{1}=2a_{1}+d_{13}-d_{23}-d_{12}, $$

$$ c_{1}(\A)=b_{1}+a_{1}-d=2a_{1}+d_{13}-d_{23}-d_{34}.$$

On calcule le reste des classes de Chern directement tout en utilisant la suite exacte suivante
$$ 0 \longrightarrow \A\longrightarrow\F\longrightarrow \B\longrightarrow 0.$$

\end{proof}

 \subsubsection{\bf Proposition }\label{4.1.5} 
{\em Soient $\lbrace T_{i}\rbrace ,   \hspace{0.2cm} \lbrace W_{i}\rbrace ,   \hspace{0.2cm} \lbrace V_{i}\rbrace  ,   \hspace{0.2cm} \lbrace U_{i} \rbrace $, $1\leq i\leq 2$ des formes homog\`enes sur $\P^3 $ de degr\'es positifs $ deg(T_{i})=t_{i} ,   \hspace{0.2cm} deg(W_{i})=w_{i} ,   \hspace{0.2cm} deg(V_{i})=v_{i} ,   \hspace{0.2cm} deg(U_{i})=u_{i}$. On consid\`ere les formes homog\`enes sur $\P^3 $ suivantes 
$$S_{1,3}=T_{1}.W_{1}+T_{2}.W_{2}   ,   \hspace{0.2cm} S_{2,4}=U_{1}.V_{1}+U_{2}.V_{2}$$
 
$$ S_{1,2}= T_{1}.U_{1}+T_{2}.U_{2} ,   \hspace{0.2cm} S_{3,4}= V_{1}.W_{1}+V_{2}.W_{2}$$

$$S_{1,4}= T_{1}.V_{2}-T_{2}.V_{1}   ,   \hspace{0.2cm} S_{2,3}= U_{2}.W_{1}-U_{1}.W_{2}$$

telles que $(S_{i,j})_{1\leq i \neq j\leq 4}$ n'ont pas de z\'ero commun sur $\P^3 $. On consid\`ere les conditions suivantes
$$t_{2}=w_{2},\hspace{0.2cm} u_{2}=v_{2},\hspace{0.2cm} t_{1}=-w_{1}+2w_{2},\hspace{0.2cm} u_{1}=w_{1}-w_{2}+v_{2},\hspace{0.2cm} v_{1}=-w_{1}+w_{2}+v_{2},$$

avec $w_{2},v_{2}>0$ et $  0\leq w_{2}-v_{2}\leq w_{1}\leq  w_{2}+v_{2}$. On consid\`ere $M$ et $N$ des  $4\times 4$-matrices antisym\'etriques suivantes
 
 $$ {M} = \left(
\begin{array}{cccc}
0 & S_{1,2}& S_{1,3} &S_{1,4} \\
S_{2,1} & 0 &  S_{2,3} &S_{2,4} \\
S_{3,1} & S_{3,2} &0  &S_{3,4} \\
S_{4,1} & S_{4,2} & S_{4,3} & 0 \\
\end{array}
\right) $$

et

$${N} = \left(
\begin{array}{cccc}
0 & S_{3,4}& -S_{2,4} &S_{3,2} \\
-S_{3,4} & 0 &  S_{1,4} &-S_{1,3} \\
S_{2,4} & -S_{1,4} &0  &S_{1,2} \\
-S_{3,2} & S_{1,3} & -S_{1,2} & 0 \\
\end{array}
\right)$$

On consid\`ere aussi les morphismes $\varphi_{0} $ et $\psi_{0} $ \ref{19} et \ref{20} associ\'es respectivement aux matrices antisym\'etriques $N$ et $M$ dans le diagramme suivant   } 

\xmat{  \F^{*}(b_{1}+a_{1}-d)  \ar[rr]^{\varphi_{0}}\ar[rd] &  & \F  
\ar[rr]^{\psi_{0}}\ar[rd]  & & \F^{*}(b_{1}+a_{1}) \ar[rd]    \ar[rr]^{\varphi_{1}} &&\F(d)    \\
&   \ka \ar[ru]  \ar[rd] & &  \kb   \ar[rd]\ar[ru] && \ka(d) \ar[ru] \ar[rd] \\
 0\ar[ru] &&   0\ar[ru]    && 0\ar[ru] && 0}

\vspace{0.4cm}
 
{\em Alors on a 

I) $M.N=N.M=0$ et $S_{1,2}.S_{3,4}-S_{1,3}.S_{2,4}+S_{1,4}S_{2,3}=0$ et $b_{1}= 3w_{2}-w_{1}+ a_{1}$ et $d=2(w_{2}+v_{2})$.
 
II) On a le fibr\'e $\F  = \ko_{\P^3}(a_{1})\oplus \ko_{\P^3}(a_{1}+2w_{2}-w_{1}-v_{2})\oplus \ko_{\P^3}(a_{1}+w_{2}-w_{1})\oplus \ko_{\P^3}(a_{1}+w_{2}-v_{2}) $, pour tout $a_{1}\in \Z$.

III) On a deux fibr\'es vectoriels $\ka=Im\varphi_{0}$ et $\kb =Im\psi_{0}$ de rang $2$ sur $\P^3 $ 
avec des classes de Chern, pour $a_{1}=0$,
$$c_{1}(\kb )= 3 w_{2}- w_{1} , \hspace{0.2cm} c_{2}(\kb )=w_{2}( 2w_{2}-w_{1}+v_{2}), $$ 
 
$$ c_{1}(\ka)= w_{2}-w_{1}-2v_{2},  \hspace{0.2cm} c_{2}(\ka)=v_{2}(w_{1}+v_{2}).$$ }

\begin{proof}
Il suffit de prendre dans l'identit\'e  de Binet-Cauchy \ref{18}
$$x_{1}=T_{1},\hspace{0.2cm} x_{2}=T_{2},\hspace{0.2cm} y_{1}=U_{1},\hspace{0.2cm} y_{2}=U_{2},\hspace{0.2cm} c_{1}=W_{1},\hspace{0.2cm}  c_{2}=W_{2} ,\hspace{0.2cm} z_{1}=V_{1},\hspace{0.2cm} z_{2}=V_{2}, $$

pour obtenir
 $$S_{1,2}.S_{3,4}-S_{1,3}.S_{2,4}+S_{1,4}S_{2,3}=0.$$

En utilisant \ref{4.1.1} et \ref{4.1.3} , on obtient (II) et (III).
 
\end{proof}

\subsection{Fibr\'e vectoriel de rang $3$ sur $\P^{4}$}

Dans cette partie, nous allons construire un exemple de fibr\'es vectoriels de 
rang $3$ sur $\P^{4}:=\P_{\K}^{4} $ \`a partir d'un fibr\'e de rang $2$ sur la vari\'et\'e 
$\P^{3}:=\P_{\K}^{3} $ o\`u $\K$ est un corps quelconque, tout en utilisant \ref{4.1} et \ref{3.1} .

\subsubsection{\bf Remarque}\label{4.2.3} 
Soient $z \in H^{0}(\P^{4},\ko_{\P^{4}}(1))\setminus \lbrace 0 \rbrace $ et $\P^{3}   \subset  \P^4 $ l'hyperplan de $\P^4$ d\'efini par l'\'equation $z=0$. Soient $\epsilon\geq 1$ un entier et $Y_{(\epsilon)}\subset \P^4$ le voisinage infinit\'esimal de l'ordre $\epsilon$ de $\P^{3} $ dans $\P^4$, d\'efini par l'\'equation $z^{\epsilon}=0$. On a les morphismes d'inclusion suivants

 \xmat{  \P^{3} :=Y_{(1)}  \ar@{^{(}->}[rr]^{b_{(2)}} &  & Y_{(2)} \ar@{^{(}->}[rr]^{b_{(3)}} &  &  \ldots \ar@{^{(}->}[rr]^{b_{(\epsilon)}} &  &   Y_{(\epsilon)}\ar@{^{(}->}[rr]^{d} &  &  \P^{4} . \\
 }

Soit $e_{0}\in \P^{4}\setminus Y_{(\epsilon)}$, on d\'efinit la projection d'un point (voir \cite{ha} page $22$) sur un sous-espace projectif $\P^{3}$ par
$$\pi:\P^{4}\setminus \{ e_{0}\}\longrightarrow \P^{3}$$
$$ e\longmapsto \pi(e)$$

o\`u $\pi(e)$ est l'intersection de $\P^{3}$ avec la droite unique passant par les points $e_{0}, e$. On consid\`ere la restriction de la projection $\pi$ sur le voisinage infinit\'esimal $Y_{(\epsilon)}$, on obtient la projection d'un point sur un sous-espace projectif $\P^{3}$
$$J_{(\epsilon)}:Y_{(\epsilon)}\longrightarrow \P^{3}.$$

On consid\`ere $\kb_{(\epsilon)}=J_{(\epsilon)}^{*}\kb$ et $\ka_{(\epsilon)}=J_{(\epsilon)}^{*}\ka$ o\`u $\kb$ et $\ka$ sont les fibr\'es vectoriels dans la proposition \ref{4.1.5}. Alors on a la suite exacte suivante
$$0\longrightarrow \ka_{(\epsilon)}\longrightarrow \F_{(\epsilon)}  \longrightarrow   \kb_{(\epsilon)}   \longrightarrow 0,$$

o\`u $\F_{(\epsilon)}=  \ko_{Y_{(\epsilon)}}\oplus \ko_{Y_{(\epsilon)}}(t_{1}-v_{2})\oplus \ko_{Y_{(\epsilon)}}(t_{1}-w_{2})\oplus \ko_{Y_{(\epsilon)}}(w_{2}-v_{2}) $.

\bigskip
\subsubsection{\bf Proposition}\label{4.2.4} 
{\em Soient $\ka_{(\epsilon)}$ et $\kb_{(\epsilon)}$ comme dans la remarque \ref{4.2.3} , et $G_{2}$ un fibr\'e vectoriel  sur $\P^{4}$ d\'efini par le diagramme commutatif suivant

\xmat{ &0\ar[d]&0\ar[d]&&\\
&\hat F(-\epsilon)\ar@{=}[r] \ar[d]& \hat F(-\epsilon)\ar[d]&&\\
0 \ar[r]& G_{2} \ar[r]\ar[d]&\hat F \ar[r]\ar[d]& \kb_{(\epsilon)} \ar[r]\ar@{=}[d]&0\\
0 \ar[r]& \ka_{(\epsilon)}\ar[r]\ar[d]&\F_{(\epsilon)} \ar[r]\ar[d]& \kb_{(\epsilon)} \ar[r]&0\\
&0&0&&}

o\`u $\hat \F := \ko_{\P^4}\oplus \ko_{\P^4}(2w_{2}-w_{1}-v_{2})\oplus \ko_{\P^4}(w_{2}-w_{1})\oplus \ko_{\P^4}(w_{2}-v_{2}) $. Dans la proposition \ref{4.1.5}, si on a que $c_{1}(\ka)= w_{2}-w_{1}-2v_{2}=-\epsilon\leq -1$ et $\epsilon$ v\'erifie une des conditions suivantes

1- $ 2v_{2}-w_{2}\leq \epsilon\leq 2v_{2}+w_{2} $ si $v_{2}\geq w_{2}>0$

2- $ v_{2} \leq \epsilon\leq 3v_{2}  $ si $0< v_{2}\leq w_{2}$

alors il existe un fibr\'e vectoriel $E_{2}$ de rang $3$ sur $\P^4$ tel que
$$0 \longrightarrow \ko_{\P^4}(-\epsilon) \longrightarrow G_{2} \longrightarrow E_{2} \longrightarrow 0,$$

qui est diff\'erent du fibr\'e de Tango pond\'er\'e provenant d'une image inverse
 g\'en\'eralis\'ee sur $\P^{4}$.}

\begin{proof} Soient $\ka, \hspace{0.2cm} \kb  $ les fibr\'es  dans la proposition \ref{4.1.5}, on a le diagramme suivant sur $\P^{3}$

\xmat{  \F^{*}( w_{2}-w_{1}-2v_{2})  \ar[rr]^{\varphi_{0}}\ar[rd] &  & \F 
\ar[rr]^{\psi_{0}}\ar[rd]  & & \F^{*}( 3 w_{2}- w_{1})\\
&   \ka \ar[ru]  \ar[rd] & &  \kb    \ar[rd]\ar[ru]   \\
 0\ar[ru] &&   0\ar[ru]    && 0 } 

o\`u $\F  = \ko_{\P^3}\oplus \ko_{\P^3}(2w_{2}-w_{1}-v_{2})\oplus \ko_{\P^3}(w_{2}-w_{1})\oplus \ko_{\P^3}(w_{2}-v_{2}) $. En prenant l'image inverse du diagramme pr\'ec\'edent par la projection d'un point sur un sous-espace projectif $\P^{3}$, \ref{4.2.3} ,
$$J_{(\epsilon)}:Y_{(\epsilon)}\longrightarrow \P^{3} ,$$

on obtient le diagramme suivant sur $Y_{(\epsilon)}$

\xmat{  \F_{(\epsilon)}^{*}( w_{2}-w_{1}-2v_{2}) \ar[rr]^{\varphi_{0}}\ar[rd] &  & \F_{(\epsilon)} 
\ar[rr]^{\psi_{0}}\ar[rd]& & \F_{(\epsilon)}^{*}( 3 w_{2}- w_{1})    \\
& \ka_{(\epsilon)} \ar[ru]\ar[rd] & &  \kb_{(\epsilon)} \ar[rd]\ar[ru]  \\
0\ar[ru] &&   0\ar[ru]    && 0 } 

o\`u $\F_{(\epsilon)} =\ko_{Y_{(\epsilon)}}\oplus \ko_{Y_{(\epsilon)}}(2w_{2}-w_{1}-v_{2})\oplus \ko_{Y_{(\epsilon)}}(w_{2}-w_{1})\oplus \ko_{Y_{(\epsilon)}}(w_{2}-v_{2})$. 
En choisissant les m\^emes notations de \ref{3.1} avec $\Phi=\varphi_{0}$ et $\Psi=\psi_{0}$ sur $\P^{4}$, tous les fibr\'es $\F_{(\epsilon)}^{*}( w_{2}-w_{1}-2v_{2}) $, $  \F_{(\epsilon)} $ et $ \F_{(\epsilon)}^{*}( 3 w_{2}- w_{1})$ peuvent se relever en fibr\'es $\hat F^{*}( w_{2}-w_{1}-2v_{2}) $, $  \hat F $ et $ \hat F^{*}( 3 w_{2}- w_{1})$ respectivement sur $\P^{4}$. On obtient donc le complexe suivant sur $\P^{4}$
$$\hat F^{*}( w_{2}-w_{1}-2v_{2}) \stackrel{ \Phi }{\longrightarrow} \hat F 
\stackrel{ \Psi }{\longrightarrow} \hat F^{*}( 3 w_{2}- w_{1}), $$

o\`u  $\Psi .\Phi  =\psi_{0}.\phi_{0}=0$. Donc on peut d\'efinir le morphisme suivant sur $\P^{4}$
$$\Delta:\hat F(-\epsilon)\oplus \hat F^{*}(  w_{2}-w_{1}-2v_{2})\longrightarrow \hat 
F\oplus \hat F^{*}( 3 w_{2}- w_{1} -\epsilon), $$

$$
\Delta = \left(
\begin{array}{c|c}
z^{\epsilon}I    &    \Phi \\
\hline
\Psi  &    0 \\
\end{array}
\right)
.$$ 

$$
\Delta = \left(
\begin{array}{cccc|cccc}
z^{\epsilon}&&&\bigzero& 0 & S_{3,4}& -S_{2,4} &S_{3,2} \\
&z^{\epsilon}&&&-S_{3,4} & 0 &  S_{1,4} &-S_{1,3} \\
\bigzero&&z^{\epsilon}&&S_{2,4} & -S_{1,4} &0  &S_{1,2} \\
&&&z^{\epsilon}&-S_{3,2} & S_{1,3} & -S_{1,2} & 0 \\
\hline
0 & S_{1,2}& S_{1,3} &S_{1,4}  &&&&\\
S_{2,1} & 0 &  S_{2,3} &S_{2,4}&&\bigzero&&\\
 S_{3,1} & S_{3,2} &0  &S_{3,4}&&&&\\
S_{4,1} & S_{4,2} & S_{4,3} & 0&&& \\
\end{array}
\right)
$$

On a
 \xmat{  \hat F(-\epsilon)\oplus \hat F^{*}( w_{2}-w_{1}-2v_{2})   
 \ar[rr]^{\Delta}\ar[rd] &  &  \hat F\oplus \hat F^{*}(  3 w_{2}- w_{1}-\epsilon)   \\
&   G_{2} \ar[ru]\ar[rd]\\
0\ar[ru] &&0 }

D'apr\`es la proposition \ref{3.1.2} et la remarque \ref{3.1.3}. Soit 
$$ (h_{1},h_{2}, h_{3} ,h_{4},h_{5},h_{6}, h_{7} ,h_{8})\in H^{0}(\ko_{\P^4}(-\epsilon))\oplus H^{0}(\ko_{\P^4}(2w_{2}-w_{1}-v_{2}-\epsilon))\oplus H^{0}(\ko_{\P^4}(w_{2}-w_{1}-\epsilon))$$
$$\oplus H^{0}(\ko_{\P^4}(w_{2}-v_{2}-\epsilon))\oplus H^{0}(\ko_{\P^4}(w_{2}-w_{1}-2v_{2}))\oplus H^{0}(\ko_{\P^4}(-w_{2}-v_{2}))$$
$$\oplus H^{0}(\ko_{\P^4}(-2v_{2}))\oplus H^{0}(\ko_{\P^4}(-w_{1}-v_{2})) . $$

On prend la somme la $1^{\text{\`ere}}$-colonne et la $5^{\text{\`eme}}$-colonne de la matrice $\Delta$ avec la condition d'homog\'en\'eit\'e suivante, sur la $1^{\text{\`ere}}$-coordonn\'ee et la $5^{\text{\`eme}}$-coordonn\'ee de

$(h_{1},h_{2}, h_{3} ,h_{4},h_{5},h_{6}, h_{7} ,h_{8})$, 
$$deg(h_{1})=deg(h_{5})$$ 
$$-\epsilon=w_{2}-w_{1}-2v_{2}.$$ 

Une telle somme existe gr\^ace \`a la condition $c_{1}(\ka)= w_{2}-w_{1}-2v_{2}=-\epsilon$. Donc on obtient une colonne $\tau$ 
$$\tau=^{T}[z^{\epsilon}, S_{1,2},S_{1,3} ,S_{1,4};0,S_{3,4},-S_{2,4} ,S_{3,2}]$$

$$\tau:\ko_{\P^{4}}(w_{2}-w_{1}-2v_{2})\longrightarrow  \hat F\oplus\hat F^{*}(4 w_{2}- 2w_{1}-2v_{2}) . $$

Supposons que toutes les formes de la matrice $\tau$ soient nulles sur $ P^{4}$. On obtient donc que la matrice $\Phi$ est nulle, ce qui est une contradiction au fait que $rg (\Phi)=2$. Donc pour tout $x\in \P^{4}$ il existe une forme des formes $z^{\epsilon}, S_{1,2},S_{1,3} ,S_{1,4}, S_{3,4},S_{2,4} ,S_{3,2}$ qui ne s'annule pas en $x$. Ce qui donne que le morphisme $\tau$ est injectif non nul sur $\P^{4}$. Mais, le morphisme suivant est repr\'esent\'e par le morphisme $\tau $
$$\ko_{\P^{4}}(w_{2}-w_{1}-2v_{2})\longrightarrow G_{2};$$
 
alors on a la suite exacte suivante
$$0 \longrightarrow \ko_{\P^4}(w_{2}-w_{1}-2v_{2}) \longrightarrow G_{2} \longrightarrow E_{2} \longrightarrow 0,$$

o\`u  $E_{2}$ est un fibr\'e vectoriel de rang $3$ sur $\P^4$. On utilise la proposition  \ref{3.1.8} pour calculer les classes de Chern du faisceau $d_{*}(\kb_{(-w_{2}+w_{1}+2v_{2})} )$. En utilisant la suite exacte pr\'ec\'edente et la suite exacte suivante
$$0\longrightarrow G_{2} \longrightarrow \hat F \longrightarrow \kb_{(-w_{2}+w_{1}+2v_{2})} 
\longrightarrow 0$$

on obtient les classes de Chern du fibr\'e $E_{2}$   
$$c_{1}(E_{2})=2(w_{2}+v_{2})-3(-w_{2}+w_{1}+2v_{2}),$$

$$c_{2}(E_{2})= 2(-w_{2}+w_{1}+2v_{2})^{2}-(-w_{2}+w_{1}+2v_{2})(3w_{2}+v_{2})+w_{2}^{2}-v_{2}^{2}+4w_{2}v_{2},$$
 
$$c_{3}(E_{2})=2v_{2}(-(-w_{2}+w_{1}+2v_{2})^{2}+ 2v_{2}(-w_{2}+w_{1}+2v_{2})+w_{2}^{2}-v_{2}^{2}).$$
  
\end{proof}

\subsubsection{\bf Lemme}\label{4.2.4.1}\cite{ok-sc-sp}
{\em Soit $E$ un fibr\'e vectoriel de rang $r$ sur $\P^{n}$ engendr\'e par ses sections globales. Si la classe de Chern $c_{r}(E)=0$, alors il existe une section non-nulle de $E$. Autrement dit, $E$ contient un sous-fibr\'e trivial de rang $1$}.

\subsubsection{\bf Th\'eor\`eme}\label{4.2.5}
{\em Soient $E_{2}$, $G_{2}$, $\ka_{(\epsilon)}$ et $\kb_{(\epsilon)}$ des fibr\'es vectoriels sur $\P^{4}$ comme dans la proposition \ref{4.2.4} avec $v_{2}=w_{2}=w_{1}=n>0,\epsilon=2n $ des entiers. Alors il existe un fibr\'e vectoriel ind\'ecomposable $L_{2}$ de rang $2$ sur $\P^4$ tel que 
$$0 \longrightarrow \ko_{\P^4} \longrightarrow E_{2} \longrightarrow L_{2} \longrightarrow 0,$$

son polyn\^ome de Chern est $c_{h}(L_{2})=1-2nh+4n^{2}h^{2}$ o\`u $h=c_{1}(\ko_{\P^4}(1))$.}

\begin{proof} D'apr\`es la proposition \ref{4.2.4}, si on consid\`ere $v_{2}=w_{2}=w_{1}=n>0,\epsilon=2n $, on obtient
le diagramme commutatif suivant
\xmat{ &0\ar[d]&0\ar[d]&&\\
& \K^{4}\otimes\ko_{\P^4}(-2n)\ar@{=}[r] \ar[d]& \K^{4}\otimes\ko_{\P^4}(-2n)\ar[d]&&\\
0 \ar[r]& G_{2} \ar[r]\ar[d]& \K^{4}\otimes\ko_{\P^4} \ar[r]\ar[d]& \kb_{(2n)} \ar[r]\ar@{=}[d]&0\\
0 \ar[r]& \ka_{(2n)}\ar[r]\ar[d]&\K^{4}\otimes\ko_{Y_{(2n)}} \ar[r]\ar[d]& \kb_{(2n)} \ar[r]&0\\
&0&0&&}

avec $E_{2}$ un fibr\'e vectoriel de rang $3$ sur $\P^4$ tel que
$$0 \longrightarrow \ko_{\P^4}(-2n) \longrightarrow G_{2} \longrightarrow E_{2} \longrightarrow 0,$$

son polyn\^ome de Chern est $c_{h}(E_{2})=1-2nh+4n^{2}h^{2}$ o\`u $h=c_{1}(\ko_{\P^4}(1))$. Donc on a 
$$h^{0}(E_{2})=h^{0}(G_{2})=h^{0}(\ka_{(2n)})=1.$$
 
On va d\'emontrer que le fibr\'e $E_{2}$ est engendr\'e par ses sections globales.\\ Pour tout $y=\K.v \in \P^4=\P(V)$ o\`u $v\in V$, on a le diagramme commutatif suivant

\xmat{ 
0 \ar[r]& H^{0}(\ka_{(2n)}) \ar[r]\ar[d]^{ev_{\ka_{(2n)}}}& (\K^{4})^{*} \ar[r]^{q}\ar[d]^{ev_{\K^{4}\otimes\ko_{Y_{(2n)}}}}& H^{0}(\kb_{(2n)}) \ar[d]^{ev_{\kb_{(2n)}}}\ar[r]& \\
0 \ar[r]& \ka_{(2n)}{}_{y}\ar[r]&\K^{4} \ar[r]^{q}\ar[d]& \kb_{(2n)}{}_{y}\ar[r]& 0\\
&&0& &}

o\`u $q=[q_{1},q_{2},q_{3},q_{4}]$ et $q_{i}\in H^{0}(\K^{4}\otimes\ko_{Y_{(2n)}})=(\K^{4})^{*},\hspace{0,2 cm} i=1,\ldots,4, $ sont des sections de $\kb_{(2n)} $. Donc les g\'en\'erateurs de $H^{0}(\ka_{2n})$ sont $e_{2}=(-q_{2},q_{1},0,0),e_{3}=(-q_{3},0,q_{1},0),e_{4}=(-q_{4},0,0,q_{1})$. Alors pour tout $X=(X_{1},X_{2},X_{3},X_{4})\in \ka_{(2n)}{}_{y}$, on a $qX=0$. On obtient que
$$X q_{1}= X_{2}e_{2}+X_{3}e_{3}+X_{4}e_{4} \in  H^{0}(\ka_{(2n)}).$$

Donc $ X\in H^{0}(\ka_{(2n)}) $, et le fibr\'e $\ka_{(2n)}$ est engendr\'e par ses sections globales. On a aussi les carr\'es commutatifs suivants
\xmat{    G_{2} \ar[r]    & \ka_{(2n)} \ar[r] & 0\\
H^{0}(G_{2})\otimes \ko_{\P^4} \ar[r]^{\backsim}\ar[u]^{ev_{G_{2}}} &  H^{0}(\ka_{2n})\otimes \ko_{\P^4}\ar[u]_{ev_{\ka_{2n}}} \\  }

et
\xmat{    G_{2} \ar[r]    & E_{2} \ar[r]& 0\\
H^{0}(G_{2})\otimes \ko_{\P^4} \ar[r]^{\backsim}\ar[u]^{ev_{G_{2}}} &  H^{0}(E_{2})\otimes \ko_{\P^4}\ar[u]_{ev_{E_{2}}} \\  }

qui nous donnent que le fibr\'e $E_{2}$ est engendr\'e par ses sections globales. D'apr\`es le lemme \ref{4.2.4.1} et comme on a $c_{3}(E_{2})=0$, on obtient alors que le fibr\'e $E_{2}$ a un sous-fibr\'e en droite. Autrement dit, on a un fibr\'e vectoriel $L_{2}$ de rang $2$ sur $\P^4$ d\'efini par la suite exacte
$$0\longrightarrow  \ko_{\P^4} \longrightarrow E_{2} \stackrel{ }{\longrightarrow} L_{2}\longrightarrow 0,$$

son polyn\^ome de Chern est $c_{h}(L_{2})=c_{h}(E_{2})=1-2nh+4n^{2}h^{2}$. Comme ce polyn\^ome de Chern est irr\'eductible dans l'anneau $H^{*}(Pn, \Z) = \Z[h] = \Z[t]/(t^{5})$, alors le fibr\'e $L_{2}$ est ind\'ecomposable.

\end{proof}

\end{document}